\magnification=1200
\input amstex
\documentstyle{amsppt}
\hoffset=-0.5pc
\nologo

\vsize=57.2truepc
\hsize=38.5truepc
\spaceskip=.5em plus.25em minus.20em

\define\Sigm{\roman S}

\define\fra{\frak}
\define\rbar{\overline{\roman B}}
\define\rcob{\overline{\Omega}}

\define\Bobb{\Bbb}
\define\Mfl{M^{\flat}}
\define\Nfl{N}
\define\atibottw{1}
\define\borelboo{2}
\define\bottone{3}
\define\bottsega{4}
\define\botshust{5}
\define\CC{A}
\define\BB{B}

\define\hhh{h}
\define\CaA{U}
\define\ebrown{6}
\define\brownez{7}
\define\cartanon{8}
\define\cartantw{9}
\define\chenone{10}
\define\chentwo{11}
\define\chenfou{12}
\define\doldone{13}
\define\doldpupp{14}
\define\drachman{15}
\define\eilmactw{16}
\define\eilmotwo{17}

\define\eilmothr{18}
\define\eilenmoo{19}
\define\frankone{20}
\define\franzone{21}

\define\franztwo{22}
\define\franpupp{23}
\define\gorkomac{24}
\define\gugenhtw{25}
\define\gugenhon{26}
\define\gulstatw{27}
\define\gugenmay{28}
\define\gugenmun{29}
\define\habili{30}
\define\perturba{31}
\define\cohomolo{32}
\define\modpcoho{33}
\define\intecoho{34}
\define\kan{35}
\define\berikas{36}
\define\minimult{37}
\define\pertlie{38}
\define\pertlitw{39}
\define\origins{40}
\define\duatwo{41}
\define\equivg{42}
\define\tornike{43}
\define\huebkade{44}
\define\huebstas{45}
\define\husmosta{46}
\define\kirwaboo{47}
\define\maclaboo{48}
\define\mooretwo{49}
\define\munkholm{50}
\define\munkhotw{51}
\define\munkhthr{52}
\define\rothstee{53}
\define\sourithr{54}
\define\stashalp{55}
\define\steeepst{56}

\define\Nsddata#1#2#3#4#5{
 ( #4
\text{
\vbox
to 1.15 pc
{
  \hbox {$@>{\,\,\, #3\,\,\,}>>$}
  \vskip-1.2pc
  \hbox {$@<<{\,\,\, #2\,\,\,}<$}
                }
      }
 #1,#5 )
}

\topmatter
\title  Homological perturbations, equivariant cohomology,\\
and Koszul duality
\endtitle
\author Johannes Huebschmann
\endauthor
\affil
Universit\'e des Sciences et Technologies de Lille
\\
UFR de Math\'ematiques
\\
CNRS-UMR 8524
\\
F-59 655 VILLENEUVE D'ASCQ C\'edex, France
\\
Johannes.Huebschmann\@math.univ-lille1.fr
\endaffil
\date{July 31, 2009}
\enddate
\abstract{Our main objective is to demonstrate how homological 
perturbation theory (HPT) results over the last
40 years immediately or with little extra work
give some of the Koszul duality results that have appeared in the
last decade.

Higher homotopies typically arise when a huge object, e.~g. a chain complex
defining various 
invariants of a certain geometric situation, is cut to a small model,
and the higher homotopies can then be dealt with concisely
in the language of sh-structures (strong homotopy structures). This amounts
to precise ways of handling the requisite additional structure encapsulating
the various coherence conditions.
Given e.~g. two augmented differential graded algebras $A_1$ and $A_2$,
an sh-map from $A_1$ to $A_2$ is a twisting cochain
from the reduced bar construction $\rbar A_1$ of $A_1$ to $A_2$
and, in this manner, the class of morphisms of augmented differential 
graded algebras is extended to that of sh-morphisms.
In the present paper, 
we explore small models for equivariant (co)homology via differential
homological algebra techniques including homological perturbation theory
which, in turn, is a standard tool to handle  sh-structures.

Koszul duality, for a finite type exterior algebra $\Lambda$  on odd
positive
degree generators, then comes down to a duality between the category of
sh-$\Lambda$-modules and that of sh-$\rbar \Lambda$-comodules. 
This kind of duality
relies on the extended functoriality of the differential graded
Tor-, Ext-, Cotor-, and Coext functors, extended to the appropriate
sh-categories.
We construct the small models as certain twisted tensor products and twisted
Hom-objects. These are chain and cochain models for the chains and cochains
on geometric bundles and are compatible with suitable additional structure.}
\endabstract
\endkeywords
\dedicatory
Dedicated to the memory of V.K.A.M. Gugenheim
\enddedicatory
\toc
\widestnumber\head{5$_*$.}
\widestnumber\subhead{3.2.1.}
\specialhead{ } Introduction \endspecialhead
\head 1. Preliminaries\endhead
\subhead 1.1. Twisting cochains and twisted tensor products\endsubhead
\subhead 1.2. Twisted Hom-objects \endsubhead
\subhead 1.3. The trivial simplicial object and the total simplicial object\endsubhead
\subhead 1.4. The geometric resolution and nerve construction\endsubhead
\subhead 1.5. The notion of model \endsubhead
\subhead 1.6. Sh-modules and sh-comodules\endsubhead
\head 2. Homological perturbations\endhead
\subhead 2.1. Definition\endsubhead
\subhead 2.2. Proposition\endsubhead
\subhead 2.3. The perturbation lemma\endsubhead
\head 3. Duality\endhead

\head 4. Small models for ordinary equivariant (co)homology\endhead
\head 5$_*$. Koszul duality for homology\endhead
\head 5$^*$. Koszul duality for cohomology\endhead
\head 6. Grand unification\endhead
\head 7. Some illustrations\endhead
\head 8. Split complexes and generalized momentum mapping\endhead
\head{ } References\endhead
\endtoc

\endtopmatter
\document

\leftheadtext{Johannes Huebschmann}
\rightheadtext{HPT, equivariant cohomology, Koszul duality}

\head {Introduction} \endhead

The purpose of this paper is to explore suitable small models for equivariant
(co)homology via differential homological algebra techniques including,
in particular, {\it homological perturbation theory\/} (HPT). 
Our main objective is to demonstrate how HPT results obtained over the last
40 years immediately or with little extra work
give some of the Koszul duality results that have appeared in the
last decade.
Indeed, we will show that these Koszul duality functors
yield small models for various differential graded  Tor-, Ext-, Cotor-,
and Coext functors  which, in turn, entail a {\it conceptual explanation
of Koszul duality in terms of the extended functoriality of these functors\/}.
Koszul duality, for a finite type exterior algebra $\Lambda$  on 
odd positive
degree generators, then comes down to a duality between the category of
sh-$\Lambda$-modules and that of sh-$\rbar \Lambda$-comodules. 
Our approach also applies to situations which, to our knowledge, have not been
addressed in the literature, e.~g. Koszul duality for infinite dimensional
groups. 
The methods we reproduce and newly develop here
provide constructive means to explore and handle the 
requisite additional operations of $\roman H_*G$  on $C_*(X)$
to recover the original action, and categories of sh-modules etc. serve
as replacements for various derived categories. We hope to convince
the reader that HPT results obtained in the last 40 years 
are still relevant today and should be better known among mathematicians
working in equivariant cohomology.

In general, when standard differential homological algebra techniques yield
a huge object calculating some (co)homological invariant and when this
object is replaced by an equivalent small model, higher homotopies arise.
These higher homotopies reflect a loss of a certain high amount of symmetry
which the huge object usually has; thus the huge object is not necessarily
\lq\lq bad\rq\rq.
HPT is a standard tool to handle such higher homotopies; we use it here,
in particular, to construct twisting cochains and contractions.
The notion of twisting cochain in differential homological algebra,
introduced in \cite\ebrown,
is intimately related to that of connection in differential geometry,
cf. \cite\cartanon, 
\cite\cartantw, \cite\chenfou, \cite\berikas,
as well as to the {\it Maurer-Cartan\/} or {\it master\/} equation, cf.
\cite\huebstas.

We now give a brief overview of the paper.
Section 1 is preliminary in character; is contains a eview of
various differential homological algebra techniques and includes in particular
a discussion of twisting cochains.
In Section 2, we explain the
requisite HPT-techniques. In Section 3 we 
discuss a general notion of duality
that arises from an acyclic twisting cochain from a coaugmented differential
graded coalgebra $C$ to an augmented differential graded algebra $A$;
the general duality is 
one between the categories of $A$-modules and $C$-comodules.
This kind of duality prepares for the conceptual development of Koszul
duality later in the paper.
The small models for (singular) equivariant
$G$-(co)homology for a topological group $G$ having as homology a strictly
exterior algebra in odd degree generators 
of finite type in the sense that each homogeneous constituent is
(necessarily free and)  finite-dimensional
will  be given in
Section 4. These models rely on the familiar fact, cf. Theorem 4.1 below
that, given a topological group $G$, the $G$-equivariant homology and
cohomology of a $G$-space are given by a {\it differential\/} Tor
and a {\it differential\/} Ext, respectively, with
reference to the chain algebra of $G$. This is an Eilenberg-Moore type result.
Our approach includes a HPT proof thereof
which involves the twisted Eilenberg-Zilber theorem, 
see Lemma 4.2 below.
A variant thereof, given as Theorem 4.1${}^*$, describes the 
$G$-equivariant
cohomology as a suitable Cotor, with reference to an appropriate
coalgebra which is a replacement for the cochains on $G$
(which do not inherit a coalgebra structure since in general the dual of a
tensor product is not the tensor product of the duals).

In Section 5 we exploit the techniques
developed or reproduced in earlier sections 
to deduce the main Koszul duality results in a conceptual manner.
Our constructions of small models extend some of the results in
\cite{\franzone,\,\franztwo,\,\gorkomac} related with {\it Koszul
duality\/} by placing the latter in the sh-context
in the sense isolated in 
the seminal paper
\cite\stashalp\ of Stasheff and  Halperin;
the theory of sh-modules 
was then exploited in \cite{\gugenmun,\, \munkhotw,\, \munkhthr}
and pushed further in our paper \cite\habili.
The \lq up to homotopy\rq\ interpretation of Koszul duality
has been known for a while in the context of operads as well
and is also behind e.~g. \cite{\franzone}.
Yet we believe that our approach in terms of HPT
clarifies the meaning and significance of Koszul duality, cf.
Section 5 below. In particular, when a group $G$ acts on a space
$X$, even when the induced action of $\roman H_*G$ on $\roman
H_*X$ lifts to an action on $C_*(X)$, in general only an sh-action
of $\roman H_*G$  on $C_*(X)$ with non-trivial higher terms will
recover the geometry of the original action. 
For example, the ordinary cohomology of a homogeneous
space $G/K$ of a compact Lie group $G$ by a closed subgroup $K$
amounts to the $K$-equivariant cohomology of $G$, and the
corresponding sh-action of $\roman H_*K$ on $C_*(G)$ {\it
cannot\/} come down to an ordinary $(\roman H_*K)$-action on
$C_*(G)$ since the homogeneous space is compact and finite
dimensional. For illustration, let $G=\roman{SU(2)}$ and let $K$
be a maximal torus which is just a circle group $S^1$. For degree
reasons, the induced action of $\roman H_*(S^1)$ on $\roman
H_*(G)$ is trivial and hence lifts to the trivial action of
$\roman H_*(S^1)$ on $C_*(G)$. However this lift cannot recover
the geometry  of the original $S^1$-action since if it did, the
homology of the homogeneous space $G/S^1$ would have a free
generator in each degree which is impossible since this
homogeneous space is just the 2-sphere. 
This example occurs in \cite\gorkomac\ (1.5).
More details and a class of examples including
the one under discussion
can be found in Example 7.2 below.

A similar illustration is given in Example 7.1 below.
This example serves, in particular, as an illustration for the 
notion of twisting cochain.
According to
a result of Frankel's \cite\frankone\ and
Kirwan's \cite\kirwaboo, a smooth compact symplectic
manifold, endowed with a hamiltonian action of a compact Lie
group, is equivariantly formal over the reals.
Both examples show that the compactness hypothesis is crucial.

What are referred to as {\it cohomology
operations\/} in \cite\gorkomac\ is really a {\it system of higher
homotopies encapsulating the requisite sh-action of\/} $\roman
H_*G$  on $C_*(X)$ or $C^*(X)$. 
Details will be given in Section 5. 
Suffice it to recall here that an sh-action of an augmented differential graded
algebra $A$ on a chain complex $V$ amounts to a twisting cochain
from the reduced bar construction $\rbar A$ of $A$ to the differential graded
algebra $\roman{End}(V)$ of endomorphisms of $V$.
Our HPT-techniques
provide the necessary algebraic machinery so that we
can handle these higher homotopies and in particular
carry out the requisite constructions of e.~g. twisting
cochains and homotopies of twisting cochains. Koszul
duality, for a finite type exterior algebra $\Lambda$  on odd
positive degree generators, comes down to a duality between the category of
sh-$\Lambda$-modules and that of sh-$\rbar \Lambda$-comodules and,
in this context, $\rbar \Lambda$ and the corresponding cofree
graded symmetric coalgebra $\Sigm'$ are equivalent. For our
purposes, the categories of sh-modules and sh-comodules serve as
{\it replacements for various derived categories\/} exploited in
\cite\gorkomac\ and elsewhere. That sh-structures may be used as
replacement for certain derived categories has been observed
already in \cite\munkhthr.

The duality, then, amounts essentially to the fact that, when 
$I\to M$ is a relatively injective 
resolution of an $\Sigm'$-comodule $M$, cf. e.~g. \cite\eilenmoo\ 
for basic notions of relative homological algebra,
this injection is a chain equivalence and that, when $P \to N$ is
a relatively projective resolution of a $\Lambda$-module $N$, 
this projection
is a chain equivalence.
We explain these chain equivalences under somewhat more general 
circumstances in Section 3.
Concerning Koszul duality,
given 
the $\Lambda$-module
$N$ and 
the $\Sigm'$-comodule
$M$, when $h$ and
$t$ refer to the two Koszul duality functors, $h(t(N))$ is a
generalized projective resolution of $N$ and $t(h(M))$ a generalized
injective resolution of $M$; see (5.3) below for details. Koszul
duality reflects the old observation that, for any $k \geq 0$, the
non-abelian left derived functor of the $k$-th exterior power functor in
the sense of \cite\doldpupp\ is the $k$-th symmetric copower on
the suspension (the invariants in the $k$-th tensor power on the
suspension with respect to the symmetric group on $k$ letters).

In Section 6 we offer some unification: 
For a group $G$ of strictly exterior type, exploiting {\it
isomorphisms\/} between $C_*G$ and $\roman H_*G$, between $\rbar
C_*G$ and $\roman H_*(BG)$, and between $\rbar^* C_*G$ and $\roman
H^*(BG)$, in suitably defined categories of sh-algebras and
sh-coalgebras as appropriate, we show that the Koszul duality
functors are equivalent to certain duality functors on the
categories of $(C_*G)$-modules, $(\rbar^* C_*G)$-modules, and $(\rbar
C_*G)$-comodules as appropriate. This equivalence amounts more or
less to the extended functoriality of the differential graded
Ext-, Tor-, and Cotor-functors in the sense of \cite\gugenmun. 
We explain this extended functoriality briefly in Section 6 below.
The
idea behind these isomorphisms goes back to the quoted paper
\cite\stashalp\ of Stasheff and  Halperin.

In Section 7 we illustrate our results with a number of examples
from equivariant cohomology.
As a final application, in Section 8,
we explore a notion of split complex, similar to that
introduced in \cite\gorkomac\ but somewhat more general and
adapted to our situation. Our approach includes an interpretation
of the splitting homotopy as a generalization of the familiar
momentum mapping in symplectic geometry.
In particular, using HPT-techniques, we isolate the difference between  
equivariant formality and the stronger property that the equivariant 
cohomology is an induced module over the cohomology of the 
classifying space.

In the present paper, spectral sequences, one of the basic tools in the
literature on Koszul duality and related topics, do not occur explicitly
save that we establish the acyclicity of a number of chain complexes
by means of an elementary spectral sequence argument. 
A long time ago, S. Mac Lane pointed out to me that, in various circumstances,
spectral sequences are somewhat too weak a tool, on the conceptual as well as
computational level. This remark prompted me to develop small models of the
kind given below via HPT techniques during the 1980's; 
the HPT-techniques include
refinements of reasoning usually carried out in the literature via spectral
sequences. In
\cite{\perturba}--\cite \intecoho, I have constructed and exploited suitable
small models encapsulating the
appropriate sh-structures in the (co)homology of a discrete group,
and by means of these small models, I did explicit numerical
calculations of group cohomology groups that still today cannot be done by
other methods. Equivariant cohomology may be viewed as being part of group
cohomology and, in this spirit, the present paper pushes further
some of the ideas developed in 
\cite{\habili}--\cite \kan\ and
in \cite \huebkade.

In a follow-up paper \cite\duatwo, we have worked out a related approach to
equivariant de Rham theory in the framework of
suitable relative derived functors, and in \cite\equivg\ we have extended
this approach to equivariant de Rham theory relative to a Lie groupoid.

I am much indebted to Jim Stasheff
and to the referee for a number of comments which
helped improve the exposition and made the paper, perhaps, 
available, to a larger audience than the original draft.

An earlier version of this 
paper was
posted to the arxiv under math.AT/0401160.
Publication of the material
has been delayed for personal (non-mathematical)
reasons.

\medskip\noindent {\bf 1. Preliminaries}
\smallskip\noindent
Let $R$ be a commutative ring with $1$, taken henceforth as ground ring
and not mentioned any more.
Graded objects will be $\Bobb Z$-graded; in the applications
we will mainly consider
only non-negatively or non-positively graded objects,
and this will then be indicated.

We will treat
chain complexes and cochain complexes on equal footing:
We will consider 
a cochain complex $(C^*,d^*)$ as a chain complex
$(C_*,d_*)$ by letting $C_j = C^{-j}$ and
$d_*=d^*\colon C_j = C^{-j} \to C^{-j+1}=C_{j-1}$,
for $j \in \Bobb Z$.
An ordinary cochain complex, concentrated in non-negative
degrees as a cochain complex, is then a chain complex which
is {\it concentrated in non-positive degrees\/}.
Differential homological algebra terminology and notation
will essentially be the same as
that in \cite\gugenmun\ and \cite\husmosta.
We will write the reduced bar and cobar functors  
as $\rbar$ and $\rcob$, respectively, 
rather than as $B$ and $\Omega$;
see e.~g. \cite\munkholm\ for explicit descriptions.
In Section 6 we will reproduce an explicit description of the reduced bar
construction.
Definitions of the
differential Cotor-functor may be found in \cite \eilenmoo\
(p.~206) and in \cite\husmosta\ (Chap. 1), and definitions of the
differential Tor and  Ext functors may be found in \cite\gugenmay\
(p.~3 and p.~11); see also \cite\eilmothr\ (p.~7).
The {\it Eilenberg-Koszul\/} sign convention 
is in force thoughout.
The degree of a homogeneous element $x$ of  a graded object
is written as $|x|$.
Given the chain complexes $U$ and $V$, the Hom-complex differential
on $\roman{Hom}(U,V)$ is as usual defined by
$$
D\varphi= d \varphi+ (-1)^{|\varphi|} \varphi d
$$
where $\varphi$ is a homogeneous, the convention being that
$\varphi\colon U_q \to V_p$ has degree $p-q$.
The identity morphism on an object is denoted by the same 
symbol as that object. The suspension operator is written as $s$;
given the chain complex $U$, the differential on the suspended object $sU$
is determined by the identity $ds+sd=0$, so that the suspension is a cycle in
the corresponding Hom-complex.
A {\it perturbation\/}
of the differential $d$ of a filtered chain complex
$Z$
is an operator
$\partial$
on $Z$ that
lowers filtration and, moreover, satisfies
$$
d \partial + \partial d + \partial \partial = 0,
$$
so that
$d_{\partial} = d + \partial$ is a new
differential on $Z$;
we shall then write $Z_{\partial}$ for the new chain complex.

Given a differential graded algebra $A$,
we will refer to a {\it differential graded\/}
left (or right) $A$-module more simply as a {\it left\/}
(or {\it right\/})
$A$-module; we shall denote the categories
of differential graded left $A$-modules by ${}_A\roman{Mod}$ and that of
differential graded
right $A$-modules by $\roman{Mod}_A$. 
We will write the multiplication map of a differential 
graded algebra $A$ as $\mu \colon A \otimes A \to A$
and the unit as $\eta \colon R \to A$
when there is a need to spell these structure maps out;
more generally, given a right $A$-module $N$, we will write the
structure map as $\mu \colon N \otimes A \to N$ and,
given a left $A$-module $N$, we will write the
structure map as $\mu \colon A \otimes N \to N$ as well.
An augmentation map for the differential graded algebra
$A$ is a morphism $\varepsilon\colon A \to R$ of differential graded algebras.
A differential graded algebra together with an augmentation map is defined 
to be {\it augmented\/}. Given the augmented differential graded algebra
$A$, with augmentation map  $\varepsilon\colon A \to R$, the {\it augmentation
ideal\/}, written as $IA$, is the kernel  of $\varepsilon$.
Given the augmented differential graded algebra $A$,
the {\it augmentation\/} filtration 
$\{\roman F^nA\}_{n \geq 0}$  given by
$$
\roman F^nA = (IA)^{\otimes n}\  (n \geq 0)
$$
turns $A$ into a filtered differential graded algebra, and
$A$ is {\it complete\/} when
the canonical map $A \longrightarrow \lim_n A\big /(IA)^n$
into the projective limit $\lim_A A\big /(IA)^n$
is an isomorphism. The augmentation filtration is descending and, 
at times, it is convenient to write it as
$\{\roman F_{n}A\}_{n \leq 0}$, where
$\roman F_{n}A= F^{-n}A$ ($n \leq 0$).
The differential graded algebra $A$ is said to be {\it connected\/}
when it is non-negative or non-positive and
when the unit $\eta\colon R \to A$ is an isomorphism onto the homogeneous
degree zero constituent $A_0$ of $A$. A connected augmented differential 
graded algebra is necessarily complete.  

Let $C$ be a differential graded coalgebra,
the diagonal map being written as $\Delta \colon C \to C \otimes C$
and the counit as $\varepsilon \colon C \to R$.
Given a left $C$-comodule $M$ we will likewise write the structure map as
 $\Delta \colon M \to C \otimes M$ and,
given a right $C$-comodule $M$, we will write the structure map as
$\Delta \colon M \to M \otimes C$.
A {\it coaugmentation map\/} for $C$ is a morphism
$\eta \colon R \to C$ of differential graded coalgebras
and a {\it coaugmented differential graded coalgebra\/}
is a differential graded coalgebra together with a coaugmentation map.
The {\it coaugmentation\/} coideal $JC$ is defined to be
the cokernel   $JC= \roman{coker}(\eta)$ of $\eta$.
Suppose that $C$ is coaugmented.
Recall that the counit $\varepsilon \colon C \to R$ and the
coaugmentation map determine a direct sum decomposition $C = R
\oplus JC$. The {\it coaugmentation\/} filtration 
$\{\roman F_nC\}_{n \geq 0}$ is as usual given by
$$
\roman F_nC = \roman{ker}(C \longrightarrow (JC)^{\otimes
(n+1)})\  (n \geq 0)
$$
where the unlabelled arrow is induced by some iterate of the
diagonal $\Delta$ of $C$. This filtration is ascending and is
well known to turn
$C$ into a {\it filtered\/} coaugmented differential graded
coalgebra; thus, in particular, $\roman F_0C = R$. We recall that
$C$ is said to be {\it cocomplete\/} when $C=\cup \roman F_nC$.
The reduced bar construction $\rbar A$ of an augmented differential graded
algebra $A$ is well known to be cocomplete.
We will refer to a {\it differential graded\/}
left (right) $C$-comodule
as a {\it left\/} ({\it right\/}) $C$-comodule;
we shall denote the categories
of left $C$-comodules by ${}_C\roman{Comod}$ and that of
right $C$-comodules by $\roman{Comod}_C$. 
The differential graded coalgebra $C$ is said to be connected when
it is non-negative or non-positive and when
the counit $\varepsilon\colon C \to R$ is an isomorphism from the 
homogeneous degree zero
constituent $C_0$ of $C$ to $R$.
A connected differential graded coalgebra is necessarily cocomplete.
A non-negative connected differential graded coalgebra $C$
is said to be {\it simply connected\/} when $C_1$ is zero.
Notice that the cobar construction $\rcob C$ on a simply connected differential graded
coalgebra $C$ is connected, and so is the
cobar construction $\rcob C$ on a non-positive differential graded
coalgebra $C$.
In the sequel, 
in a specific construction,
coaugmented differential graded coalgebras
will be cocomplete throughout but for clarity we will point out explicitly
whenever cocompleteness is needed.

Given the two chain complexes $U$ and $V$, possibly endowed with 
additional structure, a chain map from $U$ to $V$
possibly preserving additional structure
is referred to as a {\it quasi-isomorphism\/}
provided it induces an isomorphism on homology.
A chain complex which, in each degree, is
finitely generated, is said to be of {\it finite type\/}.

Let 
$(\{\roman F_p\},L):\ldots \subseteq \roman F_p \subseteq \roman F_{p+1} \subseteq \ldots \subseteq L$
be a filtered chain complex.
Recall that $(\{\roman F_p\},L)$ is
said to be {\it complete\/} when the canonical map $L \to \lim_p L/\roman F_p$
into the projective limit $\lim_p L/\roman F_p$
is an isomorphism of chain complexes, cf. \cite\eilmotwo\ (Section 4) where
the terminology $P$-complete is used and \cite\husmosta\ (I.3).
For completeness we recall that 
$(\{\roman F_p\},L)$ is
said to be {\it cocomplete\/} when the canonical map $\lim_p \roman F_p \to L$
from  the injective limit $\lim_p \roman F_p$ to $L$
is an isomorphism of chain complexes, cf. \cite\eilmotwo\ (Section 4) where
the terminology $I$-complete is used and \cite\husmosta\ (\S 0.5, I \S 3).

The {\it skeleton filtration\/} of the chain complex
$M$ is the filtration $\{\roman F_p(M)\}$ given by
$$
(\roman F_p(M))_n = \cases M_n, &\quad \text{for}\quad n \leq p ,\\
                    0, &\quad \text{for}\quad n > p ,
\endcases
$$
cf. \cite\husmosta\ (\S 0).

\head 1.1. Twisting cochains and twisted tensor products \endhead

Let $C$  be a coaugmented differential graded 
coalgebra, $A$ an augmented
differential graded algebra, $M$ a differential graded left $C$-comodule, and $N$ a
differential graded right $A$-module.
The familiar {\it cup pairing\/} $\cup$ turns
$\roman{Hom}(C,A)$ into a differential graded algebra
and 
$\roman{Hom}(M,N)$ into a differential graded {\it right\/}
$\roman{Hom}(C,A)$-module, with structure map of the kind
$$
\cup \colon 
\roman{Hom}(M,N)\otimes \roman{Hom}(C,A)\longrightarrow \roman{Hom}(M,N) .
\tag1.1.1
$$
The cup pairing assigns to $f \otimes h \in \roman{Hom}(M,N)\otimes \roman{Hom}(C,A)$ the morphism
$$
f \cup h\colon M @>{\Delta}>> 
M \otimes C @>{f \otimes h}>> 
N \otimes A
@>{\mu}>>
N.
$$
Likewise, the  {\it cap pairing\/}
$$
\cap\, \colon \roman{Hom}(C,A)\otimes N \otimes M
\longrightarrow N \otimes M
\tag1.1.2
$$
is given by the assignment to
$\varphi \in  \roman{Hom}(C,A)$ of
$$
\varphi \cap \,\cdot \,\,\colon N \otimes M 
@>{N\otimes \Delta}>> 
N \otimes C \otimes M
@>{N \otimes \varphi \otimes M}>> 
N \otimes A \otimes M
@>{\mu \otimes M}>>
N \otimes M.
$$
In the same vein, let $N$ be a differential graded left $A$-module
and $M$ a differential graded right $C$-comodule.
Under these circumstances, the {\it cap pairing\/}
$$
\cap\, \colon \roman{Hom}(C,A)\otimes M \otimes N 
\longrightarrow M \otimes N
\tag1.1.3
$$
is given by the assignment to
$\varphi \in  \roman{Hom}(C,A)$ of
$$
\varphi \cap \,\cdot \,\,\colon M \otimes N 
@>{\Delta \otimes N}>> 
M \otimes C \otimes N
@>{M \otimes \varphi \otimes N}>> 
M \otimes A \otimes N
@>{M \otimes \mu}>>
M \otimes N,
$$
For example, when $M=C$ and $N=A$, given $\varphi \in \roman{Hom}(C,A)$,
$c \in C$, and $a \in A$, when
$\Delta(c)=\sum c_j'\otimes c_j''$, in view of the Eilenberg-Koszul 
convention,
$$
\varphi \cap (c \otimes a)=\sum (-1)^{|\varphi||c_j'|} c_j'\otimes \varphi(c_j'')a.
$$
The cap pairing (1.1.3)
turns  $M \otimes N$ into a differential graded
{\it left\/} $\roman{Hom}(C,A)$-module and
the cap pairing
(1.1.2) turns  $N \otimes M$ into a differential graded
{\it left\/} $(\roman{Hom}(C,A))^{\roman{op}}$-module in the sense that,
given homogeneous $f$, $h$, $w$,
$$
(f\cup h)\cap w = (-1)^{|f||h|} h\cap f \cap w.
\tag1.1.4
$$
The last identity is an immediate consequence of the identity
$$
f\otimes h = (-1)^{|f||h|}(A\otimes h)\circ (f \otimes C)\colon 
C \otimes C \longrightarrow A \otimes A.
$$
Given the differential graded algebra
$\Cal A$, the notation $\Cal A^{\roman{op}}$ refers here to the
{\it opposite\/} differential graded algebra.
We note that,
since the underlying objects are graded,
the identity (1.1.4) does {\it not\/} say that
$N \otimes M$ acquires a right  $(\roman{Hom}(C,A))$-module structure.

A {\it twisting cochain\/} $\tau$ from $C$ to $A$
is a homogeneous morphism $\tau\colon C \to A$ of the underlying 
graded $R$-modules of degree $-1$ (i.~e. a homogeneous member of
$\roman{Hom}(C,A)$ of degree $-1$) such that
$$
D\tau=\tau \cup \tau,\ \tau \eta=0,\ \varepsilon \tau =0,
$$
where $D$ refers to the Hom-complex differential on $\roman{Hom}(C,A)$.

Let $\tau$ be a twisting cochain. Then the operator
$\partial^{\tau} =  -(\tau \cap \,\cdot \,)$ on $C \otimes A$ 
(i.~e. $\partial^{\tau}(w)=- \tau \cap w$, $w \in C \otimes A$)
is a {\it perturbation\/} of the differential on $C \otimes A$, that is,
$d +\partial^{\tau} $ is a differential on $C \otimes A$,
and the new differential is compatible with 
the differential graded left $C$-comodule and right $A$-module
structures. The resulting
differential graded left $C$-comodule and right $A$-module
is usually written as 
$C \otimes_{\tau} A$ and referred to as a {\it twisted tensor product\/}
or {\it construction\/} for $A$. More generally, given a right $C$-comodule
$M$ and a left $A$-module $N$, the twisted tensor product $M \otimes_{\tau}N$
is defined accordingly.
Likewise the operator
$\partial^{\tau} =  \tau \cap \,\cdot \,$ on $A \otimes C$ 
(i.~e. $\partial^{\tau}(w)= \tau \cap w$, $w \in A \otimes C$)
is a {\it perturbation\/} of the differential on $A \otimes C$, that is,
$d +\partial^{\tau} $ is a differential on $A \otimes C$,
and the new differential is compatible with 
the differential graded right $C$-comodule and left $A$-module
structures. The difference in sign is explained by the
structural descriptions of $C\otimes A$ and $A\otimes C$
as $(\roman{Hom}(C,A))$- vs.
$(\roman{Hom}(C,A))^{\roman{op}}$-modules explained above.

For later reference we note that,
given a differential graded left $A$-module $N$, the 
$A$-action on $N$ induces the
chain map
$$
C \otimes_{\tau} A \otimes N @>>> C \otimes_{\tau}  N
\tag1.1.5
$$
which induces an isomorphism
$$
(C \otimes_{\tau} A) \otimes_A N @>>> C \otimes_{\tau}  N
\tag1.1.6
$$
of chain complexes, cf. \cite\gugenhtw\ (2.6$_*$ Proposition).
In fact, cf. \cite\gugenhtw\ (2.4$_*$ Proposition),
the twisted differential on $C \otimes_{\tau}  N$
is the unique differential on $C \otimes N$ which makes (1.1.5)
into a chain map.

We recall that the twisting cochain $\tau$ induces a morphism 
$\overline \tau \colon \rcob C \to A$ of augmented differential graded algebras
and, when $C$ is cocomplete, a morphism
$\overline \tau \colon C \to \rbar A$, referred to as the {\it adjoints\/} of
$\tau$. The reason for this terminology is the fact that
the assignment to a twisting cochain $C \to A$
of its adjoint induces a natural bjection between the set 
$T(C,A)$ of twisting 
cochains from $C$ to $A$ and the set $\roman{Hom}_{\roman{alg}}(\rcob C,A)$
of morphisms of augmented differential graded algebras and,
when $C$ is cocomplete, 
the assignment to a twisting cochain $C \to A$
of its adjoint induces a natural bjection between $T(C,A)$ 
and the set $\roman{Hom}_{\roman{coalg}}(C,\rbar A)$
of morphisms of augmented differential graded coalgebras.
Thus the functors $\rcob$ and $\rbar$ are adjoint functors between
the category of cocomplete coaugmented differential
graded coalgebras and that of augmented differential graded algebras.
Some connectivity assumptions are necessary here to make this adjointness
precise; see e.~g. \cite\gugenmun\ for details.

We will say that the twisting cochain $\tau$ is {\it acyclic\/} when
$C \otimes _{\tau}A$ is an
acyclic complex (and hence an
acyclic construction, cf. e.~g. \cite\mooretwo\ for this notion).
Subject to mild appropriate additional
hypotheses of the kind that $C$ and/or $A$ be projective as
$R$-modules---which will always hold in the paper---the adjoints
$\overline \tau \colon C \to \rbar A$ and
$\overline \tau \colon \rcob C \to A$ of an acyclic twisting cochain
$\tau$ are chain equivalences.
We denote the universal bar construction twisting cochain 
by $\tau^{\rbar A}\colon \rbar A \to A$
and the universal cobar construction twisting cochain by
$\tau_{\rcob C}\colon C \to \rcob C$; these are acyclic.
Occasionally we write $\tau^{\rbar}$ and $\tau_{\rcob}$ rather than
$\tau^{\rbar A}$ and $\tau_{\rcob C}$, respectively.

Let $\tau_1$ and $\tau_2$ be two twisting cochains 
from $C$ to $A$. A {\it homotopy\/} of twisting cochains from
$\tau_1$ to $\tau_2$, written as $\psi\colon \tau_1 \simeq \tau_2$,
is a homogeneous morphism $\psi \colon C \to A$ of degree
zero such that
$$
D\psi = \tau_1 \cup \psi - \psi \cup \tau_2, \quad \psi \eta =
\eta, \quad \varepsilon \psi = \varepsilon
\tag1.1.7
$$
where $D$ refers to the ordinary Hom-complex differential. Here
the notation $\eta$ and $\varepsilon$ is slightly abused in the
sense that $\varepsilon$ denotes the counit of $C$ as well as the
augmentation map of $A$ and that $\eta$ denotes the unit of $A$
and the coaugmentation map of $C$.

Let $M$ be a right $C$-comodule, $N$ a left $A$-module,
and let $\psi$ be a homotopy of twisting cochains from $\tau_1$ to
$\tau_2$. Then
$$
\psi \cap \,\cdot \,\colon M \otimes _{\tau_2} N \longrightarrow 
M \otimes _{\tau_1} N
\tag1.1.8
$$
is a {\it morphism\/} of twisted tensor products.
Indeed,
$$
\left(d +\partial^{\tau_1}\right)(\psi \cap \,\cdot \,)
+
(\psi \cap \,\cdot \,)
\left(d +\partial^{\tau_2}\right)
=
\left(D\psi+ \psi \cup \tau_2-\tau_1 \cup \psi\right) \cap \,\cdot \,
= 0.
$$
Notice when $\psi=\eta\varepsilon$, the morphism $\psi \cap \,\cdot \,$ is
the identity. Suppose that $C$ is cocomplete.
Then (1.1.8) is in fact an {\it isomorphism\/}.
The inverse map is obtained by a standard procedure:
Write $\psi = \eta \varepsilon + \widetilde \psi$; 
the infinite series
$$
\psi^{-1} = \eta \varepsilon - \widetilde \psi
+\widetilde \psi \cup \widetilde \psi
-\widetilde \psi^{\cup 3} + \ldots
\tag1.1.9
$$
converges since $C$ is cocomplete, and
$$
\psi^{-1} \cap \,\cdot \,\colon M \otimes _{\tau_1} N \longrightarrow 
M \otimes _{\tau_2} N
\tag1.1.10
$$
yields the inverse for (1.1.8).

At this stage, one can build a category of twisted tensor products
with general notions of morphism and isomorphism; for our purposes,
the notion of isomorphism of the kind (1.1.8) suffices.

\head 1.2. Twisted Hom-objects and complete twisted Hom-objects\endhead

Let $C$ be a coaugmented differential graded coalgebra, $A$ an augmented
differential graded algebra, and $\tau \colon C \to A$
a twisting cochain.

Let $N$ be a differential graded right
$A$-module (we could equally well work with a left $A$-module) and
let $\delta^{\tau}$ be the operator on $\roman{Hom}(C,N)$
given, for homogeneous $f$, by $\delta^{\tau}(f) = (-1)^{|f|}f \cup \tau$.
With reference to the filtration induced by the coaugmentation
filtration of $C$, the operator $\delta^{\tau}$ is a {\it perturbation\/}
of the differential $d$ on $\roman{Hom}(C,N)$, and we write the perturbed
differential on $\roman{Hom}(C,N)$ as $d^{\tau}= d + \delta^{\tau}$.
We will refer to $\roman{Hom}^{\tau}(C,N)=(\roman{Hom}(C,N),d^{\tau})$
as a {\it twisted\/} Hom-{\it object\/}.
Under suitable circumstances,
$\roman{Hom}^{\tau}(C,N)$ calculates the differential graded
$\roman{Ext}_A(R,N)$.

\proclaim{Lemma 1.2.1} Suppose that $C$ is cocomplete.
Let $\tau_1$ and $\tau_2$ be twisting cochains from $C$ to $A$
and let $\psi$ be a homotopy of twisting cochains from $\tau_1$ to
$\tau_2$.
The assignment to a (homogeneous) $\phi \in
\roman{Hom}(C,N)$ of $\phi \cup \psi \in \roman{Hom}(C,N)$ yields 
a morphism, in fact,
isomorphism
$$
\roman{Hom}^{\tau_1}(C,N) @>>> \roman{Hom}^{\tau_2}(C,N)
\tag1.2.2
$$
of twisted $\roman{Hom}$-objects,
with inverse given by the
assignment to a (homogeneous) $\phi \in
\roman{Hom}(C,N)$ of $\phi \cup \psi^{-1} \in \roman{Hom}(C,N)$.
\endproclaim

\demo{Proof}
Indeed,
$$
\align 
d^{\tau_2}(\phi \cup \psi) &= (d\phi) \cup \psi +(-1)^{|\phi|}(\phi \cup d\psi
+\phi \cup \psi \cup \tau_2)
\\
&= (d\phi) \cup \psi +(-1)^{|\phi|}\phi \cup \tau_1 \cup \psi
\\
&=(d\phi +(-1)^{|\phi|}\phi \cup \tau_1) \cup \psi
\\
&=(d^{\tau_1}(\phi))\cup \psi
\endalign
$$
whence the assertion. \qed
\enddemo

One can now build a category of twisted Hom-objects
with general notions of morphism and isomorphism; for our purposes,
the notion of isomorphism of the kind (1.2.2) suffices.

Consider the canonical injection of
$\roman{Hom}(C,N)$ into $\roman{Hom}(C\otimes A,N)$ which assigns
$$
\Phi_{\alpha} \colon C\otimes A @>>> N,
\quad
\Phi_{\alpha} (w\otimes a)= \alpha(w) a,\ w \in C,\ a \in A,
\tag1.2.3
$$
to $\alpha \in \roman{Hom}(C,N)$; this injection plainly identifies
$\roman{Hom}(C,N)$ with the
subspace of (right) $A$-module morphisms
from $C\otimes A$ to $N$. Furthermore,
for homogeneous $\phi \in \roman{Hom}(C,A)$,
$$
\Phi_{\alpha} \circ (\phi \cap \,\cdot \, ) 
= \Phi_{\alpha \cup \phi}
\colon
C \otimes A @>>> N,
\tag1.2.4
$$
where the notation \lq\lq $\cup$\rq\rq\ in
the expression 
$\alpha \cup \phi \in \roman{Hom}(C,N) $ refers to the cup pairing (1.1.1).
Consequently
the assignment to $\alpha$ of
$\Phi_{\alpha}$ yields a morphism
$$
\roman{Hom}^{\tau}(C,N) @>>>\roman{Hom}(C\otimes_{\tau} A,N)
\tag1.2.5
$$
even of chain complexes, and this morphism identifies
$\roman{Hom}^{\tau}(C,N)$ with the subspace
$\roman{Hom}_A(C\otimes_{\tau} A,N)$
of differential graded (right) $A$-module morphisms
from $C\otimes_{\tau} A$ to $N$.

Let $M$ be  a differential graded left $A$-module.
The $R$-dual $M^*$ of $M$ inherits a canonical differential
graded right $A$-module structure, and the twisted Hom-object
$\roman{Hom}^{\tau}(C,M^*)$ is defined. 
The canonical assignment to $\Phi \colon C\otimes M \to R$
of $\alpha_{\Phi}  \colon C \to M^*$ where
$$
\alpha_{\Phi}(c)(x) = \Phi(c \otimes x),\ c \in C,\ x \in M,
$$
yields the adjointness isomorphism
$$
\roman{Hom}(C\otimes_{\tau} M,R)\longrightarrow \roman{Hom}^{\tau}(C,M^*),
\tag1.2.6
$$
compatible with the perturbed differentials as indicated. 
Indeed, the perturbation $D^{\tau}$ of the Hom-differential $D$ is given by
$$
D^{\tau}(\Phi) = (-1)^{|\Phi|+1} \Phi \circ (-\tau \cap \,\cdot\,)
 = (-1)^{|\Phi|} \Phi \circ (\tau \cap \,\cdot\,)
$$
whereas, for $\alpha \in \roman{Hom}(C,M^*)$,
$$
\delta^{\tau}(\alpha) =(-1)^{|\alpha|}\alpha \cup \tau
$$
whence, indeed, (1.2.6) is compatible with the differentials.

Under these circumstances,
when $A$ is of finite type,
its dual $A^*$ inherits a
coaugmented differential graded coalgebra structure and
$M^*$ an obvious left $A^*$-comodule structure,
and the dual $\tau^* \colon A^* \to C^*$ is a twisting cochain.
Nota bene: In view of the Eilenberg-Koszul convention, this dual is given by
$$
\tau^*(\alpha) = (-1)^{|\alpha|} \alpha \circ \tau,\ \alpha \colon A \to R.
$$
Furthermore, the canonical morphism of chain complexes
$$
C^* \otimes_{\tau^*} M^* @>>> \roman{Hom}^{\tau}(C,M^*)
\tag1.2.7
$$
is
a morphism of differential graded $C^*$-modules
and $A^*$-comodules; it is
an isomorphism when $C$ is of finite type.
The twisted Hom-object $\roman{Hom}^{\tau}(C,M^*)$ is always defined, though,
whether or not $A$ is of finite type. This twisted object
is a kind of {\it completed twisted object\/} 
associated with the data $A^*$, $C^*$, $M^*$, and $\tau^*$,
where the expression \lq\lq associated with\rq\rq\ is to be interpreted 
with a grain of salt since  
$\roman{Hom}^{\tau}(C,M^*)$ involves $C$ and $\tau$ rather than 
merely $C^*$ and $\tau^*$
and
$C$ and $\tau$ are not necessarily determined by $C^*$ and $\tau^*$.

We will now briefly discuss the notion of what we will refer to as
a {\it complete twisted Hom-object\/}.
Let $M$ be a differential graded left
$C$-comodule (we could equally well work with a right $C$-comodule)
and $N$ a differential graded right $A$-module.
Similarly as under the circumstances of
(1.2.6), the twisted differential on the right-hand side of the 
adjointness isomorphism
$$
\roman{Hom}(N,M^*) \cong \roman{Hom}(N \otimes M,R)
$$
dual to the twisted differential on $N \otimes_{\tau} M$
determines a perturbed differential $d^{\tau}$ on
$\roman{Hom}(N,M^*)$; we write the resulting twisted object as
$\roman{Hom}^{\tau}(N,M^*)$ and, similarly as before,
refer to it as a {\it twisted
Hom-object\/}. This situation can be formalized in terms of
$C^*$-modules in the following way: We define a {\it complete\/}
left $C^*$-module to be a chain complex $\Mfl$ together with
an operation $\rho\colon \roman{Hom}(C,\Mfl) \to \Mfl$
which satisfies the obvious associativity and unit constraints.
The associativity constraint says that the composite
$$
\roman{Hom}(C,\roman{Hom}(C,\Mfl))
@>{\roman{Hom}(C,\rho)}>>\roman{Hom}(C,\Mfl)
@>{\rho}>> \Mfl
$$
equals the composite
$$
\roman{Hom}(C \otimes C,\Mfl)
@>{\roman{Hom}(\Delta,\Mfl)}>>\roman{Hom}(C,\Mfl)
@>{\rho}>> \Mfl
$$
where $\roman{Hom}(C,\roman{Hom}(C,\Mfl))$ and $\roman{Hom}(C
\otimes C,\Mfl)$ are identified by adjointness in the
standard manner; here $\Delta$ is the diagonal map of $C$, and an
object and the identity morphism on that object are denoted by the
same symbol. Under these circumstances, $
\roman{Hom}(C,\Mfl)$ is a fortiori a left $C^*$-module  and,
when $C$ is of finite type, the operation $\rho$ comes down to an
ordinary left $C^*$-module structure $C^*\otimes \Mfl \to
\Mfl$. We denote the category of complete left $C^*$-modules
by ${}_{C^*}\widehat{\roman{Mod}}$. The dual $M^*$ of a
$C$-comodule $M$ is a complete $C^*$-module in an obvious manner.

Let $\Mfl$ be a complete $C^*$-module,
with structure map $\rho$, and let $N$ be a differential graded
right $A$-module. For
$\phi \colon N \to \Mfl$, define
$\delta^{\tau}(\phi)$ to be the composite of
$$
\roman{Hom}(N,\Mfl)
@>{\mu^*}>>
\roman{Hom}(N \otimes A,\Mfl)
@>>>
\roman{Hom}(N,\roman{Hom}(A,\Mfl))
$$
where the unlabelled arrow comes from adjointness, with
$$
\roman{Hom}(N,\roman{Hom}(A,\Mfl))
@>{\tau^{\sharp}}>>
\roman{Hom}(N,\roman{Hom}(C,\Mfl))
@>{\rho_*}>>
\roman{Hom}(N,\Mfl);
$$
here the notation $\tau^{\sharp}$ is intended to indicate that that morphism
is induced from $\tau$ in the obvious manner
(but we cannot simply write 
$\tau^*$ since this would conflict with notation established
before).
The operator $\delta^{\tau}$ is a
perturbation of the naive differential on
$\roman{Hom}(N,\Mfl)$, and we will denote the perturbed object by
$\roman{Hom}^{\tau}(N,\Mfl)$.
Under suitable circumstances,
$\roman{Hom}(A,\Mfl)$ calculates the differential graded
$\roman{Tor}_{C^*}(R,\Mfl)$.

The structure map $\rho$ of $\Mfl$ induces a morphism
$\rho_*$ from $\roman{Hom}(C\otimes N,\Mfl)$ to
$\roman{Hom}(N,\Mfl)$ which is, in fact, a chain map
$$
\rho_*\colon
\roman{Hom}(C\otimes_{\tau} N,\Mfl)
@>>>
\roman{Hom}^{\tau}(N,\Mfl);
\tag1.2.8
$$
the latter, in turn, passes through a chain map
$$
\roman{Hom}_{C^*}(C\otimes_{\tau} N,\Mfl)
@>>> \roman{Hom}^{\tau}(N,\Mfl).
\tag1.2.9
$$
This discussion applies, in particular,
to the complete $C^*$-module $M^*$ relative to  the obvious structure
coming from the $C$-comodule structure on $M$.
Moreover, when $A$ is of finite type, its dual $A^*$ inherits a
coaugmented differential graded coalgebra structure,
the dual $\tau^* \colon A^* \to C^*$ is a twisting cochain,
and the canonical morphism of chain complexes
$$
A^* \otimes_{\tau^*} M^* @>>> \roman{Hom}^{\tau}(A,M^*)
$$
is an isomorphism of differential graded $\roman{Hom}^{\tau}(A^*,C^*)$-modules.
When $C$ is, furthermore, of finite type,
the chain map
(1.2.8) amounts to
$$
\rho_*\colon
A^*\otimes_{\tau^*} C^* \otimes M^*
@>>>
A^* \otimes_{\tau^*} M^*
$$
and (1.2.9) cones down to an isomorphism
$$
(A^*\otimes_{\tau^*} C^*) \otimes_{C^*} M^*
@>>> A^* \otimes_{\tau^*} M^* .
$$
These are special cases of (1.1.5) and
(1.1.6) above.

The twisted Hom-object $\roman{Hom}^{\tau}(A,M^*)$ is always defined, though,
whether or not $A$ is of finite type; 
in a sense explained above, that twisted Hom-object
is then a kind of {\it completed twisted object\/} 
associated with the data $A^*$, $C^*$, $M^*$, and $\tau^*$.

\head 1.3. The trivial simplicial object and the total simplicial
object\endhead

Any object $X$ of a symmetric monoidal category 
endowed with a
diagonal---we will take the categories of spaces, of
smooth manifolds, of groups, of vector spaces, of Lie algebras,
etc.,---defines two simplicial objects in the category, the {\it trivial\/}
object which, with an abuse of notation,
we still write as $X$, and the {\it total object\/} $EX$
(\lq\lq total object\rq\rq\ 
not being standard terminology in this generality); the
trivial object $X$ has a copy of $X$ in each degree and all simplicial
operations are the identity while, for $p \geq 0$, 
the degree $p$  constituent
$EX_p$ of the total object $EX$ is a product of
$p+1$ copies of $X$ with the familiar 
face operations given by omission and degeneracy operations
given by insertion; 
somewhat more formally, {\it insertion\/} is here interpreted as
an insertion of the diagonal morphism at the appropriate place.
When $G$ is a group, this kind of construction
leads to the {\it unreduced homogeneous\/} bar construction where
the term  \lq\lq homogeneous\rq\rq\ refers to the fact that the 
formulation uses the group  structure only for the module structure
relative to the group and {\it not\/} for the simplicial structure.
Well known
categorical machinery formalizes this situation but we shall not need this
kind of formalization. For intelligibility we recall that in \cite\bottone\ 
the functor which we write as $E$ is denoted by $P$.

\head 1.4. The geometric resolution and nerve
construction\endhead

Let $G$ be a topological group.
The universal simplicial $G$-bundle $EG \to NG$
is a special case of the construction spelled out in (1.3) above; 
indeed, $EG$ acquires a
simplicial group structure, the diagonal maps combine to a morphism of
simplicial groups $G \to EG$ by means of which $EG$ 
acquires a simplicial principal right $G$-space structure, 
as a simplicial space, $NG= EG\big /G$, and the projection 
$EG \to NG$ is the obvious map.

Let $X$ a left $G$-space. The
{\it simplicial bar construction\/} $N(G,X)$ is a simplicial space whose
geometric realization $|N(G,X)|$ is the Borel construction $EG \times_GX$;
with reference to the obvious filtration, $EG \times X$ is, in particular, a
free acyclic $G$-resolution \cite\rothstee\ of $X$ in the category of left
$G$-spaces. In particular, when $X$ is a point, $N(G,X)$ comes down to the
ordinary {\it nerve\/} $NG$ of $G$, the (lean) realization of which is the
ordinary classifying space $BG$ of $G$. (Within our terminology, it would be
more  consistent to write $\overline BG$ for the classifying space but we
will stick to the standard notation $BG$.) 
We also note at this stage that, for the sake of clarity, 
we use the font $\roman B$ for the {\it algebraic\/} bar construction.
For general $X$,
the $G$-{\it equivariant
homology\/} and {\it cohomology\/} of $X$ can then be defined as the homology and
cohomology of $N(G,X)$. This description works for {\it any\/} (co)homology
theory, e.~g. singular (co)homology or de Rham cohomology: Application of the
(co)chain functor to $N(G,X)$ yields a (co)simplicial chain complex the
(co)homology of which is the $G$-equivariant theory. Thus, to define the
equivariant theory, there is {\it no need\/} 
to pass through the {\it Borel construction\/},
cf. e.~g. the discussion in \cite\bottone; in particular, this is the correct
approach to equivariant de Rham theory since the de Rham functor does not
naively apply to the Borel construction.
In Theorem 4.1 below we will recall a characterization of equivariant
(co)homology as suitable differential Tor- and Ext-functors.
When this characterization is taken as the definition,
chains, cochains and forms on $N(G,X)$ then appear as
an a posteriori object for the calculation of the
$G$-equivariant (co)homology of the $G$-space $X$.

It is worthwhile noting that the
constructions to be given below apply to the {\it two-sided\/} constructions at
no extra cost: Let $E$ be a right $G$-space.
The obvious
action of  $G$ on the simplicial space $E \times EG \times X$ is
principal, and the
{\it two-sided  simplicial bar construction\/}
$N(E,G,X)$ is the base of the resulting simplicial
principal $G$-bundle. We will
always suppose that the projection $E \to B=E/G$ is a principal bundle. It is
a classical fact that then the canonical projection from $|N(E,G,X)|$ to
$E \times_G X$  is a homotopy equivalence which is
{\it compatible with the bundle structures\/} in the appropriate sense.

\head 1.5. The notion of model\endhead 

Given a chain complex $U$, possibly with additional structure
(module over a differential graded algebra, comodule over
a differential graded coalgebra, twisted object, etc.),
we will refer to a chain complex $V$,
possibly endowed with additional structure,
as a {\it model\/} for $U$ provided
there is a chain
$U=U_0$, $U_1$, \dots, $U_m=V$ of chain complexes,
possibly endowed with the additional structure,
and a chain of morphisms
$$
U_0 \leftarrow U_1,
\
U_1 \rightarrow U_2,
\
U_2 \leftarrow U_3,
\
\ldots,
U_{m-1} \rightarrow U_m
$$
of chain complexes,
each being a quasi-isomorphism and
possibly compatible with the additional structure.
Under such circumstances, each of these various unlabelled arrows 
will occasionally be 
referred to as a {\it comparison map\/}.

\head 1.6. Sh-modules and sh-comodules\endhead

Given two augmented differential graded algebras $A_1$ and $A_2$,
in view of the structure of the reduced bar construction,
a twisting cochain $\tau\colon \rbar A_1\to A_2$
has the homogeneous constituents $\tau_j\colon (sIA)^j \to A_2$
($j \geq 1$),
and a twisting cochain $\tau$ having $\tau_j$ zero
for $j \geq 2$ amounts to an ordinary morphism of
differential graded algebras.
Non-zero higher terms are an instance of higher homotopies,
and a general twisting cochain is an instance of an sh-map,
where \lq\lq sh\rq\rq\  stands for \lq\lq strongly homotopic\rq\rq.
Details about the development of these ideas
and about the history behind can be found in \cite\origins.

Given two augmented differential graded algebras $A_1$ and $A_2$,
an sh-{\it map\/} from $A_1$ to $A_2$
is defined to be a twisting cochain $\tau \colon \rbar A_1 \to A_2$.
Likewise,
given two coaugmented differential graded coalgebras $C_1$ and $C_2$,
an sh-{\it map\/} from  $C_1$ to $C_2$
is a twisting cochain $\tau \colon  C_1 \to \rcob C_2$.

Let $C$ be a coaugmented differential graded coalgebra,
$A$ an augmented differential graded algebra, and let
$\tau \colon C \to A$ be an acyclic twisting cochain.
We will refer to
a left (right) $A'$-module $N$, for some augmented
differential graded algebra $A'$,
together with a twisting cochain
$\tau_{A'}\colon C \to A'$, as
a {\it left\/} ({\it right\/}) sh-{\it module over \/} $A$.
An ordinary $A$-module is an 
sh-module over $A$ in an obvious manner.
Thus the notion of sh-module over $A$ extends that of ordinary $A$-module.
When $(N,\tau_{A'})$ is an
sh-{\it module over \/} $A$, 
the chain complex $N$ acquires an $(\rcob C)$-module structure via
the adjoint $\overline \tau_{A'}\colon \rcob C \to A'$, and we may
in particular take $A'= \rcob C$ and $\tau_{A'} = \tau_{\rcob C}$;
thus requiring the twisting cochain $\tau_{A'}\colon C \to A'$ to be acyclic
yields an equivalent notion of sh-module over  $A$.
The given definition allows for more flexibility, though.
Likewise a {\it left\/}  ({\it right\/}) sh-{\it comodule over\/}
$C$ is a left (right) $C'$-comodule $M$, for some 
{\it cocomplete\/} coaugmented
differential graded coalgebra $C'$, together with a twisting cochain
$\tau^{C'}\colon C' \to A$.
An ordinary $C$-comodule is an 
sh-comodule over $C$ in an obvious manner.
Thus the notion of sh-comodule over $C$ extends that of ordinary 
$C$-comodule.
When $(M,\tau^{C'})$ is an
sh-{\it comodule over \/} $C$, 
the chain complex $M$ acquires a $(\rbar A)$-comodule structure via
the adjoint $\overline \tau^{C'}\colon C' \to \rbar A$, and we may
in particular take $C'= \rbar A$ and $\tau^{C'} = \tau^{\rbar A}$;
again requiring the twisting cochain $\tau^{C'}\colon C' \to A$
to be acyclic yields an equivalent notion of sh-comodule over  $C$,
and the given definition allows for more flexibility.

Given two sh-modules $(N_1,\tau_{A_1'})$ and
$(N_2,\tau_{A_2'})$ over $A$,
we {\it define\/} an {\it sh-morphism\/}
from $(N_1,\tau_{A_1'})$ to $(N_2,\tau_{A_2'})$ to be a
morphism $N_1 \to N_2$ of differential graded $(\rcob C)$-modules
where the $(\rcob C)$-module structures on $N_1$ and $N_2$ are induced
from the twisting cochains
$\tau_{A_1'}$  and $\tau_{A_2'}$, respectively.
Denote the resulting category by
${}_A\roman{Mod}^{\infty}$ or
$\roman{Mod}_A^{\infty}$, according as whether left or right modules
are involved.
An ordinary $A$-module morphism is an sh-morphism in an obvious manner.
Thus the notion of  sh-morphism extends that of ordinary
$A$-module morphism.

Likewise,
given two sh-comodules $(M_1,\tau^{C_1'})$ and
$(M_2,\tau^{C_2'})$ over $C$,
we {\it define\/} an {\it sh-morphism\/}
from
$(M_1,\tau^{C_1'})$ to $(M_2,\tau^{C_2'})$
to be a
morphism $M_1 \to M_2$ of differential graded
$(\rbar A)$-comodules
where the $(\rbar A)$-comodule structures on $M_1$ and $M_2$ are induced
from the twisting cochains
$\tau^{C_1'}$  and $\tau^{C_2'}$, respectively.
Denote the resulting category by
${}_C\roman{Comod}^{\infty}$ or
$\roman{Comod}_C^{\infty}$, according as whether left- or right comodules
are under discussion.
An ordinary morphism of $C$-comodules
is an sh-morphism in an obvious manner.
Thus the notion of sh-morphism extends that of ordinary
morphism of $C$-comodules.

As before, let $\tau\colon C \to A$ be an acyclic twisting
cochain. Further,  let $C'$ be a coaugmented differential graded
coalgebra, $A'$ an augmented differential graded algebra,
and let $\zeta_{A'} \colon C \to A'$ and $\zeta^{C'} \colon C' \to
A$ be acyclic twisting cochains.
The following is straightforward.

\proclaim{Proposition 1.6.1${}_*$} The assignment  to a left
$A'$-module $N$ of the left sh-module $(N,\zeta_{A'})$ over $A$ is
a functor from $ {}_{A'}\roman{Mod}$ to
${}_{A}\roman{Mod}^{\infty} $ and, likewise, the assignment to a
left $C'$-comodule $M$ of the left sh-comodule $(M,\zeta^{C'})$
over $C$ is a functor from $ {}_{C'}\roman{Comod}$ to
${}_{C}\roman{Comod}^{\infty}$. \qed
\endproclaim

\medskip\noindent {\bf 2. Homological perturbations}
\smallskip\noindent
{\it Homological perturbation theory\/} (HPT) is concerned with transferring
various kinds of algebraic structure through a homotopy equivalence.
Historical comments about the development of HPT may be found
in \cite\origins\ and \cite\tornike. The basic reason why HPT works
is the old observation that an exact sequence of chain complexes
which splits as an exact sequence of graded modules and which has
a contractible quotient necessarily splits in the category of
chain complexes \cite\doldone\ (2.18).

Here is an essential piece of machinery.

\noindent
{\smc Definition 2.1.}
A {\it contraction\/}
$$
\Nsddata N{\nabla}{\pi}Mh
\tag2.1.1
$$
of chain complexes consists of chain complexes $N$ and $M$, chain maps
$\pi\colon N \to M$ and $\nabla \colon M \to N$,
and a morphism $h\colon N \to N$
of the underlying graded modules of degree 1;
these data are required to satisfy the identites
$$
\align
\pi \nabla &= \roman{Id},
\quad
Dh = \nabla \pi - \roman{Id},
\tag2.1.2
\\
\pi h &= 0, \quad h \nabla = 0,\quad hh = 0.
\tag2.1.3
\endalign
$$
We will then say that $N$ {\it contracts onto\/} $M$.
If furthermore, $N$ and $M$ are filtered chain complexes, and if
$\pi$, $\nabla$ and $h$ are filtration preserving,
the contraction is said to be {\it filtered\/}.
The requirements (2.1.3) are referred to as {\it annihilation
properties\/} or {\it side conditions\/}.

The notion of contraction was introduced in \S 12 of \cite\eilmactw;
it is among the basic tools in homological perturbation theory,
cf. \cite{\huebkade} and the literature there.

\smallskip
\noindent{\smc Remark 2.1.4.} It is well known that the side 
conditions (2.1.3)
can always be achieved. This fact relies on the standard observation that
a chain complex is contractible if and only
if it is isomorphic to a cone, cf. \cite\husmosta\  (IV.1.5).
This observation is, by the way, related with that of Dold's quoted above.
Under the present circumstances,
given data of the kind (2.1.1) such
that the identities (2.1.2)  hold but not necessarily the side
conditions (2.1.3), the operator
$$
\widetilde h=(\roman{Id} -\nabla \pi)h(\roman{Id} -\nabla
\pi)d(\roman{Id} -\nabla \pi)h(\roman{Id} -\nabla \pi)
$$
satisfies the requirements (2.1.2) and (2.1.3), with
$\widetilde h$ instead of $h$;
when $h$ already satisfies (2.1.3), $\widetilde h$ coincides
with $h$.

The argument in \cite\gugenmun\ (4.1) establishes the following:

\proclaim{Proposition 2.2$_*$}
Let $C$ and $C'$ be coaugmented differential graded coalgebras,
let
${\Nsddata {C'} {\nabla}{\pi}{C}h}$
be a contraction of chain complexes
where $\nabla$ is a morphism of differential
graded coalgebras, and let
$A$ be a complete
augmented differential graded algebra.
Given a twisting cochain
$\sigma \colon C \to A$, there is a unique twisting cochain
$\xi \colon C' \to A$
with $\xi \nabla = \sigma$ and $\xi h = 0$.
This twisting cochain is given by the inductive formula
$$
\xi = \sigma \pi - (\xi\cup\xi) h
\tag2.2.1${}_*$
$$
and is natural in terms of the data. \qed
\endproclaim

Somewhat more explicitly, the formula (2.2.1${}_*$) means that
$\xi = \xi_1 + \xi_2 + \dots $ where
$$
\xi_1 =\sigma \pi, \ \xi_p = -\sum_{i+j=p} (\xi_i \cup \xi_j) h\quad (p>1).
$$
For any $n\geq 1$, the composite of $\xi_p$ with the canonical 
projection $A\to A\big /(IA)^n$ is non-zero for only finitely many $p$
whence the construction yields a twisting cochain with values in the
projective limit $\lim_n A\big /(IA)^n$. Since $A$ is supposed to be complete,
the construction actually yields a twisting cochain with values in $A$ itself.

The dual statement reads as follows.

\proclaim{Proposition 2.2$^*$} Let $A$ and $A'$ be augmented
differential graded algebras, let ${\Nsddata {A'}
{\nabla}{\pi}{A}h}$ be a contraction of chain complexes where
$\pi$ is a morphism of differential graded algebras, and let $C$
be a cocomplete coaugmented differential graded coalgebra. Given a
twisting cochain $\sigma \colon C \to A$, there is a unique
twisting cochain $\xi \colon C \to A'$ with $\pi \xi = \sigma$ and
$h\xi=0$. This twisting cochain is given by the inductive formula
$$
\xi = \nabla \sigma  - h(\xi\cup\xi)
\tag2.2.1$^*$
$$
and is natural in terms of the data.  \qed
\endproclaim

Here, for each $n\geq 1$, on the degree $n$ 
coaugmentation filtration constituent $\roman F_nC$ of $C$,
only finitely many of the recursive terms of $\xi$ are non-zero.
Since $C$ is cocomplete, each element of $C$ actually lies in some
$\roman F_nC$ whence, given a particular element $c$ of $C$, the convergence of
$\xi(c)$ is naive.

\proclaim{Lemma 2.3}{\rm [Perturbation lemma]} Let $\Nsddata N{\nabla}gMh$ be a
filtered contraction, let $\partial$ be a perturbation of the
differential on $N$, and let
$$
\align
\Cal D &= \sum_{n\geq 0} g\partial (h\partial)^n\nabla =
\sum_{n\geq 0} g(\partial h)^n\partial\nabla
\tag2.3.1
\\
\nabla_{\partial}&= \sum_{n\geq 0} (h\partial)^n\nabla
\tag2.3.2
\\
g_{\partial}&= \sum_{n\geq 0} g(\partial h)^n
\tag2.3.3
\\
h_{\partial}&=\sum_{n\geq 0} (h\partial)^n h =\sum_{n\geq 0} h(\partial h)^n .
\tag2.3.4
\endalign
$$
If the filtrations on $M$ and $N$ are complete,
these infinite series converge,
$\Cal D$ is a perturbation of the differential on
$M$ and, if we write
$N_{\partial}$ and
$M_{\Cal D}$
for the new chain complexes,
$$\Nsddata {N_{\partial}}{\nabla_{\partial}}{g_{\partial}}{M_{\Cal D}}
{h_{\partial}}
\tag2.3.5
$$
constitute a new filtered contraction
that is natural in terms of the given data.
\endproclaim
\demo{Proof} Details may be found in \cite{\brownez}\  and
in Lemmata 3.1 and 3.2 of
\cite{\gugenhtw}.  \qed
\enddemo

\noindent
{\smc Remark 2.4.\/} Given an arbitrary augmented
differential graded algebra $A$,
the construction in \cite\munkholm\ 
(2.14 Proposition, 2.15 Corollary) yields a contraction
$$
\Nsddata{A}{\rcob \overline \tau}{\nabla}{\rcob \rbar A}h
$$
that is natural in $A$ in such a way that, by construction,
$\rcob \overline \tau$ is a morphism of differential graded algebras.
See also \cite\husmosta\ (II.4.4 Theorem p.~148).
The corresponding result for a connected coaugmented differential
graded coalgebra is given in \cite\husmosta\ (II.4.5 Theorem p.~148).
The special cases for a connected algebra $A$ 
and a simply connected coalgebra $C$ can be found already in \cite\drachman.

\beginsection 3. Duality

Let $C$ be a coaugmented differential graded coalgebra,
$A$ an augmented differential graded algebra, and
$\tau \colon C \to A$ an acyclic twisting cochain.
Suppose there be given an explicit contracting homotopy 
$s\colon C \otimes_{\tau}A \to C \otimes_{\tau}A$ on
$C \otimes_{\tau}A$ as well, so that 
$\Nsddata{C\otimes_\tau A}{\eta}{\varepsilon}{R}s$
is a contraction, i.~e.
$Ds = \roman{Id}-\eta \varepsilon$, 
$ss=0$, etc. In view of Remark 2.1.4, once a contracting
for $C \otimes_{\tau}A$ is given, the requirement
$ss=0$  can always be arranged for.
With a slight abuse of notation, we will
denote the corresponding contracting homotopy of
$A \otimes_{\tau}C$ by $s$ as well.
These contracting homotopies will be used in Theorems 3.5$_*$ and
3.5$^*$ below.

We suppose that the cohomology of $C$, that is, the homology of $C^*$,
is of finite type and that, likewise, the homology of $A^*$ is of 
finite type as well.
These hypotheses are not independent, see Remark 3.6 below.
Consider the functors
$$
\align
t &\colon {}_A\roman{Mod} @>>> {}_C\roman{Comod},
\quad
t(N) = C \otimes_{\tau} N
\tag3.1${}_*$
\\
h &\colon  {}_C\roman{Comod} @>>> {}_A\roman{Mod},
\quad h(M) =  A \otimes_{\tau}M.
\tag3.2${}_*$
\endalign
$$
Unfortunately the notation $h$ conflicts with the earlier notation
for a contracting homotopy but, for intelligibility, we prefer
to keep the notation $h$ for both objects;
in the sequel, which one is intended will always be clear from
the context.

Given the left $C$-comodule $M$,
the
composite 
$$
M @>>>  C \otimes_{\tau} A \otimes_{\tau} M = t(h(M))
\tag3.3${}_*$
$$
of the comodule structure map
$M @>>>  C \otimes M$ with the
canonical injection
$$C \otimes M \to  C \otimes A \otimes M$$
is a morphism of left $C$-comodules that is natural in $M$.
Henceforth the numbering  (3.3${}_*$)
will also refer to the corresponding natural transformation
$\Cal I \to th$ where $\Cal I$ refers to the identity functor.

Likewise, given the left  $A$-module $N$,
the composite 
$$
h(t(N)) = A \otimes_{\tau} C \otimes_{\tau} N @>>> N
\tag3.4${}_*$
$$
of the canonical
projection $A \otimes C \otimes N \to A\otimes N$ with the $A$-module
structure map $A\otimes N \to N$ of $N$ 
is a morphism of left $A$-modules that is natural in $N$.
Henceforth the numbering  (3.4${}_*$)
will also refer to the corresponding natural transformation
$ht \to \Cal I$.

The functors $t$ and $h$ are chain homotopy inverse to each other
in a very precise sense which we now explain.

Let $\Cal F \colon {}_C\roman{Comod} \to {}_C\roman{Comod}$
and 
$\Cal G \colon {}_A\roman{Mod} \to {}_A\roman{Mod}$
be functors, and denote the identity functor
by $\Cal I$. We define a {\it left\/} $C$-{\it comodule
functor contraction of \/} $\Cal F$ {\it to the identity\/} to be a contraction
$$
\Nsddata{\Cal F}{\nabla}{\pi}{\Cal I}{s};
$$
here $\nabla$, $\pi$, $s$ are natural transformations in such a way that,
for every left $C$-comodule $M$, the data
$$
\Nsddata{\Cal F(M)}{\nabla_M}{\pi_M}{M}{s_M}
$$
constitute a contraction of chain complexes with
$\nabla_M$ a morphism of $C$-comodules.
Likewise we define a {\it left\/} $A$-{\it module
functor contraction of \/} $\Cal G$ {\it to the identity\/} to be a contraction
$$
\Nsddata{\Cal G}{\nabla}{\pi}{\Cal I}{s};
$$
here $\nabla$, $\pi$, $s$ are natural transformations in such a way that,
for every left $A$-module $N$, the data
$$
\Nsddata{\Cal G(N)}{\nabla_N}{\pi_N}{N}{s_N}
$$
constitute a contraction of chain complexes with
$\pi_N$ a morphism of $A$-modules.
Accordingly we define the notion of
{\it right\/} $C$-{\it comodule
functor contraction to the identity\/}
and that of
{\it right\/} $A$-{\it module
functor contraction to the identity\/}.

\proclaim{Theorem 3.5${}_*$}
The functors $t$ and $h$ are chain homotopy inverse to each other
in the sense of {\rm (1)} and {\rm (2)} below:
\roster
\item
The natural transformation
{\rm (3.3$_*$)} of left $C$-comodules
extends to a left $C$-comodule contraction to the identity
$$
\Nsddata{th}{ (3.3_*)}{\phantom{91}\pi\phantom{91} }{\Cal I}{s}
$$
of the endofunctor $th$ on ${}_C\roman{Comod}$; 
\item
the natural transformation {\rm (3.4$_*$)} 
of left $A$-modules
extends to a left $A$-module 
contraction to the identity
$$
\Nsddata{ht}{\phantom{91} \nabla \phantom{91}} { (3.4_*)}{\Cal I}{s}
$$
of the endofunctor $ht$ on ${}_A\roman{Mod}$.
\endroster
\endproclaim

\demo{Proof} 
Let $M$ be a left $C$-comodule. The projection
$$
\varepsilon_M =
\varepsilon \otimes \varepsilon \otimes \roman{Id}
\colon t(h(M))=  C \otimes_{\tau} A \otimes_{\tau} M @>>> M
\tag3.5.1$_*$
$$
is a chain map (beware: not a morphism of $C$-comodules) that is 
a retraction for (3.3$_*$). Furthermore, 
$\varepsilon_M$ and (3.3$_*$) are also chain maps when
$C \otimes_{\tau} A \otimes M$
is substituted for $C \otimes_{\tau} A \otimes_{\tau} M$, and 
the morphism
$$
s^0_M=s \otimes \roman{Id}
\colon  C \otimes_{\tau} A \otimes M 
\longrightarrow
C \otimes_{\tau} A \otimes M 
\tag3.5.2$_*$
$$
yields a chain homotopy between 
the identity and the composite $(\roman{3.3}_*)\,\varepsilon_M$.
Thus the data
$$
\Nsddata{C \otimes_{\tau} A \otimes M}
{(3.3_*)}
{\phantom{1}\varepsilon_M \phantom{91}}
M
{s^0_M}
\tag3.5.3$_*$
$$
constitute a contraction.
If the side conditions (2.1.3) are not satisfied we can modify $s^0_M$ if need be,
cf. Remark 2.1.4, and we suppose that this has already been arranged for.
We do not indicate this in notation. By construction,
$$
\varepsilon_M s^0_M=0,
\ 
s^0_M (3.3_*)=0.
\tag3.5.4$_*$
$$
Write the differential on 
$C \otimes_{\tau} A \otimes_{\tau} M
$ 
as $d+ \partial^{\tau}$
where $d$ refers to the differential on
$C \otimes_{\tau} A \otimes M$.
Relative to the filtration
of $C \otimes_{\tau} A \otimes_{\tau} M$
induced by the skeletal filtration of $M$,
the operator $\partial^{\tau}$ lowers filtration,
and this filtration is complete.
Application of the perturbation lemma (Lemma 2.3)
yields a new contraction
$$
\Nsddata{C \otimes_{\tau} A \otimes_{\tau} M}{(3.3_*)}
{\phantom{1}\varepsilon_M \phantom{91}}
M
{s_M}.
\tag3.5.5$_*$
$$
In view of (3.5.4$_*$), the perturbation modifies only the homotopy $s^0_M$.
The morphisms $(3.3_*)$, $\varepsilon_M$ and $s_M$ are plainly natural in
$M$.

Likewise, let $N$ be a left $A$-module. The injection
$$
\eta_N =
\eta \otimes \eta \otimes \roman{Id}\colon
N @>>>  A \otimes_{\tau} C \otimes_{\tau} N = h(t(N)) 
\tag3.5.6$_*$
$$
is a chain map (beware: not a morphism of $A$-modules) that is 
a section for (3.4$_*$). 
Furthermore, 
$\eta_N$ and (3.4$_*$) are also chain maps when
$A \otimes_{\tau} C \otimes N$
is substituted for $A \otimes_{\tau} C \otimes_{\tau} N$, and 
the morphism
$$
s^0_N=s \otimes \roman{Id}
\colon  A \otimes_{\tau} C \otimes N 
\longrightarrow
A \otimes_{\tau} C \otimes N 
\tag3.5.7$_*$
$$
yields a chain homotopy between 
the identity and the composite $\eta_N(\roman{3.4}_*)$.
Thus the data
$$
\Nsddata{A \otimes_{\tau} C \otimes N}
{\phantom{91}\eta_N \phantom{1}}
{(3.4_*)}N
{s^0_N}
\tag3.5.8$_*$
$$
constitute a contraction.
Again, if the side conditions (2.1.3) 
are not satisfied we can modify $s^0_N$ if need be,
cf. Remark 2.1.4, and we suppose that this has already been arranged for.
We do not indicate this in notation. By construction,
$$
s^0_N \eta_N=0,
\ 
(3.3_*) s^0_N=0.
\tag3.5.9$_*$
$$
Write the differential on 
$A \otimes_{\tau} C \otimes_{\tau} N
$ 
as $d+ \partial^{\tau}$
where $d$ refers to the differential on
$A \otimes_{\tau} C \otimes N$.
Relative to the filtration of
$A \otimes_{\tau} C \otimes_{\tau} N$
coming from the skeletal filtration of
$A \otimes_{\tau} C$,
the operator $\partial^{\tau}$ lowers filtration, and the 
filtration is complete. 
Application of the perturbation lemma (Lemma 2.3)
yields a new contraction
$$
\Nsddata{A \otimes_{\tau} C \otimes_{\tau} N}
{\phantom{91}\eta_N\phantom{1} }
{(3.4_*)}N
{s_N}
\tag3.5.10$_*$
$$
In view of (3.5.9$_*$), the perturbation modifies only the homotopy $s^0_N$.
The morphisms $(3.4_*)$, $\eta_N$ and $s_N$ are plainly natural in
$N$. \qed\enddemo

The theorem implies the following:
Suppose that $C$ and $A$ satisfy certain mild arithmetical hypotheses
such as e.~g. that they are projective over the ground ring as graded modules.
Then, when $C$ is 
connected, 
the twisted object 
$t(h(M))$ serves as a replacement for a relatively injective resolution
of $M$ in the category of left $C$-comodules,
that is, given the right $C$-comodule $M'$,
the homology of the cotensor product
$$
 M'\square_C t(h(M)) = M'\square_C  (C \otimes_{\tau} A \otimes_{\tau} M)
\cong  M' \otimes_{\tau} A \otimes_{\tau} M
$$
yields the differential graded 
$\roman{Cotor}^{C}(M',M)$.
When $C$ is, furthermore, cocomplete,
the twisted object $h(t(N))$ serves as a replacement for
a relatively projective resolution 
of $N$ in the category of left $A$-modules, that is,
given the right $A$-comodule $N'$,
the homology of the tensor product
$$
N'\otimes_A h(t(N)) = N'\otimes_A (A \otimes_{\tau} C \otimes_{\tau} N)
\cong
 N'\otimes_{\tau} C \otimes_{\tau} N
$$
yields the differential graded
$\roman{Tor}_{A}(N',N)$.
In particular, the functor 
$t$ assigns to a left $A$-module $N$ a twisted tensor product
which calculates the
differential graded $\roman{Tor}^{A}(R,N)$,
and the functor $h$
assigns to a left $C$-comodule $M$ a 
twisted tensor product
which,
when $C$ is 
connected, 
calculates the differential
graded $\roman{Cotor}^{C}(R,M)$.
The reason for the connectedness hypothesis is that, in the definition of
Cotor, the injective resolution must be assembled by the product
and not by the coproduct that is used in the definition of both
$\rbar$ and $\rcob$.

\smallskip
The dual $C^*$ of $C$ is an augmented differential graded algebra.
We remind the reader that  ${}_{C^*}\widehat{\roman{Mod}}$
refers to the category of $C^*$-modules which are complete
in a sense explained in (1.2.4) above.
Relative to the differential graded algebra $C^*$,
the functor
$$
t^* \colon
\roman{Mod}_A @>>> {}_{C^*}\widehat{\roman{Mod}},
\quad
t^*(\Nfl) = \roman{Hom}^{\tau}(C,\Nfl)
\tag3.1${}^*$
$$
assigns to a right $A$-module $\Nfl$ the twisted Hom-object
$t^*(\Nfl)$ which calculates the differential graded
$\roman{Ext}_{A}(R,\Nfl)$, and the functor
$$
h^* \colon
{}_{C^*}\widehat{\roman{Mod}} @>>> \roman{Mod}_A,
\quad
h^*(\Mfl) = \roman{Hom}^{\tau}(A,\Mfl)
\tag3.2${}^*$
$$
assigns to a (complete)
left $C^*$-module $\Mfl$ the twisted Hom-object
$h^*(\Mfl)$ which, since the homology of $C^*$  is of finite type,
calculates the differential graded $\roman{Tor}^{C^*}(R,\Mfl)$.
When $A$ is itself of finite type, 
$A^*$ is a differential graded coalgebra, the dual
$\tau^* \colon A^* \to C^*$ 
(given by $\tau^*(\alpha)=(-1)^{|\alpha|}\alpha \circ \tau$)
is a twisting cochain, and the obvious morphism is
an isomorphism
$$
A^* \otimes_{\tau^*} \Mfl @>>> \roman{Hom}^{\tau}(A,\Mfl)
$$
of (differential graded) left $A^*$-comodules. Whether or not $A$
is of finite type, 
given the complete
left $C^*$-module $\Mfl$,
the composite of the structure map
$\roman{Hom}(C, \Mfl) \to  \Mfl$ with the morphism
$$
\roman{Hom}(C, \roman{Hom}(A,\Mfl)) @>>> \roman{Hom}(C,
\Mfl)
$$
induced by the evaluation map
$$
\roman{Hom}(A,\Mfl) @>>> \Mfl,\quad \phi \mapsto \phi(1),
$$
yields a morphism
$$
 \Mfl @<<< \roman{Hom}^{\tau}(C \otimes_{\tau} A,\Mfl)
\cong \roman{Hom}^{\tau}(C, \roman{Hom}^{\tau}(A,\Mfl))
=t^*(h^*(\Mfl)) \tag3.3${}^*$
$$
of left $C^*$-modules which is also a chain equivalence that is
natural in $\Mfl$,
and 
$t^*(h^*(\Mfl))$ is
relatively projective when $C$ is of finite type. Likewise, 
let $\Nfl$ be a right $A$-module.
The
composite of the adjoint $\Nfl \to \roman{Hom}(A,\Nfl)$
of the structure map $\Nfl \otimes A \to \Nfl$ with the
canonical injection
$$
\varepsilon^*\colon \roman{Hom}(A,\Nfl) @>>>
\roman{Hom}(A,\roman{Hom}(C,\Nfl))
$$
induced by the counit $\varepsilon$ of $C$ 
yields a morphism
$$
\Nfl @>>> \roman{Hom}^{\tau}(A \otimes_{\tau} C,\Nfl)
\cong \roman{Hom}^{\tau}(A,\roman{Hom}^{\tau}(C,\Nfl))
=h^*(t^*(\Nfl))
\tag3.4${}^*$
$$
in the category of right $A$-modules
which is also a chain equivalence that is natural in $\Nfl$.

We define a {\it complete left\/} $C^*$-{\it module
functor contraction to the identity\/} of the functor 
$\Cal G\colon {}_{C^*}\widehat{\roman{Mod}} \to
{}_{C^*}\widehat{\roman{Mod}}$  to be a contraction
$$
\Nsddata{\Cal G}{\nabla}{\pi}{\Cal I}{s}
$$
as defined before for $A$-modules save that
now the constituents of the contraction refer
to complete $C^*$-modules.
Henceforth the numberings  
(3.3${}^*$) and
(3.4${}^*$)
will also refer to the corresponding natural transformations.

\proclaim{Theorem 3.5${}^*$}
The functors $t^*$ and $h^*$ are chain homotopy inverse to each other
in the sense of {\rm (1)} and {\rm (2)} below:
\roster
\item
The natural transformation
{\rm (3.3$^*$)} of complete
left  $C^*$-modules
extends to a complete left $C^*$-module contraction
to the identity
$$
\Nsddata{t^*h^*}{\phantom{91}\nabla\phantom{91} }{ (3.3^*)}{\Cal I}{s}
$$
of the endofunctor $t^*h^*$ on ${}_{C^*}\widehat{\roman{Mod}}$; 
\item
the natural transformation {\rm (3.4$^*$)} of $A$-modules
extends to a 
contraction to the identity
$$
\Nsddata{h^*t^*} { (3.4^*)}{\phantom{91} \pi \phantom{91}}{\Cal I}{s}
$$
of the endofunctor $h^*t^*$ on $\roman{Mod}_A$.
\endroster
\endproclaim

\demo{Proof}
Let $\Mfl$ be a complete left $C^*$-module. The obvious injection
$$
\eta_{\Mfl}\colon \Mfl @>>> 
\roman{Hom}^{\tau}(C \otimes_{\tau} A,\Mfl)
\cong \roman{Hom}^{\tau}(C, \roman{Hom}^{\tau}(A,\Mfl))
=t^*(h^*(\Mfl))
\tag3.5.1${}^*$
$$
is a chain map (beware: not a morphism of $C^*$-modules) that is 
a section for (3.3$^*$). 
Furthermore, 
$\eta_{\Mfl}$ and (3.3$^*$) are also chain maps when
$\roman{Hom}(C \otimes_{\tau} A,\Mfl)$
is substituted for 
$\roman{Hom}^{\tau}(C \otimes_{\tau} A,\Mfl)$, 
and 
the morphism
$$
s^0_{\Mfl}\colon  
\roman{Hom}(C \otimes_{\tau} A,\Mfl)
\longrightarrow
\roman{Hom}(C \otimes_{\tau} A,\Mfl),
\ s^0_{\Mfl}(\beta)=(-1)^{|\beta|}\beta \circ s
\tag3.5.2${}^*$
$$
yields a chain homotopy between 
the identity and the composite $\eta_{\Mfl}(\roman{3.3}^*)$.
Thus the data
$$
\Nsddata{\roman{Hom}(C \otimes_{\tau} A,\Mfl)}
{\phantom{1}\eta_{\Mfl}\phantom{4}}
{(3.3^*)}{\Mfl}
{s^0_{\Mfl}}
\tag3.5.3${}^*$
$$
constitute a contraction.
If the side conditions (2.1.3) are not satisfied we can modify $s^0_{\Mfl}$ 
if need be,
cf. Remark 2.1.4, and we suppose that this has already been arranged for.
We do not indicate this in notation. By construction,
$$
s^0_{\Mfl} \eta_{\Mfl}=0,
\ 
(3.3^*)s^0_{\Mfl} =0.
\tag3.5.4$^*$
$$
Write the differential on 
$\roman{Hom}^{\tau}(C \otimes_{\tau} A,\Mfl)$ 
as $d+ \delta^{\tau}$
where $d$ refers to the differential on
$\roman{Hom}(C \otimes_{\tau} A,\Mfl)$. 
Relative to the filtration of
$\roman{Hom}^{\tau}(C \otimes_{\tau} A,\Mfl)$
induced by the skeletal filtration of
$\Mfl$, 
the operator $\delta^{\tau}$ lowers filtration,
and the filtration of
$\roman{Hom}^{\tau}(C \otimes_{\tau} A,\Mfl)$
is complete.
Application of the perturbation lemma (Lemma 2.3)
yields a new contraction
$$
\Nsddata{\roman{Hom}^{\tau}(C \otimes_{\tau} A,\Mfl)}
{\phantom{1}\eta_{\Mfl}\phantom{4}}
{(3.3^*)}{\Mfl}
{s_{\Mfl}}
\tag3.5.5${}^*$
$$
In view of (3.5.4$^*$), the perturbation modifies only the homotopy 
$s^0_{\Mfl}$.
The morphisms $(3.3^*)$, 
$\eta_{\Mfl}$ and $s_{\Mfl}$
are plainly natural in
$M$.

Likewise, let $\Nfl$ be a right $A$-module. The obvious projection
$$
\varepsilon_{\Nfl} 
\colon
h^*(t^*(\Nfl))= \roman{Hom}^{\tau}(A,\roman{Hom}^{\tau}(C,\Nfl))
@>>> \Nfl 
\tag3.5.6${}^*$
$$
is a chain map (beware: not a morphism of $A$-modules) that is 
a retraction for (3.4$^*$). 
Furthermore, 
$\varepsilon_N$ and (3.4${}^*$) are also chain maps when
$\roman{Hom}^{\tau}(A,\roman{Hom}(C,\Nfl))$
is substituted for 
$\roman{Hom}^{\tau}(A,\roman{Hom}^{\tau}(C,\Nfl))$.
By adjointness,
$$
 \roman{Hom}^{\tau}(A,\roman{Hom}(C,\Nfl))
\cong \roman{Hom}(A \otimes_{\tau} C,\Nfl)
$$
and 
the morphism
$$
s^0_N
\colon  \roman{Hom}(A \otimes_{\tau} C,\Nfl) 
\longrightarrow
\roman{Hom}(A \otimes_{\tau} C,\Nfl),
\ s^0_N(\beta)=(-1)^{|\beta|}\beta \circ s
\tag3.5.7${}^*$
$$
yields a chain homotopy between 
the identity and the composite $(\roman{3.4}^*)\varepsilon_N$.
Thus the data
$$
\Nsddata{\roman{Hom}(A \otimes_{\tau} C,\Nfl)}
{\phantom{91}\varepsilon_N \phantom{4}}
{(3.4^*)}N
{s^0_N}
$$
constitute a contraction and, for the sake of clarity of exposition,
we write this contraction as
$$
\Nsddata{ \roman{Hom}^{\tau}(A,\roman{Hom}(C,\Nfl))}
{\phantom{91}\varepsilon_N \phantom{4}}
{(3.4^*)}N
{s^0_N} ,
\tag3.5.8${}^*$
$$
with a slight abuse of the notation $s^0_N$.
Again, if the side conditions (2.1.3) 
are not satisfied we can modify $s^0_N$ if need be,
cf. Remark 2.1.4, and we suppose that this has already been arranged for.
We do not indicate this in notation. By construction,
$$
 \varepsilon_Ns^0_N=0,
\ 
s^0_N(3.3^*)=0.
\tag3.5.9${}^*$
$$
Write the differential on 
$ \roman{Hom}^{\tau}(A,\roman{Hom}^{\tau}(C,\Nfl))$ 
as $d+ \delta^{\tau}$
where $d$ refers to the differential on
$\roman{Hom}^{\tau}(A,\roman{Hom}(C,\Nfl))$.
Relative to the filtration
of 
$$
 \roman{Hom}(A,\roman{Hom}(C,\Nfl))
\cong  \roman{Hom}(A \otimes C,\Nfl)
$$
coming from the skeletal filtration of
$A \otimes C$,
the operator $\delta^{\tau}$ lowers filtration,
and the filtration of  
$\roman{Hom}(A,\roman{Hom}(C,\Nfl))$
is complete.
Application of the perturbation lemma (Lemma 2.3)
yields a new contraction
$$
\Nsddata{ 
 \roman{Hom}^{\tau}(A,\roman{Hom}^{\tau}(C,\Nfl))
}
{\phantom{91}\varepsilon_N\phantom{1} }
{(3.4^*)}N
{s_N}
\tag3.5.10${}^*$
$$
In view of (3.5.9${}^*$), the perturbation modifies only the homotopy $s^0_N$.
The morphisms $\eta_N$ and $s_N$ are plainly natural in
$N$. \qed\enddemo

\smallskip
\noindent {\smc Remark 3.6\/.} Suppose that $C$
is cocomplete. Then the adjoint $\overline \tau \colon C \to \rbar A$
is a morphism of coalgebras and a chain equivalence. Indeed,
the canonical comparison yields a chain inverse
$$
\rbar A \otimes_{\tau^{\roman B}} A \longrightarrow C \otimes_{\tau}A
$$
which is in fact a morphism in the category of differential graded
$A$-modules. This morphism descends to a chain map
$\rbar A \to C$, not necessarily compatible
with the coalgebra structures but a chain homotopy inverse
for $\overline \tau$.
Consequently the induced map
$$
\rcob C \longrightarrow \rcob \rbar A \longrightarrow A
$$
is a morphism of differential graded algebras  and  a chain equivalence
whence the induced chain map
$$
A^* \longrightarrow (\rcob C)^*=\roman{Hom}(\rcob C,R)
$$
is a chain equivalence.
Furthermore, since $\roman H^*C$ is supposed to be of finite type,
the canonical map
$$
\rbar C^* \longrightarrow \roman{Hom}(\rcob C,R)
$$
is a chain equivalence.
Since the homology of $\rbar C^*$ is the differential torsion product
$\roman{Tor}_{C^*}(R,R)$, we conclude that the homology of $A^*$, that is,
the cohomology $\roman H^*(A)$ of $A$, is naturally isomorphic to
$\roman{Tor}_{C^*}(R,R)$.
Since $\roman H^*(C)$ is of finite type, so is
$\roman{Tor}_{C^*}(R,R)$ whence 
$\roman H^*(A)$ is necessarily of finite type.

\smallskip
\noindent
{\smc Remark 3.7.\/} Suppose that $C$ is itself of finite type,
and let $\Mfl$ be a left $C^*$-module.
By adjointness, the left $C^*$-module structure is then induced from a
left $C$-comodule structure on $\Mfl$,
and the twisted Hom-object $h^*(\Mfl)$
calculates the differential graded
$\roman{Coext}_C(R,\Mfl)$.

\smallskip
Let $N$ be a differential graded left $A$-module, and consider
$N^*$, viewed as a differential graded right $A$-module.
The twisted Hom-object 
$$
\roman{Hom}^{\tau}(C,N^*) \cong \roman{Hom}_A(C\otimes_{\tau}A,N^*)
$$
acquires a canonical $C^*$-module structure and, furthermore,
calculates the differential $\roman{Ext}_A(R,N^*)$.
The bar construction $\roman B(R,C^*,\roman{Hom}^{\tau}(C,N^*))$ can be written
as the twisted object
$$
\roman B(R,C^*,\roman{Hom}^{\tau}(C,N^*)) 
=
(\rbar C^*)\otimes_{\tau^{\rbar}} \roman{Hom}^{\tau}(C,N^*) 
$$
and the forgetful map $\roman{Hom}^{\tau}(C,N^*) \to N^*$
is a chain map and extends canonically to a projection
$$
\pi_{N^*} \colon
\roman B(R,C^*,\roman{Hom}^{\tau}(C,N^*)) 
\longrightarrow N^* .
\tag3.8
$$
For later reference, we will now explore the question under what circumstances
this projection is a chain equivalence.

The projection (3.8) is formally similar to the projection
(3.5.1$_*$) above.
Furthermore, it is a kind of formal dual of a map 
$$
N @>>> \rcob C \otimes_{\tau_{\rcob}} C \otimes_{\tau}N
\tag3.9
$$
of the kind (3.5.6$_*$) above. 
By Theorem 3.5$_*$, the map (3.5.4), in turn, fits into an 
$(\rcob C)$-module contraction of the kind
$$
\Nsddata N {(3.4_*)}{(3.5.4)}
{\rcob C \otimes_{\tau_{\rcob}} C \otimes_{\tau}N}{h_N}
\tag3.10
$$
where, relative to (3.4$_*$), we have substituted
$\rcob C$ for $A$.
Dualizing this contraction, we arrive at the contraction
$$
\Nsddata {N^*}{(3.5.4)^*} {(3.4_*)^*}
{\roman{Hom}(\rcob C \otimes_{\tau_{\rcob}} C \otimes_{\tau}N,R)}{h^*_N}
\tag3.11
$$
having the feature that $(3.4_*)^*$ is compatible with the additional structure.
The object $\roman{Hom}(\rcob C \otimes_{\tau_{\rcob}} C \otimes_{\tau}N,R)$
can be written as the complete twisted Hom-object
$$
\roman{Hom}^{\tau_{\rcob}}(\rcob C,\roman{Hom}^{\tau}(C,N^*))
$$

Consider the canonical injection
$$
\roman B(R,C^*,\roman{Hom}^{\tau}(C,N^*)) 
@>>>
\roman{Hom}(\rcob C \otimes _{\tau_{\rcob}} C \otimes_{\tau}N,R),
\tag3.12
$$
or, equivalently,
$$
\rbar C^* \otimes_{\tau^{\rbar}}\roman{Hom}^{\tau}(C,N^*) 
@>>>
\roman{Hom}^{\tau_{\rcob}}(\rcob C,\roman{Hom}^{\tau}(C,N^*)).
\tag3.13
$$
The composite of (3.12)
with the canonical forgetful map
$$
\roman{Hom}(\rcob C \otimes _{\tau_{\rcob}} C \otimes_{\tau}N,R)
@>>>
N^*
$$
coincides with the projection (3.8).

\proclaim{Lemma 3.14}
Suppose that $C$ is non-negative, simply connected and free as a graded
$R$-module.
When the homology $\roman H(C)$ of $C$ is of finite type,
the injection {\rm (3.13)} is a chain equivalence 
whence the projection {\rm (3.8)} is then a chain equivalence.
\endproclaim

\demo{Proof} We will exploit the construction in \cite\huebkade\ (4.6):
To adjust to the notation in the quoted reference,
let $V=s^{-1}JC$ and $\roman H=\roman H(s^{-1}JC)$,
the differential on $s^{-1}JC$ being that induced from the
differential on $C$. The differential graded algebra $\rcob C$
is one of the kind $(\roman T[V],d)$, the differential $d$ being 
the ordinary coalgebra perturbation of the differential on $V=s^{-1}JC$.
As in \cite\huebkade\ (4.6), let $(F\roman H, \delta)$ be a free resolution
of $\roman H$ which we take here to be of finite type.
The construction in \cite\huebkade\ (4.6) yields 
a perturbed differential $d'$ on the graded tensor algebra
$\roman T[F\roman H]$ (perturbing the ordinary tensor algebra differential) and
a filtered chain equivalence
$$
\Nsddata {(\roman T[V],d)}
{f}{g}{\mu',(\roman T[F\roman H],d')}
{\nu'}
\tag3.14.1
$$
of augmented differential graded algebras. 
By construction, 
the augmented differential graded algebra $\Cal A=(\roman T[F\roman H],d')$
is of finite type and, since (3.14.1) is a filtered chain equivalence
of augmented differential graded algebras,
$f\colon \Cal A \to \rcob C$ and $g \colon \rcob C \to A$
are morphisms of differential graded algebras which are also chain
equivalences,
indeed, the chain homotopies $\mu'$ and $\nu'$ are even compatible with
the algebra structures.

Since $\Cal A$ is of finite type, the dual $\Cal A^*$ is a coaugmented
differential graded coalgebra.
Moreover, the composite 
$$
\tau_g\colon C @>{\tau_{\rcob}}>> \rcob C @>g>> \Cal A
$$
is a twisting cochain, and so is the dual $\tau^*_g\colon \Cal A^* \to C^*$.

Instead of (3.13), we can now consider
the canonical chain map
$$
\Cal A^* \otimes_{\tau^*_g}\roman{Hom}^{\tau}(C,N^*) 
@>>>
\roman{Hom}^{\tau_g}(\Cal A,\roman{Hom}^{\tau}(C,N^*))\cong
\roman{Hom}(\Cal A \otimes_{\tau_g} C \otimes_{\tau} N,R))
\tag3.14.2
$$
Since $\Cal A$ is of finite type, this chain map is an isomorphism.
The left-hand side of (3.14.2) is chain homotopic to
$\rbar C^* \otimes_{\tau^{\rbar}}\roman{Hom}^{\tau}(C,N^*)$
compatibly with the bundle structure and the right-hand side 
of (3.14.2)
is chain homotopic to
$\roman{Hom}^{\tau_{\rcob}}(\rcob C,\roman{Hom}^{\tau}(C,N^*))$
compatibly with the bundle structures. Indeed, the adjoint
$\overline \tau^*_g\colon \Cal A^* \to \rbar C^*$ of
$\tau^*_g$ yields the chain map
$$
\Cal A^* \otimes_{\tau^*_g}\roman{Hom}^{\tau}(C,N^*) 
\longrightarrow
\rbar C^* \otimes_{\tau^{\rbar}}\roman{Hom}^{\tau}(C,N^*),
$$
necessarily compatible with the bundle structures and a chain equivalence.
Likewise the morphism
$g \colon \rcob C \to \Cal A$
of augmented differential graded algebras
induces the morphism
$$
g \otimes C \otimes N\colon\rcob C \otimes_{\tau_{\rcob}} C \otimes_{\tau} N
\longrightarrow
\Cal A \otimes_{\tau_g} C \otimes_{\tau} N
$$
of twisted objects, necessarily a chain equivalence.
This morphism of bundles, in turn, dualizes to the chain equivalence
$$
\roman{Hom}(\rcob C \otimes_{\tau_{\rcob}} C \otimes_{\tau} N,R)
\longleftarrow
\roman{Hom}(\Cal A \otimes_{\tau_g} C \otimes_{\tau} N,R).
$$
Putting the various items together we obtain the
asserted chain equivalence. \qed
\enddemo

\medskip\noindent{\bf 4. Small models for ordinary equivariant (co)homology}
\smallskip\noindent
Given the two simplicial sets $S$ and $T$, the Eilenberg-Zilber theorem
(cf. \cite\eilmactw\ II, Theorem 2.1) yields a contraction
$\Nsddata {S\times T} \nabla g {S\otimes T }h$;
here we do not distinguish in notation between a simplicial set
and its normalized chain complex.
See e.~g. \cite\eilmactw\ (I \S 5 and II \S 2), 
\cite\gugenhtw\ (section 4, in particular p.~409),
\cite\gugenmay\ (Proposition A.3) for details.
The morphism $\nabla$ is referred to as the {\it shuffle map\/}
and $g$  as the {\it Alexander-Whitney map\/}.
In \cite\eilmactw\ (II, Theorem 2.1) the vanishing of the square
$hh$ of the homotopy $h$ constructed there is not established
but, in view of Remark 2.1.4 above, this vanishing can always be achieved.
The fact that the shuffle map is a morphism of differential
graded coalgebras is given in \cite\eilenmoo\ (\S 17),
cf. also Lemma 4.4 in \cite\gugenmay. A somewhat more general form of the
Eilenberg-Zilber theorem, due to Cartier, can be found in 
\cite\doldpupp\ (2.9). Below we will apply the 
Eilenberg-Zilber theorem 
to the singular set $SY$ of a topological space $Y$ and to the simplicial 
bar construction. The Eilenberg-Zilber theorem 
is well known to imply that the normalized chain complex $C_*S$ 
of a simplicial set
$S$ acquires a differential graded coalgebra structure; 
moreover, a choice of base point determines a coaugmentation
for $C_*S$, 
and $C_*S$ has a unique coaugmentation and is cocomplete
if and only if $S$ is {\it reduced\/} in the sense that $S_0$ consists of 
a single point.
Now, given a pathwise and locally pathwise connected topological space $Y$,
the injection of the 
first Eilenberg subcomplex $S^1Y$ relative to a chosen base point $o$
(which consists of all singular simplices
mapping each vertex to the base point $o$)
into the singular set $SY$ associated to $Y$ is an injection of differential
graded coalgebras and a chain equivalence;
thus we can then take $C_*(Y)$ to be a cocomplete differential 
graded coalgebra.

Let $G$ be a topological group.
Before we proceed further we note the following.
Suppose that the group $G$ is connected.
Then we can replace the
singular set $SG$ associated to $G$---this is a simplicial group---with
the first Eilenberg subcomplex $S^1G$ 
relative to the neutral element of $G$; 
the first Eilenberg subcomplex $S^1G$ is a simplicial subgroup
and the injection of $S^1G$ into $SG$ is a morphism of simplicial groups 
which is a chain equivalence. The normalized chain algebra of $S^1G$
is connected. Thus, given a connected topological group $G$,
we can take $C_*G$ to be this normalized chain algebra.

Let now $G$ be a general topological group,
$X$  a left  $G$-space and $E$ a
right $G$-space. It is well known that the Eilenberg-Zilber theorem
has the following consequences:
(i)~{\sl~The chain complex $C_*G$ of
normalized singular chains on $G$ acquires a
differential graded Hopf algebra structure\/};
(ii) {\sl this Hopf algebra $C_*G$ acts from the left on the
normalized singular chains $C_*(X)$ and from the right on the
normalized singular cochains $C^*(X)$ on $X$\/}. 

The construction of the small models
for ordinary equivariant (co)homology
relies on the following
folk-lore result, cf. \cite\gugenmay\ (Theorem 3.9),
the significance of which was already commented on
in the introduction.

\proclaim{Theorem 4.1}
The $G$-equivariant homology 
$\roman H^G_*(X,R)$
of $X$ is canonically isomorphic to 
the differential
$\roman{Tor}_*^{C_*G}(R,C_*X)$ and, likewise,
the $G$-equivariant cohomology $\roman H_G^*(X,R)$
of $X$ is canonically isomorphic to
the differential $\roman{Ext}^*_{C_*G}(R,C^*X)$.
The isomorphisms are natural in the data.
\endproclaim

Consider the $(G\times G)$-action on $X$ which is the given action on 
the first copy of $G$ and the trivial action on the second one.
The diagonal map
$G \to G \times G$ induces the morphisms
$$
\align
\roman H^G_*(X,R)&\longrightarrow \roman H^{G\times G}_*(X,R)
\tag4.1.1
\\
\roman H_G^*(X,R) \otimes \roman H^*(BG,R)
&\longrightarrow
\roman H_G^*(X,R)
\tag4.1.2
\endalign
$$
of graded $R$-modules that satisfy the standard associativity constraints
so that, in particular,
$\roman H_G^*(X,R)$ acquires a graded $\roman H^*(BG,R)$-module structure
and that, when 
$\roman H_*(BG,R)$ is projective as a graded $R$-module,
$\roman H_*(BG,R)$ acquires a graded cocommutative
coalgebra structure and
$\roman H_G^*(X,R)$  a graded comodule structure over $\roman H_*(BG,R)$.
Likewise, the diagonal map
of $G$ induces the morphisms
$$
\align
\roman{Tor}_*^{C_*G}(R,C_*X)&\longrightarrow \roman{Tor}_*^{C_*(G\times G)}(R,C_*X)
\tag4.1.3
\\
\roman{Ext}^*_{C_*G}(R,C^*X)\otimes \roman{Ext}^*_{C_*G}(R,R)
&\longrightarrow
\roman{Ext}^*_{C_*G}(R,C^*X)
\tag4.1.4
\endalign
$$
of graded $R$-modules that satisfy the standard associativity constraints
so that, in particular,  $\roman{Ext}^*_{C_*G}(R,R)\cong \roman H^*(BG)$
is a graded commutative algebra  
and that $\roman{Ext}^*_{C_*G}(R,C^*X)$
acquires a graded  $\roman{Ext}^*_{C_*G}(R,R)$-module structure
via (4.1.4).
Moreover, when 
$\roman{Tor}_*^{C_*G}(R,R)\cong \roman H_*(BG)$
is projective as a graded $R$-module,
$\roman{Tor}_*^{C_*G}(R,R)$ acquires a graded cocommutative
coalgebra structure and
$\roman{Tor}_*^{C_*G}(R,C_*X)$ acquires a
graded comodule structure over the graded coalgebra
$\roman{Tor}_*^{C_*G}(R,R)$ via (4.1.3).

\proclaim{Addendum 1}
The naturality of the constructions implies that the
isomorphisms
between
$\roman H^G_*(X,R)$
and $\roman{Tor}_*^{C_*G}(R,C_*X)$ and
between  $\roman H_G^*(X,R)$
and $\roman{Ext}^*_{C_*G}(R,C^*X)$
are compatible with the additional structure just explained.
\endproclaim

\proclaim{Addendum 2}
More generally, when $E$ is a principal $G$-space,
there are canonical isomorphims of the kind
$$
\roman H_*(E \times _GX,R)
\cong \roman{Tor}_*^{C_*G}(C_*E,C_*X),
\quad
\roman H^*(E \times _GX,R)
\cong\roman{Ext}^*_{C_*G}(C_*E,C^*X) 
$$
which are natural in the data.
\endproclaim

We now explain briefly how Theorem 4.1 is established by means of
certain HPT-techniques needed later in the paper anyway.
To this end,
we denote by $\tau_G \colon \rbar C_*G \to C_*G$
the universal twisting cochain (which is acyclic) for the
(reduced normalized) bar construction $\rbar C_*G$
of the augmented differential graded algebra $C_*G$.
As explained earlier, by taking the first Eilenberg subcomplex,
we will 
henceforth take  $C_*(BG)$ to be connected, so that
$C_*(BG)$ acquires a cocomplete coaugmented differential 
graded coalgebra structure.
We remind the reader that
$\rbar C_*G$ is automatically cocomplete.

\proclaim{Lemma 4.2}
The data determine an acyclic twisting cochain
$$
\tau\colon C_*(BG) \to C_*G.
\tag4.2.1
$$
Furthermore, the data determine a chain map
from
$C_*(EG\times_GX)$
to $(\rbar C_*G) \otimes_{\tau_G}C_*X$ which
is a chain equivalence and a
morphism of differential graded
$(\rbar C_*G)$-comodules
via the adjoint $\overline \tau \colon C_*(BG) \to \rbar C_*G$
of $\tau$.
This chain map
is natural
in $X$, $G$, and
the
$G$-action.
\endproclaim

\demo{Proof}
The twisted {\it Eilenberg-Zilber\/} theorem \cite\brownez,
\ \cite\gugenhtw\ yields a twisting cochain 
\linebreak
$\tau \colon C_*BG \to C_*G$
and a contraction
$$
\Nsddata {C_*(EG \times_GX)}{\nabla}{\pi}
{(C_*BG) \otimes_{\tau}C_*X}h
\tag4.2.2
$$
of
$C_*(EG \times_GX)$ onto $(C_*BG) \otimes_{\tau}C_*X$
in such a way that $\pi$
is a morphism of $(C_*BG)$-comodules. 
This contraction is, furthermore, natural in the data.
The adjoint $\overline \tau \colon C_*BG \to \rbar C_*G$
of $\tau$ is
a morphism of differential graded coalgebras 
as well as a chain equivalence; it induces, in turn,
the chain map 
$$
\overline \tau  \otimes \roman{Id}
\colon (C_*BG)\otimes_{\tau}C_*X @>>> (\rbar C_*G) \otimes_{\tau_G}C_*X
\tag4.2.3
$$
which, via $\overline
\tau$, 
 is a morphism of differential graded
$(\rbar C_*G)$-comodules. A standard spectral sequence comparison argument
shows that (4.2.3) is a chain equivalence.
Indeed, the filtrations
of 
$C_*BG$ and $\rbar C_*G$ relative to the ordinary degree
yield Serre spectral sequences
and
(4.2.3) induces
a morphism of spectral sequences.
Both spectral sequences have $E_2$ isomorphic to
$\roman H_*(BG,\roman H_*(X))$
and (4.2.3) induces
an isomorphism
of spectral sequences from $E_2$ on whence 
(4.2.3) is an isomorphism on homology.
Since the chain complexes on both sides of
(4.2.3) are free over the ground ring,
(4.2.3) is a chain equivalence.

The following reasoning avoids spectral
sequences: The twisted objects \linebreak
$(C_*BG) \otimes_{\tau}C_*G$ and
$(\rbar C_*G) \otimes_{\tau_G}C_*G$ are both acyclic constructions
for the differential graded algebra $C_*G$, cf. e.~g.
\cite\mooretwo\ for the notion of acyclic construction, 
and choices of contracting
homotopies for these constructions determine canonical
$(C_*G)$-linear comparison maps
$$
\alpha\colon (C_*BG) \otimes_{\tau}C_*G @>>> (\rbar C_*G) \otimes_{\tau_G}C_*G,
\quad
\beta \colon (\rbar C_*G) \otimes_{\tau_G}C_*G @>>> (C_*BG) \otimes_{\tau}C_*G
$$
together with 
chain homotopies $\alpha \beta \cong
\roman  {Id}$ and $ \beta \alpha \cong \roman  {Id}$, cf. e.~g.
\cite \mooretwo\ for details. Things may in fact be arranged in
such a way that
$$
\alpha = \overline \tau  \otimes \roman{Id}
\colon (C_*BG)\otimes_{\tau}C_*G @>>> (\rbar C_*G) \otimes_{\tau_G}C_*G.
$$
Now  $\beta$ induces a morphism
$$
\beta^{\sharp} \colon (\rbar C_*G) \otimes_{\tau_G}C_*X @>>> (C_*BG)\otimes_{\tau}C_*X,
$$
and the two chain homotopies induce the requisite chain homotopies
so that $\beta^{\sharp}$ and (4.2.3)
are chain homotopy inverse to each other. 
Notice that
the chain inverse $\beta^{\sharp}$ need  not be compatible with the
$(\rbar C_*G)$-comodule structures, though.

The two chain equivalences $\pi$ and (4.2.3)
combine to the chain equivalence
$$
C_*(EG \times_GX) 
@>{\pi}>>
(C_*BG)\otimes_{\tau}C_*X 
@>{\overline \tau  \otimes \roman{Id}}>> 
(\rbar C_*G) \otimes_{\tau_G}C_*X
$$
which, by construction, is a morphism of
differential graded
$(\rbar C_*G)$-comodules
via the adjoint $\overline \tau \colon C_*(BG) \to \rbar C_*G$
of $\tau$.
The construction is natural in $X$, $G$, and the $G$-action. \qed
\enddemo

We note that the proof of the twisted Eilenberg-Zilber theorem
in \cite{\gugenhtw}\ involves the
{\it perturbation lemma\/}.
Thus Lemma 4.2 relies on HPT as well.

Let $\xi\colon P \to B$
be a principal right $G$-bundle, let $F$ be a left $G$-space, and consider
the associated
fiber bundle $P\times_GF \to B$.
Let $C$ be a differential graded coalgebra, let $\phi \colon C_*B
\to C$ be a morphism of differential graded coalgebras which is as
well a chain equivalence,  and let $\vartheta \colon C \to C_*G$
be a twisting cochain. The twisted tensor product $C \otimes
_{\vartheta} C_* F$ together with a chain map
$$
C_*(P\times_GF) @>>> C \otimes _{\vartheta} C_* F
$$
which is a chain equivalence and a morphism of $C$-comodules via
$\phi$ is a model for the chains of $P\times_GF $
that is, furthermore, compatible with the induced
$C$-comodule structures. Likewise,
 the twisted Hom-object
$\roman{Hom}^{\vartheta}(C,C_* F)$ together with a (co)chain map
$$
\roman{Hom}^{\vartheta}(C,C_* F)
@>>>
C^*(P\times_GF)
$$
which is a chain equivalence and a morphism of $C^*$-modules via
$\phi^*\colon C^* \to C^*B $ is a model
for the cochains of $P\times_GF $ that is, furthermore, compatible with the 
induced $C^*$-module structures.

We now spell out a consequence of Lemma 4.2 which, in turn,
immediately implies the statement of Theorem 4.1.

\proclaim{Corollary 4.3}
The twisted tensor product $(\rbar C_*G) \otimes_{\tau_G}C_*X$,
together with {\rm (4.2.3)} and either of the morphisms
$\pi$ or $\nabla$ in {\rm (4.2.2)},
is a model for the chains of the Borel construction $EG \times_GX$.
Likewise the twisted {\rm Hom}-object
${\roman{Hom}^{\tau_G}(\rbar C_*(G), C^*X)}$,
together with the corresponding isomorphism
of the kind {\rm (1.2.6)}
and the duals of the comparison maps between 
$(\rbar C_*G) \otimes_{\tau_G}C_*X$ and
$C_*(EG \times_GX)$ just spelled out
is a model for the cochains of the Borel construction $EG \times_GX$.
These models are natural in $X$, $G$, and the $G$-action. \qed
\endproclaim

A variant of Theorem 4.1 merely for cohomology is obtained in the following
manner: Since in general the dual of a tensor product is not the tensor 
product of the duals, 
the cochains $C^*G$ on $G$ do not inherit a coalgebra structure;
a replacement for the missing comodule structure on $C^*X$
(dual to the $C_*G$-module structure on $C_*X$ exploited above)
is provided by the following construction:
The dual $\rbar^* (C_*G)$ of the reduced bar construction $\rbar (C_*G)$
acquires a differential graded algebra structure.
The cochains
$C^*(EG \times_G X)$ inherit a canonical $C^*(BG)$-module structure and hence,
via the dual of the adjoint of the twisting cochain (4.2.1), a canonical
$\rbar^* (C_*G)$-module structure. The forgetful map from
$C^*(EG \times_G X)$ to $C^*X$ (which forgets the $G$-action) is a chain map
which extends canonically to a projection
$$
\pi \colon \roman B(R,\rbar^* (C_*G),C^*(EG \times_G X)) @>>>
C^*(X). \tag4.1.1${}^*$
$$
When $G$ is connected and when
the homology $\roman H_*(BG)$ of $BG$ is of finite type,
in view of Lemma 3.14,
$\pi$ is an isomorphism on homology.
We will refer to $\roman
B(R,\rbar^* (C_*G),C^*(EG \times_G X))$ as the {\it extended
cochain complex\/} of $X$, with reference to $G$. For
intelligibility, we recall that, by the construction of 
the extended cochain complex, this complex
inherits
a canonical (differential graded) $\rbar\, \rbar^*
(C_*G)$-comodule structure. 
We note that, for the present constructions, 
we cannot naively work with
a reduced cobar construction
on $C^*G$ since the latter differential graded
algebra does not acquire a coalgebra structure.

The following is a variant of Theorem
4.1, in the realm of ordinary singular {\it co\/}homology.

\proclaim{Theorem 4.1${}^*$} Suppose that $G$ is connected and that
the homology
$\roman H_*(BG)$ is of finite type.
Then the differential graded coalgebra $\rbar\, \rbar^* (C_*G)$
is a model for the cochains on $G$
in such a way that the diagonal map of
$\rbar\, \rbar^* (C_*G)$
serves as a replacement for the missing coalgebra structure on $C^*G$ and
that the $G$-equivariant cohomology of $X$ is canonically
isomorphic to the differential graded
$$
\roman{Cotor}^{\rbar\, \rbar^* (C_*G)} (R, 
\roman B(R,\rbar^* (C_*G),C^*(EG \times_G X))) . \tag4.1.2${}^*$
$$
\endproclaim

\demo{Proof} With the notation $A=\rbar^* (C_*G)$ and
$\tau^{\rbar}\colon \rbar A \to A$ for the universal twisting
cochain for the (reduced normalized) bar construction on $A$, the
differential graded Cotor is the homology of the twisted object $
A \otimes_{\tau^{\rbar}} \roman B(R,A,C^*(EG \times_G X)) $ which
can be written as
$$
A \otimes_{\tau^{\rbar}} \rbar A
\otimes_{\tau^{\rbar}}  C^*(EG \times_G X),
$$
and the morphism
(3.4$_*$)
yields  a projection onto $C^*(EG \times_G X)$
which is in fact a relatively projective
resolution of $C^*(EG \times_G X)$ in the category of (left)
$\rbar^* (C_*G)$-modules. \qed
\enddemo

We now proceed towards the construction of small models for equi\-variant
(co)ho\-mo\-lo\-gy.
The group $G$ is said to be of {\it strictly exterior type\/} 
(over the ground ring $R$) 
provided,
as a graded Hopf algebra, its homology $\roman H_*(G)$ is
the exterior algebra $\Lambda[x_1,\dots ]$ in odd degree primitive generators
$x_1,\dots$; when the number of generators is infinite, we will assume
throughout that $\roman H_*(G)$ is of finite type. 
With this assumption, in each degree,
$\roman H_*(G)$ is plainly a free $R$-module of finite dimension.
A detailed discussion of the property of being of strictly exterior type
can be found in \cite\husmosta\ (IV.7).
When $G$ is of strictly exterior type, it is
necessarily connected, and $\Sigm'=\roman H_*(BG)$ is a 
{\it simply connected\/}
graded cocommutative
coalgebra, indeed, the graded symmetric coalgebra on the suspension of the
free $R$-module generated by $x_1,\dots $ (the module of indecomposables).
Henceforth we suppose that $G$ is of strictly exterior type
in such a way that the duals of the primitive generators
are universally transgressive. The group $G$ being of strictly exterior
type, it is necessarily connected and,
in view of an observation made earlier,
we can take $C_*G$ to be a
connected differential graded algebra if need be.

Recall for intelligibility that, given the fibration $f \colon E\to B$,
with typical fiber $F$,
the differential 
$
d_p\colon \roman E_p^{0,p} \longrightarrow \roman E_p^{p+1,0}
$
($p \geq 2$) in the cohomology spectral sequence $(\roman E_*, d_*)$
of the fibration
determines an additive relation
$$
\tau\colon \roman E_2^{0,p} \rightharpoonup \roman E_2^{p+1,0}
$$
referred to as {\it transgression\/}, cf. \cite\maclaboo\  (XI.4 p.~332);
this is not the original definition of transgression but it is
an equivalent notion. The elements of 
$\roman E_p^{0,p}\subseteq  \roman E_2^{0,p}\cong \roman H^p(F)$
are then referred to as {\it transgressive\/}.
When $f$ is the universal bundle of a Lie group, the transgressive elements
are said to be {\it universally transgressive\/}.  

The following lemma yields an sh-morphism from 
$C_*G$ to  $\roman H_*G$ which is a quasi-isomorphism.

\proclaim{Lemma 4.4${}_*$} For each $x_j$, choose a cocycle of
$\rbar C_*(G)$ such that the dual of $x_j$ in $\roman H^*G$
transgresses to the class of this cocycle. This choice determines
an acyclic twisting cochain
$$
\zeta^B\colon \rbar C_*G @>>> \roman H_*G.
\tag4.4.1${}_*$
$$
Consequently, for any differential graded right $(\rbar C_*G)$-comodule
$M$, 
the morphism {\rm (3.3$_*$)} (for right comodules
rather than left ones) yields an injection
$$
M @>>> M \otimes_{\zeta^B} (\roman H_*G) \otimes_{\zeta^B} (\rbar C_*G)
\tag4.4.2{}$_*$
$$
of $M$ into the
relatively injective twisted object
$M \otimes_{\zeta^B} (\roman H_*G) \otimes_{\zeta^B} (\rbar C_*G)$
which is a morphism of  right $(\rbar C_*G)$-comodules
and a chain equivalence
and hence serves as a replacement for a
relatively injective resolution
of $M$ in the category of differential graded right $(\rbar C_*G)$-comodules.
\endproclaim

A version of this lemma  may be found in Section 4 of \cite\habili.
Alternative constructions for a twisting cochain of the kind (4.4.1${}_*$)
may be found in \cite\husmosta\ (IV.7.3) or may be deduced from a
corresponding construction in \cite\gugenmay. Since \cite\habili\ is not
easily available, and for later reference, we sketch a proof of the lemma.

\demo{Proof} The Eilenberg-Zilber theorem yields a contraction
$$
\Nsddata
{\rbar (C_*G\otimes C_*G)}{\nabla}{\pi}
{\rbar C_*G\otimes \rbar C_*G}h,
\tag4.4.3
$$
and, as noted earlier,
the shuffle map $\nabla$ is well known to be a morphism of differential
graded coalgebras. Write $\Lambda =\roman H_*(G)= \Lambda[x_1,\dots ]$
and, for each $j \geq 1$, write $\Lambda^j = \Lambda[x_1,\dots, x_j]$.
For each $x_j$, the chosen cocycle amounts to a twisting cochain
$\zeta_j \colon \rbar C_*(G) \to \Lambda [x_j]$. By induction, we construct a
sequence of twisting cochains $\zeta^j \colon  \rbar C_*(G) \to \Lambda^j$
($j \geq 1$).

The induction starts with $\zeta^1=\zeta_1$. Let $j \geq 1$, and suppose,
by induction, that the twisting cochain
$\zeta^j \colon \rbar C_*(G) \to \Lambda^j$
has been constructed. Then
$$
\zeta^j \otimes \eta \varepsilon
+ \eta \varepsilon \otimes \zeta_{j+1}
\colon
\rbar C_*(G)
\otimes
\rbar C_*(G)
@>>>
\Lambda^j \otimes \Lambda[x_{j+1}] =\Lambda^{j+1}
$$
is a twisting cochain. The construction (2.2.1$_*$), applied to
$\sigma =\zeta^j \otimes \eta \varepsilon
+ \eta \varepsilon \otimes \zeta_{j+1}$ and (4.4.3), yields a  twisting cochain
from $\rbar (C_*G\otimes C_*G)$ to $\Lambda^{j+1}$  which, combined with
$$
\rbar \Delta \colon \rbar C_*G
\to \rbar (C_*G\otimes C_*G),
$$
yields the twisting cochain
$\zeta^{j+1} \colon \rbar C_*G @>>> \Lambda^{j+1}$.
This completes the inductive step. When the generators
of $\roman H_*G$ constitute a finite set, the construction
stops after finitely many steps and yields the twisting cochain $\zeta^B$;
otherwise, in view of the assumption that $\roman H_*G$ be of finite type,
the $\zeta^j$'s converge to a twisting cochain $\zeta^B$;
in fact, the convergence is naive in the sense that, in each
degree, the limit is achieved after finitely may steps.

An elementary spectral sequence argument shows that the construction
\linebreak
$(\rbar C_*G) \otimes_{\zeta^B} \roman H_*G$ is acyclic. Hence the twisting
cochain $\zeta^B$ is acyclic. \qed
\enddemo

Under the circumstances of Lemma 4.4, extend the adjoint
$\overline \zeta^B\colon \rcob\, \rbar C_*(G) \to \roman H_*G$ of
$\zeta^B$ to a contraction
$$
\Nsddata {\rcob\, \rbar C_*(G)}
{\nabla}{\overline \zeta^B}{\roman H_*G}h
\tag4.5.1
$$
where the notation $\nabla$ and $h$ is abused somewhat.
Indeed, by construction, $\overline \zeta^B$ is 
surjective and an isomorphism on homology and
$\rcob\, \rbar C_*(G)$ and ${\roman H_*G}$ are both free
as graded modules over the ground ring. Consequently the kernel
of $\overline \zeta^B$ has zero homology and is even contractible,
and we can extend $\overline \zeta^B$ to data of the kind (2.1.1)
satisfying (2.1.2). As explained in Remark 2.1.4, a 
suitable modification of the homotopy if need be then
yields a homotopy such that the side conditions (2.1.3)
are satisfied as well.

The following lemma yields an sh-inverse of the sh-morphism from 
$C_*G$ to  $\roman H_*G$ constructed in Lemma 4.4$_*$ above.

\proclaim{Lemma 4.5}
Under the circumstances of Lemma {\rm 4.4}, the choice of cocycle of
$\rbar C_*(G)$ for each $x_j$ and the contraction {\rm (4.5.1)}
determine an acyclic twisting cochain
$$
\vartheta \colon \rbar \roman H_*(G)
@>>>
C_*G
\tag4.5.2
$$
such that the adjoints $\rbar C_*G @>>> \rbar \roman H_*(G)$ and
$\rbar \roman H_*(G)@>>> \rbar C_*G$ of {\rm (4.4.1${}_*$)} and {\rm (4.5.2)},
respectively, are mutually chain homotopy inverse to each other.
In particular, the composite
$$
\roman H_*(BG) @>>> \rbar C_*G \tag4.5.3
$$
of the canonical injection $\roman H_*(BG) \to \rbar \roman
H_*(G)$ with the adjoint $\rbar \roman H_*(G)@>>> \rbar C_*G$ of
{\rm (4.5.2)} is a morphism of differential graded coalgebras
inducing an isomorphism on homology.
\endproclaim

\demo{Proof}
The construction (2.2.1$^*$), applied to
$\sigma =\tau_{\rbar \roman H_*G} \colon \rbar \roman H_*G \to \roman H_*G$
(the universal twisting cochain for $\roman H_*G$) and (4.5.1), yields a
twisting cochain $\rbar \roman H_*G  @>>> \rcob\, \rbar C_*(G)$.
The composite of this twisting cochain with the universal morphism
$\rcob\, \rbar C_*(G) \to C_*(G)$ of differential graded algebras yields
the twisting cochain (4.5.2).

A spectral sequence argument shows that the construction
$\rbar \roman H_*(G) \otimes_{\vartheta} \roman C_*G$
is acyclic. Hence the twisting
cochain (4.5.2) is acyclic.

Standard comparison arguments
involving the two acyclic constructions
\linebreak
$(\rbar C_*G) \otimes_{\zeta^B} \roman H_*G$
and
$(\rbar \roman H_*(G)) \otimes_{\vartheta} \roman C_*G$
show that the adjoints of
{\rm (4.4.1${}_*$)}
and {\rm (4.5.2)} are mutually chain homotopy inverse to each
other. \qed
\enddemo

We will denote the twisting cochain which is the composite of (4.5.3) with
the universal twisting cochain $\tau_G$ of $C_*G$ by
$$
\zeta_G \colon \roman H_*(BG) @>>> C_*G .
\tag4.4.1${}^*$
$$
This twisting cochain is plainly acyclic.
It is logically the adjoint (or dual) of the twisting cochain
(4.4.1${}_*$) whence the numbering.
The description of (4.4.1${}^*$) involves the injection (4.5.3), though,
whence there seems no way to avoid this kind of numbering.

The statement of the following lemma is now immediate; yet we spell it out
to bring the duality between $G$ and $BG$ to the fore.

\proclaim{Lemma 4.4${}^*$}
For any differential graded left $(C_*G)$-module
$N$, the 
morphism {\rm (3.4$_*$)}
yields a $(C_*G)$-linear
projection
$$
(C_*G) \otimes_{\zeta_G} (\roman H_*BG) \otimes_{\zeta_G} N @>>> N
\tag4.4.2${}^*$
$$
from
the relatively projective twisted object
$(C_*G) \otimes_{\zeta_G} (\roman H_*BG) \otimes_{\zeta_G} N$
onto $N$ which is a chain equivalence
and hence serves as a replacement for a relatively projective resolution
of $N$ in the category of differential graded
left $(C_*G)$-modules. \qed
\endproclaim

\proclaim{Theorem 4.6}
The morphism
$$
\overline{\zeta_G} \otimes \roman{Id}\colon
\roman H_*(BG) \otimes_{\zeta_G} C_*X
@>>> (\rbar C_*(G)) \otimes_{\tau_G} C_*X
\tag4.6.1${}_*$
$$
of twisted tensor products induces an isomorphim on homology,
and so does the morphism
$$
\roman{Hom}^{\zeta_G}(\roman H_*(BG), C^*X)
@<<<
\roman{Hom}^{\tau_G}(\rbar C_*(G), C^*X).
\tag4.6.1${}^*$
$$
Consequently the twisted tensor product
$$
\roman H_*(BG) \otimes_{\zeta_G} C_*X
\tag4.6.2${}_*$
$$
is a model for chains of the Borel construction of $X$ in the sense that it
calculates the $G$-equivariant homology of $X$
as an $\roman H_*(BG)$-comodule, and the twisted Hom-object
$$
\roman{Hom}^{\zeta_G}(\roman H_*(BG), C^*X)
\cong \roman H^*(BG)\otimes_{\zeta_G^*} C^*X
\tag4.6.2${}^*$
$$
is a model for cochains of the Borel construction of $X$ in the sense that it
calculates the $G$-equivariant cohomology of $X$
as an $\roman H^*(BG)$-module.
\endproclaim

\demo{Proof} Standard spectral sequence comparison arguments establish the
first claim.
The following reasoning avoids spectral
sequences.

Extend the adjoint
$\overline{\zeta_G} \colon \roman H_*(BG) @>>>
\rbar C_*G$
of the acyclic twisting cochain $(4.4.1^*)$
to a contraction
$$
\Nsddata { \rbar C_*(G)}
{\overline{\zeta_G}}{\phantom{a}\pi \phantom{a}}{\roman H_*BG}h .
\tag4.6.3
$$
This contraction, in turn, induces  the contraction
$$
\Nsddata { (\rbar C_*(G))\otimes C_*X}
{\overline{\zeta_G}\otimes \roman{Id}}{\pi\otimes \roman{Id}}
{(\roman H_*BG)\otimes C_*X}
{h\otimes \roman{Id}} .
\tag4.6.4
$$
The twisted differentials on both sides of
$(4.6.1_*)$ are perturbed differentials, and an application of the perturbation
lemma yields a contraction of the kind
$$
\Nsddata { (\rbar C_*(G))\otimes_{\tau_G} C_*X}
{\overline{\zeta_G}\otimes \roman{Id}}
{\phantom{\overline{\zeta_G}}\widetilde \pi\phantom{\roman{Id}} }
{(\roman H_*BG)
\otimes_{\zeta_G} C_*X}
{\widetilde h} .
\tag4.6.5
$$
Hence $(4.6.1_*)$ is in particular a chain equivalence.
Passing to the appropriate duals, we see that
$(4.6.1^*)$ is likewise a chain equivalence. \qed
\enddemo

\noindent{\smc Remark.\/} Theorem 4.6 is implicit in \cite\habili\ 
but this reference is not generally available. At the time when 
\cite\habili\ was written, there was little interest in this kind 
of result.

\beginsection 5${}_*$. Koszul duality for homology 

Let $V$ be a free graded $R$-module of finite type concentrated
in odd positive degrees, let $\Lambda = \Lambda[V]$,
the exterior $R$-algebra on $V$,
and let $\Sigm'=\Sigm'[sV]$,
the symmetric coalgebra on the suspension $sV$;
thus 
$\Lambda$
is a strictly exterior $R$-algebra of finite type, 
$\Sigm'$ is the symmetric coalgebra on the suspension $sV$
 of the indecomposables of $\Lambda$, and
the canonical injection
$\Sigm'\to \rbar \Lambda$ is a homology isomorphism.
Notice that both
$\Sigm'$ and $\rbar \Lambda$ are simply connected.
To adjust to the notation established in 
Subsection 1.6 and in Section 3
above, let
$$
(C,A,C',A',\zeta^{C'},\zeta_{A'})=
(\Sigm', \Lambda, \rbar \Lambda, \rcob \Sigm',
\tau^{\rbar \Lambda},\tau_{\rcob \Sigm'})
$$
and let $\tau\colon \Sigm' \to \Lambda$ be the universal
(acyclic) twisting cochain determined by the desuspension map $sV \to V$.
For simplicity, write $\tau_{\rcob}= \tau_{\rcob \Sigm'}$
and $\tau^{\rbar} =\tau^{\rbar \Lambda}$.

Non-negatively graded objects will be indicated with the notation
${}_{\geq 0}$ if need be. Consider the functors
$$
\alignat2
t^{\infty} &\colon {}_{\Lambda}\roman{Mod}^{\infty}
@>>> {}_{\Sigm'}\roman{Comod},
\quad
t^{\infty}(N,\tau_{\rcob})
&&= \Sigm' \otimes_{\tau_{\rcob}} N
\tag5.1${}_*$
\\
h^{\infty} &\colon {}_{\Sigm'}\roman{Comod}^{\infty}
@>>> {}_{\Lambda}\roman{Mod},
\quad
h^{\infty}(M,\tau^{\rbar})
&&= \Lambda \otimes_{\tau^{\rbar}} M.
\tag5.2${}_*$
\endalignat
$$
Plainly, these functors restrict to functors
$$
t^{\infty} \colon {}_{\Lambda}\roman{Mod}^{\infty}_{\geq 0}
@>>> {}_{\Sigm'}\roman{Comod}_{\geq 0},
\quad
h^{\infty} \colon {}_{\Sigm'}\roman{Comod}^{\infty}_{\geq 0}
@>>> {}_{\Lambda}\roman{Mod}_{\geq 0}.
$$

\proclaim{Proposition 5.3${}_*$} The functors $t^{\infty}$ and
$h^{\infty}$ are quasi-inverse to each other and yield a
quasi-equivalence
 of categories between
${}_{\Lambda}\roman{Mod}^{\infty}$ and
${}_{\Sigm'}\roman{Comod}^{\infty}$
which restricts to a quasi-equivalence
of categories between
${}_{\Lambda}\roman{Mod}^{\infty}_{\geq 0}$
and
${}_{\Sigm'}\roman{Comod}^{\infty}_{\geq 0}$.
\endproclaim

We now explain the meaning of this proposition and in particular that
of the terms \lq\lq quasi-inverse\rq\rq\ and \lq\lq quasi-equivalence\rq\rq: 
For any left
sh-comodule $(M,\tau^{\rbar})$ over $\Sigm'$, the relatively
injective object (3.3)${}_*$ with $C= \rbar \Lambda$
has the form
$$
M @>>> \rbar \Lambda\otimes_{\tau^{\rbar}} \Lambda
\otimes_{\tau^{\rbar}} M;
$$
this injection
and the obvious injection
$$
t(h^{\infty}(M,\tau^{\rbar})) = \Sigm' \otimes_{\tau} \Lambda
\otimes_{\tau^{\rbar}} M @>>> \rbar \Lambda\otimes_{\tau^{\rbar}}
\Lambda \otimes_{\tau^{\rbar}} M
$$
are both morphisms of left $\rbar \Lambda$-comodules, that is,
morphisms of sh-comodules over $\Sigm'$, both morphisms are
quasi-isomorphisms, and these are natural in the sh-comodule
$(M,\tau^{\rbar})$ over $\Sigm'$.
Likewise, for any left sh-module $(N,\tau_{\rcob})$ over
$\Lambda$, the relatively projective object (3.4)${}_*$ with
$A= \rcob \Sigm'$  has the form
$$
\rcob \Sigm'\otimes_{\tau_{\rcob}} \Sigm' \otimes_{\tau_{\rcob}} N
@>>> N;
$$
this map and the obvious surjection
$$
\rcob \Sigm'\otimes_{\tau_{\rcob}} \Sigm' \otimes_{\tau_{\rcob}} N
@>>> \Lambda \otimes_{\tau} \Sigm' \otimes_{\tau_{\rcob}} N
=h(t^{\infty}(N,\tau_{\rcob}))
$$
are both morphisms of left $(\rcob \Sigm')$-modules, that is,
morphisms of sh-modules over $\Lambda$, both morphisms are
quasi-isomorphisms, and these are natural in the sh-module
$(N,\tau_{\rcob})$ over $\Lambda$.
It is in this sense that the functors $t^{\infty}$ and
$h^{\infty}$ are {\it quasi-equivalences which are quasi-inverse to
each other\/}. These observations imply a proof of
Proposition 5.3${}_*$. 

\smallskip
\noindent {\smc 5.4${}_*$. Koszul duality in ordinary
equivariant homology,\/} cf. \cite\franzone, \cite\franztwo,
\cite\gorkomac. Let $G$ be a topological group of strictly
exterior type such that the
duals of the exterior generators
are universally transgressive and let $\Lambda = \roman H_*G$ and $\Sigm' = \roman
H_*(BG)$. {\sl Any left $(C_*G)$-module is a left sh-module over
$\roman H_*G$ via the twisting cochain\/} $\zeta_G\colon \roman
H_*(BG) \to C_*G$ given in (4.4.1${}^*$) and  {\sl any left
$(\rbar C_*G)$-comodule is a left sh-comodule over $\roman
H_*(BG)$\/} via the twisting cochain $\zeta^B\colon \rbar C_*G \to
\roman H_*G$ given in (4.4.1${}_*$). In particular, given a left
$G$-space $X$, the chain complex
$C_*X$ is a left sh-module over $\roman H_*G$ via
$\zeta_G$ and, for any space $Y$ over $BG$, the chains $C_*Y$
constitute a left sh-comodule over $\roman H_*(BG)$ via $\zeta^B$
and the twisting cochain (4.2.1). Thus the functors (5.1${}_*$)
and (5.2${}_*$) are defined, and the {\sl statement of
Proposition\/} 5.3${}_*$ {\sl entails that ordinary and
equivariant homology are related by Koszul duality in the
following sense: On the category of (left) $G$-spaces, the functor
$h\circ t^{\infty} \circ C_*$ is chain-equivalent to the functor
$C_*$ as sh-module functors over $\Lambda = \roman H_*G$; and on
the category of spaces over $BG$, the functor $t\circ h^{\infty}
\circ C_*$ is chain-equivalent to the functor $C_*$ as sh-comodule
functors over $\Sigm' = \roman H_*(BG)$.\/} Thus,
given the left $G$-space $X$, the $\Sigm'$-comodule
$t^{\infty}(C_*(X))$, being a model for the chains of 
the Borel construction
$N(G,X)$,  cf. Theorem 4.6, calculates the $G$-equivariant homology of 
$X$; 
application of the functor $h$ to the twisted object $t^{\infty}
(C_*X)$ calculating the $G$-equivariant homology of $X$ then
reproduces
an object calculating the ordinary homology of $X$. 
Likewise,
given the space $Y \to BG$ over $BG$,
the twisted object $h^{\infty} (C_*Y)$ is a model for the chains
of the total space of the $G$-bundle over $Y$ induced 
via  $Y \to BG$ and
application of the functor $t$ to  $h^{\infty}
(C_*Y)$ reproduces an object calculating the ordinary homology of
$Y$.

\beginsection 5${}^*$. Koszul duality for cohomology

We maintain the circumstances of the previous section.
Let $\Sigm = \roman {Hom}(\Sigm',\Bobb R)$; since $V$ is of finite type,
$\Sigm$ amounts to $\Sigm[s^{-1}V^*]$,
the symmetric algebra on the desuspension
of the graded dual $V^*$ of $V$.
Likewise,
$\Lambda' = \roman {Hom}(\Lambda,\Bobb R)$
amounts to the exterior coalgebra cogenerated by $V^*$,
and the induced twisting cochain
$\tau^* \colon \Lambda' \to \Sigm$ is acyclic.
Non-positively graded objects will be indicated with the notation
${}_{\leq 0}$  if need be.
For intelligibility, we note that
{\sl a (left) sh-module over $\Sigm$ is a (left) $(\rcob \Lambda')$-module\/},
and that
{\sl a (right) sh-module over $\Lambda$ is a (right) $(\rcob \Sigm')$-module.\/}
Consider the functors
$$
\alignat2
t^*_{\infty} &\colon \roman{Mod}^{\infty}_{\Lambda}
@>>>{}_{\Sigm}\roman{Mod},
\quad
t^*_{\infty}(N^{\sharp},\tau_{\rcob})
&&=\roman{Hom}^{\tau_{\rcob}}(\Sigm', N^{\sharp}),
\tag5.1${}^*$
\\
h^*_{\infty} &\colon {}_{\Sigm}\roman{Mod}^{\infty}
@>>> \roman{Mod}_{\Lambda},
h^*_{\infty}(M^{\sharp},\tau_{\rcob \Lambda'})
&&=   \Lambda' \otimes_{\tau_{\rcob \Lambda'}}M^{\sharp}
=\roman{Hom}^{\tau^{\rbar \Lambda}}(\Lambda, M^{\sharp}),
\tag5.2${}^*$
\endalignat
$$
where $\Lambda'$ is a right
$\Lambda$-module via the left
$\Lambda$-module structure on itself,
that is,
$\Lambda'$ is a right
$\Lambda$-module via the operation of contraction.
These functors restrict to functors
$$
t^*_{\infty} \colon {}_{\leq 0}\roman{Mod}^{\infty}_{\Lambda}
@>>>{}_{\Sigm}\roman{Mod}_{\leq 0},
\quad
h^*_{\infty} \colon {}_{\Sigm}\roman{Mod}^{\infty}_{\leq 0}
@>>> {}_{\leq 0}\roman{Mod}_{\Lambda}.
$$

\proclaim{Proposition 5.3${}^*$} The functors $t^*_{\infty}$ and
$h^*_{\infty}$ are quasi-inverse to each other and yield a
quasi-equivalence of categories between
$\roman{Mod}^{\infty}_{\Lambda}$ and
${}_{\Sigm}\roman{Mod}^{\infty}$ which restricts to a
quasi-equivalence of categories between ${}_{\leq
0}\roman{Mod}^{\infty}_{\Lambda}$ and
${}_{\Sigm}\roman{Mod}^{\infty}_{\leq 0}$.
\endproclaim

The meaning of this proposition, the significance
of the terms \lq\lq quasi-inverse\rq\rq\ and 
\lq\lq quasi-equivalence\rq\rq, and the proof of the proposition  
are somewhat similar to those for Proposition 5.3$_*$.
We refrain from spelling out details.

\smallskip
\noindent {\smc 5.4${}^*$. Koszul duality in ordinary
equivariant cohomology,\/} cf. \cite\franzone, \cite\franztwo,
\cite\gorkomac. As before, let $G$ be a topological group of
strictly exterior type such that the duals of the primitive generators
are universally transgressive, let $\Lambda = \roman H_*G$, $\Sigm' =
\roman H_*(BG)$, $\Lambda' = \roman H^*G$, and $\Sigm = \roman
H^*(BG)$. {\sl Any right $(C^*G)$-comodule is a right sh-comodule
over $\roman H^*G$ via the twisting cochain\/} $\zeta^*_G\colon
C^*G \to \roman H^*(BG)$, cf. (4.4.1${}^*$), and {\sl any (left)
$(\rcob C^*G)$-module is an sh-module over\/} $\roman H^*(BG)$ via
the twisting cochain $(\zeta^B)^*\colon   \roman H^*G \to \rcob
C^*G$, cf. (4.4.1${}_*$). In particular, $C^*X$ is a right
sh-module over $\roman H_*G$ via $\zeta^*_G$ and, for any space
$Y$ over $BG$, the cochains $C^*Y$ constitute a (left) sh-module
over $\roman H^*(BG)$ via $(\zeta^B)^*$ and the dual of the
twisting cochain (4.2.1). Thus the 
functors 
(5.1${}^*$) and
(5.2${}^*$) 
are defined, and the {\sl statement of
Proposition\/} 5.3${}^*$ {\sl entails that ordinary and
equivariant cohomology are related by Koszul duality in the
following sense: On the category of left $G$-spaces, the functor
$h^*\circ t^*_{\infty} \circ C^*$ is chain-equivalent to the
functor $C^*$ as sh-comodule functors over $\Lambda'= \roman
H^*G$; and on the category of spaces over $BG$, the functor
$t^*\circ h^*_{\infty} \circ C^*$ is chain-equivalent to the
functor $C^*$ as sh-module functors over\/} $\Sigm = \roman
H^*(BG)$. Thus, given the $G$-space $X$,
the twisted object $t^*_{\infty} (C^*X)$, being a model for 
the Borel construction, calculates the
$G$-equivariant cohomology of $X$; 
application of the functor $h^*$ to 
$t^*_{\infty} (C^*X)$
then reproduces an object calculating
the ordinary cohomology of $X$.
Likewise, the twisted object
$h^*_{\infty} (C^*Y)$ is a model for the cochains of the total
space of the induced $G$-bundle over $Y$ and
application of the functor
$t^*$ to $h^*_{\infty} (C^*Y)$ reproduces an
object calculating the ordinary cohomology of $Y$. 

\medskip\noindent
{\bf 6. Grand unification}
\smallskip\noindent
Given two augmented differential graded algebras $A_1$ and $A_2$
and two sh-maps from $A_1$ to $A_2$, that is,
twisting cochains $\tau_1,\tau_2 \colon \rbar A_1 \to A_2$,
a {\it homotopy of sh-maps\/}  from
$\tau_1$ to $\tau_2$ is a homotopy of twisting cochains
 from
$\tau_1$ to $\tau_2$.
Likewise,
given two coaugmented differential graded coalgebras $C_1$ and $C_2$
and two sh-maps from  $C_1$ to $C_2$, that is,
twisting cochains $\tau_1,\tau_2 \colon  C_1 \to \rcob C_2$,
a {\it homotopy of sh-maps\/}  from
$\tau_1$ to $\tau_2$ is a homotopy of twisting cochains
 from
$\tau_1$ to $\tau_2$.

Recall that the category $\roman {DASH}_{\roman h}$ has as its objects 
augmented  differential
non-negatively graded algebras with homotopy classes of sh-maps as
morphisms;
likewise the category $\roman
{DCSH}_{\roman h}$ has as its objects 
connected differential non-negatively graded
coalgebras, necessarily coaugmented  in a unique manner,
with homotopy classes of sh-maps as morphisms.
For these matters, cf. \cite\munkhotw.
In that reference, the
ground ring is a field, but that is not strictly necessary, under
suitable additional hypotheses; for example, it suffices to
require that the graded algebras and graded coalgebras under
discussion be free as modules over the ground ring. Below, the requisite
hypotheses will always tacitly be assumed to hold.
The differential Tor can then be defined
as a functor, to be written as shTor,
on $\roman {DASH}_{\roman h}$ 
and the
differential Cotor can be defined as a functor on $\roman {DCSH}_{\roman h}$, 
to be written as shCotor;
likewise, the differential Ext can be defined as a functor,
to be written as shExt,  on
$\roman {DASH}_{\roman h}$. 
In other words, with the appropriate interpretation of the
term \lq extension\rq,
this yields functorial extensions of the originalt Tor-, Cotor-, and Ext-functors, 
and the functoriality of those extensions is then referred to as
{\it extended functoriality\/}  \cite\gugenmun.

We maintain the hypothesis that $G$ is of strictly exterior type
in such a way that the duals of the exterior homology generators
are universally transgressive.
In the category $\roman {DASH}_{\roman h}$ of differential
non-negatively graded algebras with homotopy classes of sh-maps as
morphisms, the twisting cochain $\zeta_G\colon \roman H_*(BG) \to
C_*G $ given as (4.4.1${}^*$) yields an isomorphism from $\roman
H_*G$ onto $C_*G$. Likewise, in the category $\roman
{DCSH}_{\roman h}$ of connected differential non-negatively graded
coalgebras with homotopy classes of sh-maps as morphisms, the
twisting cochain $\zeta^B \colon \rbar C_*G \to \roman H_*G$ given
as (4.4.1${}_*$) above yields an isomorphism from $\rbar C_*G$
(notice that this coalgebra is connected!) onto $\roman H_*(BG)$.
In view of the extended functoriality of the Tor- and
Cotor-functors, the Koszul duality functor
$t^{\infty}$ 
defined on the category
${}_{\Lambda}\roman{Mod}^{\infty}$ of left sh-modules over
$\Lambda$, cf. (5.1${}_*$) above,
is then equivalent to a corresponding functor
$t$
to be defined shortly on the category 
${}_{C_*G}\roman{Mod}$
of ordinary differential graded left $(C_*G)$-modules, cf. (6.1$_*$) below,
in the sense that the isomorphism in 
$\roman {DASH}_{\roman h}$
between 
$C_*G$ and $\Lambda$---natural in 
$G$---and the appropriate isomorphism between the target objects
induce an equivalence between these two functors;
likewise, the Koszul duality functor $h^{\infty}$ defined
on the category ${}_{\Sigm'}\roman{Comod}^{\infty}$
of left sh-comodules over
$\Sigm'$, cf. (5.2${}_*$) above,
is then equivalent to a corresponding functor $h$
to be defined shortly on the category 
${}_{\rbar C_*G}\roman{Comod}$
of ordinary differential graded left 
$(\rbar C_*G)$-comodules, cf. (6.2$_*$) below,
in the sense that the isomorphism in 
$\roman {DCSH}_{\roman h}$
between  $\rbar C_*G$ and $\Sigm'$---natural in 
$G$---and the appropriate isomorphism between the target objects
induce an equivalence between these two functors.
Furthermore, in view of the
extended functoriality of the Ext- and Tor-functors, the 
Koszul duality functor $t^*_{\infty}$
defined on the category
$\roman{Mod}^{\infty}_{\Lambda}$ of right sh-modules over $\Lambda$, cf. (5.1${}^*$),
is then equivalent to a corresponding functor
$t^*$ to be defined shortly on the category 
$\roman{Mod}_{C_*G}$ of ordinary differential graded right $(C_*G)$-modules,
cf. (6.1$^*$) below,
the notion of equivalence of functors being interpreted in formally 
the same way as explained above, with the 
requisite modifications being taken into account;
likewise
the Koszul duality functor $h^*_{\infty}$ 
defined on the category ${}_{\Sigm}\roman{Mod}^{\infty}$
of left sh-modules over $\Sigm$,  cf. (5.2${}^*$),
is then equivalent to a corresponding functor
$h^*$ to be defined shortly on the category 
${}_{\rbar^* C_*G}\widehat{\roman{Mod}}$ 
of complete differential graded left 
$(\rbar^* C_*G)$-modules,
cf. (6.2$^*$) below,
again the notion of equivalence of functors  being interpreted in formally 
the same way as explained above, with the 
requisite modifications being taken into account.

We will now
make this precise. The definitions of the shTor, shExt and shCotor
will be given in (6.5)--(6.8) below.

In ordinary equivariant singular (co)homology, 
the situation of Section 3 arises in the following manner:
Let $A= C_*G$, $C=\rbar C_*G$, and
$\tau=\tau_G \colon \rbar C_*G \to C_*G$.
The functors $t$, $h$, $t^*$, $h^*$ now take the form
$$
\alignat1
t &\colon {}_{C_*G}\roman{Mod} @>>> {}_{\rbar C_*G}\roman{Comod},
\quad
t(N) = \rbar C_*G \otimes_{\tau} N
\tag6.1${}_*$
\\
h &\colon  {}_{\rbar C_*G}\roman{Comod} @>>> {}_{C_*G}\roman{Mod},
\quad h(M) =  C_*G \otimes_{\tau}M
\tag6.2${}_*$
\\
t^* &\colon
\roman{Mod}_{C_*G} @>>>
{}_{\rbar^* C_*G}\widehat{\roman{Mod}},
\quad
t^*(\Nfl)
= \roman{Hom}^{\tau}(\rbar C_*G, \Nfl)
\tag6.1${}^*$
\\
h^* &\colon
{}_{\rbar^* C_*G}\widehat{\roman{Mod}}
@>>>
\roman{Mod}_{C_*G},
\quad
h^*(\Mfl)
= \roman{Hom}^{\tau}(C_*G, \Mfl) .
\tag6.2${}^*$
\endalignat
$$
In particular, the twisted object $t(C_*X)$ calculates the
$G$-equivariant homology of $X$, and $t^*(C^*X)$ calculates the
$G$-equivariant cohomology of $X$. Likewise, given a space $Y \to
BG$ over $BG$, its normalized chain complex $C_*Y$ inherits a left
$(\rbar C_*G)$-comodule structure via the adjoint $\overline \tau$
of $\tau$
in a canonical manner, and the normalized cochain complex $C^*Y$
inherits a  complete $(\rbar^* C_*G)$-module structure. The twisted
object $h(C_*Y)$ calculates the differential graded
$$
\roman{Cotor}^{\rbar C_*G}(R,C_*Y) \tag6.3${}_*$
$$
and,
in view of the discussion in (1.2) above,
under appropriate finiteness assumptions,
the twisted object $h^*(C^*Y)$ calculates
the differential graded
$$
\roman{Tor}_{\rbar^* C_*G}(R,C^*Y). \tag6.3${}^*$
$$
By the Eilenberg-Moore theorem \cite\eilenmoo, with
reference to the  fiber square
$$
\CD
E_Y @>>> EG
\\
@VVV
@VVV
\\
Y @>>> BG,
\endCD
\tag6.4
$$
$\roman{Cotor}^{\rbar C_*G}(R,C_*Y)$ calculates
the homology of $E_Y$ and
$\roman{Tor}_{\rbar^* C_*G}(R,C^*Y)$
calculates
the cohomology of $E_Y$, compatibly with the bundle structures.
This situation is {\it dual\/} to that of Theorem 4.1.
The original approach of Eilenberg and Moore involves appropriate resolutions
in the differential graded category.

Let $N$ be a $(C_*G)$-module; then $(N,\zeta_G)$ is an sh-module
over $\roman H_*G$ and, for the sake of consistency with the
definitions in (5) above, we will write this sh-module as
$(N,\tau_{\rcob})$ so that $N$ is considered as a differential
graded $(\rcob\roman H_*(BG))$-module via $\zeta_G\colon \roman
H_*(BG) \to C_*G $, cf. 
(4.4.1${}^*$). The small model
$$
t^{\infty}(N,\tau_{\rcob})=(\roman H_*(BG)) \otimes_{\zeta_G} N
\tag6.5
$$
{\it defines\/} the
$\roman {shTor}^{\roman H_*G}(R,(N,\tau_{\rcob}))$,
and the morphism
$$
\overline \zeta_G \otimes \roman{Id}\colon
(\roman H_*(BG)) \otimes_{\zeta_G} N @>>> (\rbar C_*(G)) \otimes_{\tau_G} N
$$
of twisted objects which, for $N=C_*X$, comes down to (4.6.1${}_*$) above,
makes explicit the isomorphism between
$\roman {shTor}^{\roman H_*G}(R,(N,\tau_{\rcob}))$
and
$\roman {Tor}^{C_*G}(R,N)$
induced by the isomorphism
in $\roman {DASH}_{\roman h}$
from $\roman H_*G$ onto
$C_*G$.

Likewise, let $M$ be a $(\rbar C_*G)$-comodule; then $(M,\zeta^B)$
is an sh-comodule over $\roman H_*(BG)$. Write this sh-comodule as
$(M,\tau^{\rbar})$ so that $M$ is considered as a (differential
graded) $(\rbar\roman H_*(BG))$-comodule via $\zeta^B \colon \rbar
C_*G \to \roman H_*G$, cf. 
(4.4.1${}_*$). The twisted object
$$
h^{\infty}(M,\tau^{\rbar}) =(\roman H_*G) \otimes_{\zeta^B} M
\tag6.6
$$
{\it defines\/} the
$\roman {shCotor}^{\roman H_*(BG)}(R,(M,\tau^{\rbar}))$,
and the two morphisms
$$
\align
\overline \zeta^B \otimes \roman{Id}
\colon
(\rcob\, \rbar C_*G) \otimes_{\tau_G} M
&@>>> (\roman H_*G) \otimes_{\zeta^B} M
\\
\overline \tau_{C_*G} \otimes \roman{Id}
\colon
(\rcob\, \rbar C_*G) \otimes_{\tau_G} M
&
@>>>
(C_*G) \otimes_{\tau_G} M
\endalign
$$
of twisted objects make explicit the isomorphism
$$
\roman {Cotor}^{\rbar C_*G}(R,M)
@>>>
\roman {shCotor}^{\roman H_*(BG)}(R,(M,\tau^{\rbar}))
$$
induced by the isomorphism
in $\roman {DCSH}_{\roman h}$
from
$\rbar C_*G$  onto $\roman H_*(BG)$.
Thus, in view of the extended functoriality of the Tor- and Cotor-functors,
the two Koszul duality functors $h^{\infty}$ and $t^{\infty}$
are equivalent to the functors $h$ and $t$, respectively.
The duality between the functors $t$ and $h$ described in
Theorem 3.5${}_*$
implies the
duality for the functors $h^{\infty}$ and $t^{\infty}$
spelled out in Proposition 5.3${}_*$
and hence the Koszul duality
in equivariant singular homology as given in (5.4${}_*$) above.

In the same vein:
Let $\Nfl$ be a $(C_*G)$-module; then $(\Nfl,\zeta_G)$ is an sh-module
over $\roman H_*G$.
Write this  sh-module as $(\Nfl,\tau_{\rcob})$
so that $\Nfl$ is considered as
a (differential graded) $(\rcob\roman H_*(BG))$-module via
$\zeta_G$.
The small model
$$
t^*_{\infty}(\Nfl,\tau_{\rcob})
=\roman{Hom}^{\tau_{\rcob}}(\Sigm', \Nfl)=
\roman{Hom}^{\zeta_G}(\roman H_*(BG),\Nfl)
\tag6.7
$$
{\it defines\/} the
$\roman {shExt}_{\roman H_*G}(R,(\Nfl,\tau_{\rcob}))$,
and the morphism
$$
\roman{Hom}^{\zeta_G}(\roman H_*(BG), \Nfl)
@<<<
\roman{Hom}^{\tau_G}(\rbar C_*(G),\Nfl)
$$
of twisted Hom-objects
which for $\Nfl=C^*X$ comes down to (4.6.1${}^*$) above,
makes explicit the isomorphism
between
$\roman {shExt}_{\roman H_*G}(R,(\Nfl,\tau_{\rcob}))$
and
$\roman {Ext}_{C_*G}(R,\Nfl)$
induced by the isomorphism
in $\roman {DASH}_{\roman h}$
from $\roman H_*G$ onto
$C_*G$.
Likewise,
let $\Mfl$ be a complete $(\rbar^* C_*G)$-module; then $(M,\zeta^B)$ is an
sh-module over $\roman H^*(BG)$.
Write this  sh-module as $(\Mfl,\tau_{\rcob \Lambda'})$
so that $\Mfl$ is considered as
a (differential graded) $(\rcob\roman H^*G)$-module via
$\zeta^B$. More precisely, the adjoint
$$
\overline {(\zeta^B)^*} \colon
\rcob \roman H^*G @>>> \rbar^* C_*G
$$
of the dual $(\zeta^B)^* \colon \roman H^*G @>>> \rbar^* C_*G$
makes
$\Mfl$ into a (differential graded) $(\rcob\roman H^*G)$-module.
The twisted object
$$
h^*_{\infty}(\Mfl,\tau_{\rcob \Lambda'})
=   \Lambda' \otimes_{\tau_{\rcob \Lambda'}}\Mfl
=\roman{Hom}^{\tau^{\rbar \Lambda}}(\Lambda, \Mfl)
=\roman{Hom}^{\zeta^B}(\roman H_*G, \Mfl)
\tag6.8
$$
{\it defines\/}
$\roman {shTor}_{\roman H^*(BG)}(R,(\Mfl,\tau_{\rcob \Lambda'}))$
(since $\roman H^*(BG)$ is of finite type),
and the two morphisms
$$
\align
\roman{Hom}^{\tau_G}(\rcob\, \rbar C_*G, \Mfl)
&@<<< \roman{Hom}^{\zeta^B}(\roman H_*G, \Mfl),
\\
\roman{Hom}^{\tau_G}(\rcob\, \rbar C_*G,  \Mfl)
&
@<<<
\roman{Hom}^{\tau_G}(C_*G, \Mfl)
\endalign
$$
of twisted Hom-objects,
the former being induced by
$\zeta^B$ and the
latter being the obvious one,
make explicit the isomorphism between
$\roman {Tor}_{\rbar^* C_*G}(R,\Mfl)$
and
$\roman {shTor}_{\roman H^*(BG)}(R,(\Mfl,\tau_{\rcob \Lambda'}))$
induced by the isomorphism
in $\roman {DASH}_{\roman h}$
from
$\rbar^* C_*G$  onto $\roman H^*(BG)$.
Thus, in view of the extended functoriality
of the Ext- and Tor-functors,
the two Koszul duality functors
$h^*_{\infty}$ and $t^*_{\infty}$
are equivalent to the functors
$h^*$ and $t^*$, respectively.
The duality between the functors $t^*$ and $h^*$ described in
Theorem 3.5${}^*$ implies the
duality for the functors $h^*_{\infty}$ and $t^*_{\infty}$
spelled out in Proposition 5.3${}^*$
and hence the Koszul duality
in equivariant singular cohomology as given in (5.4${}^*$) above.

\medskip\noindent
{\bf 7. Some illustrations}
\smallskip\noindent
Let $G$ be a group, possibly endowed with additional structure
(e.~g. a topological group or a Lie group).
Recall that a $G$-space is called {\it equivariantly formal\/}
over the ground ring $R$ when the spectral sequence from ordinary cohomology to
$G$-equivariant cohomology collapses, 
cf. \cite\borelboo\ and \cite\frankone. 
Thus when the $G$-space $X$ is equivariantly formal (over $R$),
the graded object $\roman E_{\infty}(X)$ associated with  the
$G$-equivariant cohomology $\roman H^*_G(X)$ relative to the Serre filtration
(the filtration  coming from $\roman H^*(BG)$-degree)
takes the form of 
$\roman H^*(BG \times X)$ but this does not imply that
the equivariant cohomology is, as a graded $\roman H^*(BG)$-module,
an induced module of the kind 
$\roman H^*(BG, \roman H^*(X))$
(or of the kind  $\roman H^*(BG) \otimes \roman H^*(X)$
when $\roman H^*(BG)$ is free over the ground ring)
unless the additive extension problem is trivial,
for example when the ground ring is a field.
The property that the equivariant cohomology is such an induced module
is strictly stronger than equivariant formality.
See \cite\franpupp\ Theorem 1.1 for a discussion and, in particular,
Example 5.2 in that paper. Below we shall come back to the difference between
the two properties.

According to
an observation in \cite\kirwaboo\ (proof of Proposition 6.8), 
a smooth compact symplectic
manifold, endowed with a hamiltonian action of a compact Lie
group, is equivariantly formal over the reals. 
This fact is also an
immediate consequence of the result in \cite\frankone\ which says that,
given a torus $T$ and a compact hamiltonian $T$-manifold $X$,
the real $T$-equivariant cohomology of $X$ is an extended
$\roman H^*(BT)$-module of the kind
$$
\oplus \,\roman H^*(BT)\otimes \roman H^{*-\lambda_k}(F_k)
$$
where $F_1,\ldots, F_s$ are the finitely many components of
the fixed point set $X^T$ and where $\lambda_1,\ldots, \lambda_s$
are even natural numbers.
Here the
compactness of the manifold is crucial as the following example
shows, which also serves as an illustration for the notions of
twisting cochain etc. exploited above:

\head Example 7.1 \endhead 

The familiar $S^1$-action on $X=\Bobb
C^n\setminus {0}$ is free and hence cannot be equivariantly
formal. This action is plainly hamiltonian. The induced action of
$\roman H_*(S^1)$ on $\roman H_*(X)$ is manifestly trivial and
hence lifts to the trivial action of $\roman H_*(S^1)$ on
$C_*(X)$. Over the integers $\Bobb Z$ as ground ring, we now
describe the corresponding sh-action of $\roman H_*(S^1)$ on
$C_*(X)$ which has to be non-trivial since the $S^1$-action on $X$
cannot be equivariantly formal. This example is, admittedly,
nearly a toy example, but it illustrates how under such
circumstances homological perturbation theory works in general. In
particular, it illustrates the fact spelled out already in the
introduction that, for a general Lie group $G$ and a general
$G$-space $Y$, even when the induced action of $\roman H_*(G)$ on
$\roman H_*(Y)$ lifts to an action of $\roman H_*(G)$ on $C_*(Y)$,
in general only an sh-action of $\roman H_*(G)$ on $C_*(Y)$
involving higher degree terms will recover the geometry of the
action.

As a coalgebra, the homology $\roman
H_*(BS^1)$ is the symmetric coalgebra cogenerated by a single
generator $u$ of degree 2. This is the coalgebra which underlies
the divided polynomial Hopf algebra $\Gamma = \Gamma [u]$  on $u$.
This algebra, in turn, has a generator $\gamma_i(u)$ of degree
$|\gamma_i(u)| = 2i$ in each degree $i\geq 1$ where $\gamma_1(u)
=u$, subject to the relations
$$
\gamma_{i}(u)\gamma_{j}(u)=\tbinom{i+j}j  \gamma_{i+j}(u), \quad
i,j \geq 1,
$$
and the diagonal map $\Delta$ is given by the formula
$$
\Delta\gamma_i(u)=\sum_{j+k=i}\gamma_j(u)\otimes \gamma_k(u),
\quad i \geq 1.
$$
Let
$$
\tau \colon \Gamma [u] @>>> \roman{End}(\roman H_*(S^{2n-1}))
$$
be the twisting cochain which is zero on $\gamma_j(u)$ for $j\ne
n$ and whose value $ \tau(\gamma_n(u))$ is the endomorphism of
$\roman H_*(S^{2n-1})$ which sends $1 \in \roman H_0(S^{2n-1})$ to
one of the two generators of $\roman H_{2n-1}(S^{2n-1})\cong \Bobb
Z$. We assert that $\roman H_*(S^{2n-1})$ may be taken as a model
for the chains on $X$ and that the requisite sh-action of $\roman
H_*(S^1)$ on $C_*X$ then amounts to the twisting cochain $\tau$.

In order to justify this claim, let $v$ be a  generator of $\roman
H_1(S^1)$, write the homology algebra $\roman H_*(S^1)$ as the
exterior algebra $\Lambda [v]$, and let $\vartheta \colon \Gamma
[u]\to \Lambda [v]$ be the obvious acyclic twisting cochain which
sends $u$ to $v$ and is zero on $\gamma_j(u)$ for $j \ne 1$. Let
$w$ be an indeterminate of degree 2 and denote by $\Gamma_{n-1}
[w]$ the subcoalgebra of $\Gamma [w]$ which is additively
generated by $\gamma_1(w),\dots, \gamma_{n-1}(w)$. 
When $n=1$, this coalgebra comes down to the ground ring which is here 
taken to be that of the integers.
The coalgebra
$\Gamma_{n-1} [w]$ amounts to the homology of complex projective
$(n-1)$-space $\Bbb C P^{n-1}$,
in fact, can be taken as a model for the chains of  $\Bbb C P^{n-1}$:
Take  $C_*(\Bbb C P^{n-1})$ to be the normalized chain complex of the first 
Eilenberg subcomplex of  $\Bbb C P^{n-1}$.
A choice $\zeta\colon C_*(\Bbb C P^{n-1})\to \Bbb Z$ of 2-cocycle
representing the generator of $\roman H^2(\Bbb C P^{n-1})$
determines a twisting cochain 
$\tau\colon  C_*(\Bbb C P^{n-1})\longrightarrow \Lambda [v]$
such that, for degree reasons, the values of adjoint 
$\overline \tau\colon  C_*(\Bbb C P^{n-1})\longrightarrow \rbar\Lambda [v] \cong \Gamma [u]$
of $\tau$ lie in $\Gamma_{n-1}[u]$, that is, 
the adjoint takes the form
$$
\overline \tau\colon  C_*(\Bbb C P^{n-1})\longrightarrow \Gamma_{n-1} [u],
$$
necessarily a chain equivalence, since $\overline \tau$ induces an isomorphism
on homology and since the two chain complexes are free over the ground ring.
Alternatively, the fact that the
homology of complex projective
$(n-1)$-space $\Bbb C P^{n-1}$,
can be taken as a model for its chains results from the familiar cell decomposition
of this space.
Likewise the standard cell decomposition of the $(2n-1)$-sphere relative to
the $S^1$-action has a single cell in each degree in such a way that
the resulting cellular chain complex is exactly the chain complex which underlies
the twisted tensor product
$$
N=\Lambda[v] \otimes_{\vartheta}\Gamma_{n-1} [w]
$$
where $\Gamma_{n-1} [w]$ is viewed as a right 
$\Gamma [w]$-comodule in the obvious manner; 
this chain complex is an equivariant model for the
chains of the $(2n-1)$-sphere $S^{2n-1}$ and hence of $X$ as a free left
$S^1$-space. We remind the reader that $S^{2n-1}$ is the
{\it Stiefel manifold\/} $G(n,n-1)$ of unitary complex 1-frames in $\Bbb C^n$;
cell decompositions of general Stiefel manifolds can be found in \cite\steeepst\  (Chap. 4).

Somewhat more formally: Pick a 1-cocycle representing the generator of $\roman H^1(S^1)$
and view this 1-cocycle as a morphism $\vartheta \colon C_*(S^1) \to \Lambda[v]$
of differential graded algebras, necessarily inducing an isomorphism 
on homology. The twisted Eilenberg-Zilber theorem, applied to the the principal
$S^1$-bundle $S^{2n-1} \to \Bbb C P^{n-1}$, yields
a twisting cochain 
$\tau^{\Bbb C P^{n-1}} \colon C_*(\Bbb C P^{n-1}) \to C_*(S^1)$
and a contraction
$$
\Nsddata {C_*(S^{2n-1})}{\nabla}{\pi}
{ C_*(S^1)\otimes_{\tau^{\Bbb C P^{n-1}}} C_*(\Bbb C P^{n-1}) }h
$$
of
$C_*(S^{2n-1})$ onto $ C_*(S^1)\otimes_{\tau^{\Bbb C P^{n-1}}} C_*(\Bbb C P^{n-1}) $
that is compatible with the 
$(C_*(S^1))$-module structures.
Hence 
$ C_*(S^1)\otimes_{\tau^{\Bbb C P^{n-1}}} C_*(\Bbb C P^{n-1}) $
is a $(C_*(S^1))$-equivariant model for the chains of $S^{2n-1}$.
The composite  
$$
\vartheta \tau^{\Bbb C P^{n-1}}
\colon \Bbb C P^{n-1} \longrightarrow \Lambda[v]
$$
is plainly a twisting cochain in such a way that
$$ 
\Lambda[v] \otimes_{\vartheta \tau^{\Bbb C P^{n-1}}} C_*(\Bbb C P^{n-1})
$$
is a $(C_*(S^1))$-equivariant model for the chains of $S^{2n-1}$.
The twisting cochain $\vartheta \tau^{\Bbb C P^{n-1}}$ can be interpreted
as a choice of 2-cocycle on $C_*(\Bbb C P^{n-1})$ of the kind $\zeta$ made above
and, with this choice of $\zeta$, the morphism
$$
\roman{Id} \otimes \overline \tau \colon
\Lambda[v] \otimes_{\vartheta \tau^{\Bbb C P^{n-1}}} C_*(\Bbb C P^{n-1})
\longrightarrow
\Lambda[v] \otimes_{\vartheta}\Gamma_{n-1} [w]
$$
of twisted tensor products, necessarily a chain equivalence and compatible
with the requisite additional bundle structure, justifies
the claim that
$N=\Lambda[v] \otimes_{\vartheta}\Gamma_{n-1} [w]$
is an equivariant model for the
chains of the $(2n-1)$-sphere $S^{2n-1}$.

As a chain complex, $N$
decomposes canonically as $N = B \oplus \roman H_*(S^{2n-1})$
where $B$ is contractible. Indeed, $N$ has the additive generators
$ \gamma_k(w)$ and $v\otimes\gamma_k(w)$ where $0\leq k \leq n-1$
and
$$
d\gamma_k(w) = v \otimes \gamma_{k-1}(w),\quad dv \otimes
\gamma_k(w) =0,\quad 0\leq k \leq n-1.
$$
Now $B$ is the span of
$$
v, w, v \otimes w, \gamma_2(w),\dots,v \otimes \gamma_{n-2}(w),
\gamma_{n-1}(w)
$$
and $\roman H_*(S^{2n-1})$ is additively generated by
$1=\gamma_0(w)$ and $v \otimes \gamma_{n-1}(w)$.

Consider, then, the twisted tensor product
$$
\Gamma [u] \otimes_\vartheta N = \Gamma [u] \otimes_\vartheta
\Lambda[v] \otimes_{\vartheta}\Gamma_{n-1} [w].
$$
It calculates the $S^1$-equivariant cohomology of $S^{2n-1}$ and,
indeed, the homology of this twisted tensor product
plainly amounts to
$\Gamma_{n-1} [w]$, the homology of complex projective
$(n-1)$-space.

The obvious surjection from $N$ onto $\roman H_*(S^{2n-1})$
extends to a contraction
$$
\Nsddata {N}{\nabla}{\pi}{\roman H_*(S^{2n-1})}h
$$
in an obvious manner. Tensoring this contraction with $\Gamma
[u]$, we obtain the contraction
$$
\Nsddata {\Gamma [u] \otimes N}{\roman{Id}\otimes
\nabla}{\roman{Id}\otimes\pi}{\Gamma [u] \otimes\roman
H_*(S^{2n-1})}{\roman{Id}\otimes h}.
$$
Write the differential on $\Gamma [u] \otimes_{\vartheta} N$ as
$d+ \partial$ where the perturbation $\partial$ arises from the
twisting cochain $\vartheta$, that is, this perturbation equals
the operator $\vartheta \cap \,\cdot\,$. Application of the
perturbation lemma (Lemma 2.3) yields a perturbation $\Cal D$ of the zero
differential together with a contraction
$$\Nsddata {\Gamma [u] \otimes_{\vartheta} N}{\nabla_{\partial}}{g_{\partial}}
{\Gamma [u] \otimes \roman H_*(S^{2n-1}),\Cal D )} {h_{\partial}}.
$$
Evaluating the formula (2.3.1) under the present circumstances we
see that the perturbation $\Cal D$ is the twisted differential
associated with the twisting cochain $\tau$, that is,
$$(\Gamma [u] \otimes \roman
H_*(S^{2n-1}),\Cal D ) = \Gamma [u] \otimes_{\tau} \roman
H_*(S^{2n-1}) .
$$
Hence $\roman H_*(S^{2n-1})$ may be taken as a model for the
chains on $X$, and the sh-action of $\roman H_*(S^1)$ on $C_*X$
then amounts to the twisting cochain $\tau$ from $\Gamma [u]$ to
$\roman{End}(\roman H_*(S^{2n-1}))$.

\head Example 7.2 \endhead

Consider a homogeneous space $G/K$ of a
compact Lie group $G$ by a closed subgroup $K$. As noted earlier,
the cohomology of $G/K$ equals the $K$-equivariant cohomology of
$G$. Let $R$ be a ground ring such that $G$ and $K$ are both of
strictly exterior type over $R$. The familiar model
$$
\roman H_*(BK) \otimes_{\tau}\roman H_*G
$$
for the homology of $G/K$, cf. e.~g. \cite\husmosta, incorporates
the requisite sh-action of $\roman H_*K$ on $\roman H_*G$ via the
composite of the universal transgression for the universal
$G$-bundle with the  induced morphism from $\roman H_*(BK)$ to
$\roman H_*(BG)$. In particular, the induced $K$-action on the
total space $\roman T^*G$ of the cotangent bundle on $G$ is
hamiltonian, and the reduced space is the total space $\roman
T^*(G/K)$ of the cotangent bundle on $G/K$. Thus, the
(hamiltonian) $K$-action on $\roman T^*G$ cannot be equivariantly
formal. This shows once again that, in order for a hamiltonian
action of a compact Lie group to be equivariantly formal, the
compactness of the manifold is crucial.

An explicit (admittedly toy) example arising from $G= \roman{SU}(2)$ and 
$K=T$, a maximal torus which here comes down to the circle group $S^1$,
has been spelled out in the introduction and in Example 7.1 above.

\head Example 7.3 \endhead

 Let $K$ be a connected and simply connected
compact Lie group, and let $G=\roman{Map}^0(S^2,K)$, 
otherwise known as $\Omega^2 K$,
the group of
based smooth maps from $S^2$ to $K$, with pointwise
multiplication. Endowed with the induced differentiable structure
\cite\chenone, \cite\chentwo, $G$ is a Lie group. Let
$x_1,\dots,x_n$ be odd degree real cohomology classes of $K$  so
that the real cohomology $\roman H^*K$ is the exterior algebra on
these classes. Then there are real cohomology classes
$\xi_1,\dots,\xi_n$ of $G$ with $|\xi_j|=|x_j|-2$ for $1\leq j\leq
n$ so that, when $x_1,\dots,x_n$ are interpreted as cohomology
classes of $G$ in the obvious manner, the real cohomology $\roman
H^*G$ is the exterior algebra on
$x_1,\dots,x_n,\xi_1,\dots,\xi_n$;
see e.~g. \cite\atibottw\ or \cite\kan\ 
where such generators are constructed for 
a group of the kind $\roman{Map}^0(\Sigma,K)$
where $\Sigma$ is a closed surface of arbitrary genus 
rather than just the 2-sphere $S^2$ (for genus higher than zero
additional generators are necessary).
Our approach to Koszul duality
applies to the example $G=\roman{Map}^0(S^2,K)$.

\head Example 7.4 \endhead

Let $G$ be a Lie group of strictly exterior type over the
ground ring $R$ in such a way that
the duals of the exterior generators are universally transgressive. 
Then, with respect to the conjugation action 
of $G$ on itself, $G$ is {\sl equivariantly formal\/},
in fact, its equivariant cohomology is an induced module of the kind
$\roman H^*(G)\otimes \roman H^*(BG)$. 
Indeed, consider the path fibration $\Omega G \to PG \to G$; with
pointwise multiplication on $\Omega G$ and $PG$, the path
fibration is an extension of groups.
Moreover, with reference to the $G$-action on $\Omega G$ induced
by conjugation, as a group, the free loop group decomposes as the
semi-direct product $\Lambda G \cong \Omega G \rtimes G$. Now
$\Omega G \rtimes G$ acts on $PG$ where $\Omega G$ acts via
translation and $G$ by conjugation, and this action extends to a
free action on $PG \times EG$ through the projection to $G$ whence
the Borel construction $G\times_G EG$ is a classifying space for
$\Lambda G$. This classifying space, in turn, amounts to the free
loop space $\Lambda BG$
on the classifying space $BG$. 
Over the cohomology $\roman H^*(BG)$ of the classifying space $BG$,
the cohomology of $\Lambda BG$ is an induced module of the kind
$\roman H^*(G)\otimes \roman H^*(BG)$.

\beginsection 8. Split complexes and
generalized momentum mapping 

Let $C$ be a coaugmented
differential graded coalgebra and $A$ an augmented differential
graded algebra.

We will refer to a twisting cochain which is homotopic to the zero
twisting cochain as an {\it exact\/} twisting cochain; thus $\tau$
being exact means that there is a morphism $h\colon C \to A$ of
degree zero  such that
$$
Dh = \tau \cup h, \quad h \eta = \eta, \quad \varepsilon h =
\varepsilon .
\tag8.1
$$
Such a morphism  $h$ will then be referred to
as a {\it splitting\/} for $\tau$.

Let $\tau\colon C \to A$ be a twisting cochain, and let $N$ be a
differential graded right $A$-module. Consider the twisted
Hom-object $\roman{Hom}^{\tau}(C,N)$. We will say that this object
is {\it split\/} when the composite of the twisting cochain $\tau$
with the action $A \to \roman{End}(N)^{\roman{op}}$ is exact, that
is, when there is a morphism $h \colon C \to
\roman{End}(N)^{\roman{op}}$ of degree zero satisfying (8.1),
with $\roman{End}(N)^{\roman{op}}$ substituted for $A$.
We then refer to $h$ as a {\it splitting homotopy\/} for
$\roman{Hom}^{\tau}(C,N)$.
The notion of split object arises by abstraction
from the property of a $G$-space 
having the property that its equivariant cohomology
is an induced module over $\roman H^*(BG)$,
that is, makes precise the idea that the action
behaves like the trivial action, as far as equivariant cohomology
is concerned.

\proclaim{Proposition 8.2} Suppose that $C$ is cocomplete. Given  a
splitting homotopy $h$ for the twisted Hom-object
$\roman{Hom}^{\tau}(C,N)$, the assignment to a (homogeneous) $\phi \in
\roman{Hom}(C,N)$ of $\phi \cup h \in \roman{Hom}(C,N)$ yields an
isomorphism
$$
\roman{Hom}^{\tau}(C,N) @>>> \roman{Hom}(C,N)
$$
of twisted $\roman{Hom}$-objects, the target of the isomorphism
being untwisted.
\endproclaim

\demo{Proof} This is a special case of Lemma 1.2.1 above. \qed
\enddemo

Thus a twisted Hom-object which is split is isomorphic to
the corresponding untwisted object.

Suppose that, as a graded $R$-modules, $C$ is free
and that $N$ is the dual of a free graded $R$-module.
Consider a general twisted Hom-object
$\roman{Hom}^{\tau}(C,N)$.
In the standard manner, the filtration by $C$-degree yields a
spectral sequence
$(\roman E_r, d_r)\  (r \geq 0)$ which has
$$
(\roman E_0, d_0) \cong
(\roman{Hom}(C,N), d_N),
\
(\roman E_1, d_1) \cong
(\roman{Hom}(C,\roman H^* N), d_C),
\
\roman E_2 \cong
\roman H^*(C,\roman H^* N).
$$

\proclaim{Corollary 8.3} When the twisted Hom-object
$\roman{Hom}^{\tau}(C,N)$ is split the spectral
sequence $(\roman E_r, d_r)$ $(r \geq 0)$ collapses from $\roman
E_2$. 
\endproclaim

\demo{Proof} This is an immediate consequence of
Proposition 8.2. \qed
\enddemo

Let $N^\sharp$ be a differential graded left $A$-module 
whose dual $(N^\sharp)^*$ coincides with the differential graded 
right $A$-module $N$ and suppose that
$N^\sharp$ is non-negative and connected.
Suppose that there is a contraction
$\Nsddata{N^\sharp}{\nabla^\sharp}{\pi^\sharp}{\roman
H_*(N^\sharp)}{h^\sharp}$
of chain complexes.
Notice that, when the ground ring $R$ is a field, 
such a contraction necessarily exists. 

\proclaim{Theorem 8.4} Under these circumstances,
the twisted Hom-object
$\roman{Hom}^{\tau}(C,N)$ is split if and only if the spectral
sequence $(\roman E_r, d_r)$ $(r \geq 0)$ collapses from $\roman
E_2$.
\endproclaim

Thus, over a field, the collapse of the spectral sequence from $\roman E_2$
suffices to guarantee that the twisted Hom-object is split.

We begin with the preparations for the proof of Theorem 8.4.
Consider the differential graded algebra
$\roman{End}(N^\sharp)$ of homogeneous endomorphisms of $N^\sharp$, endowed
with the obvious differential graded algebra structure,
and let
$\Cal A=\roman{End}(N^\sharp)_{\geq 0}$
be the differential graded subalgebra of homogeneous 
endomorphisms of non-negative degree. The assignment to a degree zero 
endomorphism of $N^\sharp$ of its degree zero constituent
$N^\sharp_0 \to N^\sharp_0$ yields a multiplicative augmentation map
$\varepsilon \colon \Cal A \to R$ for $\Cal A$.
Likewise 
consider the graded algebra
$\roman{End}(\roman H_*(N^\sharp))$ of homogeneous
endomorphisms of $\roman H_*(N^\sharp)$, endowed with the obvious graded algebra structure,
and let
$\Cal B=\roman{End}(\roman H_*(N^\sharp))_{\geq 0}$
be the graded subalgebra of homogeneous 
endomorphisms of non-negative degree. Similarly as before, the assignment to a degree zero 
endomorphism of $\roman H_*(N^\sharp)$ of its degree zero constituent
$\roman H_0(N^\sharp) \to \roman H_0(N^\sharp)$ yields a multiplicative augmentation map
$\varepsilon \colon \Cal B \to R$ for $\Cal B$. Moreover, the assigment to
$\alpha \in \Cal A$ of $\Pi(\alpha)=\pi^\sharp \alpha \nabla^\sharp \in \Cal
B$
and to $\beta \in \Cal B$ of $\nabla(\beta)=\nabla^\sharp \beta \pi^\sharp \in
\Cal A$ yields chain maps
$$
\Pi\colon \Cal A @>>> \Cal B,\ \nabla \colon \Cal B @>>>  \Cal A
$$
such that $\Pi \nabla = \roman{Id}_{\Cal B}$ but neither $\Pi$ nor $\nabla$
are morphisms of differential graded algebras unless the differential of
$N^\sharp$ is zero.

Thus suppose that the differential  of
$N^\sharp$ is non-zero. While the projection $\Pi$ is not a 
morphism of differential graded algebras, 
using HPT, we will now extend $\Pi$ to an sh-morphism
from $\Cal A$ to $\Cal B$, that is, to a twisting cochain
$\tau_{\Pi}\colon \rbar \Cal A \to \Cal B$.

For greater clarity, we proceed somewhat more generally. Thus let $\CaA$ be
a general augmented differential graded algebra.
Recall that, by construction,
as a graded coalgebra, $\rbar \CaA$ is the graded tensor coalgebra
$$
\rbar \CaA = \roman T^{\roman c}[sI\CaA]=\oplus (sI\CaA)^{\otimes j}.
$$
The differential $d_{\rbar \CaA}$ on $\rbar \CaA$
takes the form
$$
d_{\rbar \CaA}=d_{\CaA}+ \partial.
$$
Here $d_{\CaA}$ refers to the differential induced by the differential
on $\CaA$, that is,  for $j \geq 1$, the operator 
$d_{\rbar \CaA}$ is defined on each homogeneous constituent 
$(sI\CaA)^{\otimes j}$ and then takes the form
$$
d_{\CaA}\colon  (sI\CaA)^{\otimes j}\longrightarrow (sI\CaA)^{\otimes j}
$$
of the tensor product differential on the $j$-fold tensor product 
$(sI\CaA)^{\otimes j}$ of the
suspension $sI\CaA$; thus, with the aid of the standard notation
$$
\left [\alpha_1|\alpha_2|\ldots|\alpha_j \right]
=(s \alpha_1) \otimes \ldots \otimes (s \alpha_j) \ (\alpha_{\nu} \in \CaA)
$$
and
$\overline \alpha= (-1)^{|\alpha|+1}=\alpha$  ($\alpha \in \CaA$),
the values of $d_{\CaA}$ are given by
$$
\align
d_{\CaA}[\alpha] &= -  \left[d (\alpha)\right]
\\
 d_{\CaA}\left[\alpha_1|\ldots|\alpha_j\right] &=
- \sum\left[
\overline \alpha_1|\ldots|\overline \alpha_{\nu -1}|
d \alpha_{\nu}|\alpha_{\nu +1}| \ldots|\alpha_{j}\right] .
\endalign
$$
Moreover,
$\partial$ is the coalgebra perturbation (relative to the
bar construction or coaugmentation filtration) determined by the
multiplication
in $\CaA$, that is, $\partial$ is the coderivation 
of  $\roman T^{\roman c}[sI\CaA]$ determined by 
the requirement that the diagram
$$
\CD
(sI\CaA) \otimes (sI\CaA) @>{\partial}>> (sI\CaA)
\\
@V{s^{-1}\otimes s^{-1}}VV
@A{s}AA
\\
(I\CaA) \otimes (I\CaA) @>{\mu}>> (I\CaA)
\endCD
$$
be commutative. Explicitly,  $\partial$ is zero
on $s I \CaA$ and, for $ j \geq 2$, the values of $\partial$ are given by
the identity
$$
 \partial\left[\alpha_1|\ldots|\alpha_j\right] =
\sum\left[
\overline \alpha_1|\ldots|\overline \alpha_{\nu -1}|
\overline \alpha_{\nu}\alpha_{\nu +1}|\alpha_{\nu +2}| 
\ldots|\alpha_{j}\right] .
$$

Let $\tau_1= \tau^{\rbar \CaA}\colon \rbar \CaA \to \CaA$, the universal 
bar construction twisting cochain for $\CaA$. Recall that 
$\tau_1$ is the desuspension on
$sI\CaA$ and zero elsewhere.
Henceforth we will occasionally write the multiplication in $\CaA$
as $\,\cdot\,\colon \CaA \otimes \CaA \to \CaA$, and
we denote by $D$ the unperturbed Hom-differential on
$\roman{Hom}(\rbar \CaA, \CaA)$, that is, the differential 
that relies only on the differential on $\CaA$; then the perturbed differential
on $\roman{Hom}(\rbar \CaA, \CaA)$
which takes into account the multiplication 
$\,\cdot\,$ on $\CaA$ takes the form
$D+\delta$ where $\delta(\varphi)=(-1)^{|\varphi|+1}\varphi \partial$ 
($\varphi \in \roman{Hom}(\rbar \CaA, \CaA)$).

By construction, $\tau_1$ is a cycle relative to the differential 
$d_{\CaA}$, that is,
$$
D\tau_1= d\tau_1 + \tau_1 d_{\CaA} = 0.
\tag8.5
$$
Furthermore, the twisting cochain property of $\tau_1$ 
is equivalent to the identity
$$
\tau_1 \partial =\tau_1\cup \tau_1
\colon \rbar \CaA \longrightarrow \CaA.
\tag8.6
$$
Let $p\in \CaA$ be an idempotent, necessarily of degree zero.
Using the idempotent $p$ of $\CaA$, we introduce a new composition
$\,\circ\,\colon \CaA \otimes \CaA \to \CaA$ in $\CaA$ by the assignment to $(\alpha,\beta)$ of
$$
\alpha \circ \beta = \alpha p \beta,\ \alpha,\beta \in \CaA.
$$
A little thought reveals that this composition is
associative, that is, given $\alpha,\beta,\gamma \in \CaA$,
$$
(\alpha \circ \beta)\circ \gamma=\alpha \circ (\beta\circ \gamma).
$$
Let $\odot$ denote the cup pairing in $\roman{Hom}(\rbar \CaA,\CaA)$
relative to the original coalgebra structure on 
$\rbar \CaA$ but relative to the new composition $\circ$ on
$\CaA$. Thus, given $\alpha,\beta \in \CaA$, the element
$\alpha \odot \beta$ of  $\roman{Hom}(\rbar \CaA,\CaA)$
is given as the composite
$$
\alpha \odot \beta\colon \rbar \CaA @>{\Delta}>> 
\rbar \CaA \otimes \rbar \CaA @>{\alpha \otimes \beta}>> 
\CaA \otimes \CaA
@>{\circ}>>
\CaA.
$$

\proclaim {Lemma 8.7} Given the idempotent $p$ of $\CaA$, let 
$\hhh\in \CaA$ be an element of degree {\rm 1\/} such that
$d \hhh = 1 - p$.
Let $\tau_1= \tau^{\rbar \CaA}\colon \rbar \CaA \to \CaA$ and,
for $j \geq 2$,
define
$\tau_j \colon (sI\CaA)^{\otimes j} \to \CaA$ by
$$
\tau_j\left[\alpha_1|\ldots|\alpha_j\right]=\alpha_1 \cdot \hhh
\cdot \alpha_2\cdot \ldots \cdot \alpha_{j-1}\cdot \hhh \cdot 
\alpha_j \quad \text{($j-1$\ copies of \ $\hhh$)}.
$$
Suppose that $\CaA$ is complete. Then
$$
\tau^{\hhh}= \tau_1+\tau_2 + \ldots \colon \rbar \CaA \to (\CaA,\,\circ\,)
$$
is a twisting cochain with respect  to the original coalgebra structure on 
$\rbar \CaA$ and the new composition $\,\circ\,$ on $\CaA$, that is, $\tau$ satisfies
the identity
$$
(D+\delta)(\tau) (=D (\tau) + \tau \partial) = \tau \odot \tau.
$$
\endproclaim

\demo{Proof}
Let $j \geq 2$. By construction
$$
\align
(D\tau_j)\left[\alpha_1|\ldots|\alpha_j\right]
&=d(\alpha_1 \cdot \hhh
\cdot \ldots \cdot \alpha_j) 
+ \tau_j  d_{\CaA}\left[\alpha_1|\ldots|\alpha_j\right]
\\
&=(d\alpha_1) \cdot \hhh \cdot \ldots \cdot \alpha_j +\ldots 
\\
&\quad -\overline \alpha_1 (1-p)\cdot \alpha_2\cdot \ldots \cdot \hhh \cdot 
\alpha_j 
\\
&\quad -
\overline \alpha_1 \cdot \hhh\cdot \overline \alpha_2\cdot (1-p) \ldots \cdot \hhh \cdot 
\alpha_j +\ldots
\\
&\quad -
\overline \alpha_1 \cdot \hhh\cdot \overline \alpha_2\cdot \hhh \ldots \cdot \hhh \cdot 
 \overline \alpha_{j-1}\cdot (1-p) \alpha_j
\\
&\quad -\tau_j 
 \sum\left[
\overline \alpha_1|\ldots|\overline \alpha_{\nu -1}|
d \alpha_{\nu}|\alpha_{\nu +1}| \ldots|\alpha_{n}\right]
\\
&=\sum_{\nu =1}^{j-1}
\overline \alpha_1 \cdot \hhh \cdot \ldots \cdot \hhh \cdot \overline \alpha_{\nu}\cdot (p-1)
\cdot \alpha_{\nu+1}\cdot  \hhh \cdot 
\ldots \cdot \hhh \cdot 
 \alpha_{j-1}\cdot \hhh \cdot \alpha_j
\\
&=\sum_{\nu =1}^{j-1}
\left(\overline \alpha_1 \cdot \hhh \cdot \ldots \cdot  \hhh \cdot \overline \alpha_{\nu}\right) \circ  
\left(\alpha_{\nu+1}\cdot   \hhh \cdot   
\ldots \cdot \hhh \cdot 
 \alpha_{j-1}\cdot \hhh \cdot \alpha_j\right )
\\
&\quad -
\sum_{\nu =1}^{j-1}
\left(\overline \alpha_1 \cdot \hhh \cdot \ldots \cdot  \hhh \cdot \overline \alpha_{\nu}\right) \cdot  
\left(\alpha_{\nu+1}\cdot  \hhh \cdot 
\ldots \cdot \hhh \cdot 
 \alpha_{j-1}\cdot \hhh \cdot \alpha_j\right ) .
\endalign
$$
Likewise, let $1 \leq \nu \leq j-1$. Since
$$
\align
(\tau_{\nu}\cup \tau_{j-\nu})
\left[\alpha_1|\ldots|\alpha_j\right]
&=
\left(\tau_{\nu} \left[\overline \alpha_1|\ldots|\overline \alpha_{\nu}\right]\right)
\cdot
\left(\tau_{j-\nu}\left[ \alpha_{n+1}|\ldots| \alpha_j\right]\right)
\\
(\tau_{\nu}\odot \tau_{j-\nu})
\left[\alpha_1|\ldots|\alpha_j\right]
&=
\left(\tau_{\nu} \left[\overline \alpha_1|\ldots|\overline \alpha_{\nu}\right]\right)
\circ
\left(\tau_{j-\nu}\left[\alpha_{n+1}|\ldots| \alpha_j\right]\right)
\endalign
$$
we conclude
$$
D\tau_j=\sum_{\nu =1}^{j-1}\tau_{\nu}\odot \tau_{j-\nu}
-\sum_{\nu =1}^{j-1}\tau_{\nu}\cup \tau_{j-\nu} .
\tag8.7.1
$$
The same kind of calculation yields
$$
\tau_{j-1}\partial =
\sum_{\nu =1}^{j-1}\tau_{\nu}\cup \tau_{j-\nu}
\tag8.7.2
$$
Indeed,
$$
\align 
\tau_{j-1}\partial \left[\alpha_1|\ldots|\alpha_j\right]&=
\sum
\tau_{j-1}\left[
\overline \alpha_1|\ldots|\overline \alpha_{\nu -1}|
\overline \alpha_{\nu}\alpha_{\nu +1}|\alpha_{\nu +2}| 
\ldots|\alpha_{j}\right]
\\
&=
\sum
\left(\overline \alpha_1\cdot \hhh \cdot\ldots\cdot \hhh \cdot\overline \alpha_{\nu -1}\cdot \hhh \cdot
\overline \alpha_{\nu}\right)\left(\alpha_{\nu +1}\cdot \hhh \cdot\alpha_{\nu +2}\cdot \hhh \cdot 
\ldots\cdot \hhh \cdot\alpha_{j}\right)
\\
&=
\sum_{\nu =1}^{j-1}\tau_{\nu}\cup \tau_{j-\nu}\left[\alpha_1|\ldots|\alpha_j\right] .
\endalign
$$
Combining the identities (8.7.1) and (8.7.2)
with the twisting cochain property of $\tau_1$,
cf. (8.5) and (8.6), we conclude
$$
D(\tau_1+\tau_2 +\ldots)+(\tau_1+\tau_2+\ldots)\partial
= (\tau_1+\tau_2 +\ldots) \odot (\tau_1+\tau_2 +\ldots)
$$
as asserted. \qed
\enddemo

We now apply Lemma 8.7 to the differential graded algebra $\Cal A$ together with the idempotent
$p=\nabla^{\sharp}\pi^{\sharp}\in \Cal A$ and the homotopy $h^{\sharp}\in \Cal A$.
By construction,
$$
\tau_{\Pi}= \Pi\tau^{h^{\sharp}} \colon \rbar \Cal A \to \Cal B
$$
is a twisting cochain which extends the projection $\Pi \colon \Cal A \to \Cal B$. 

The original $A$-action on $N^\sharp$ amounts to a
morphism $A\to \Cal A$ of differential graded algebras, and the composite
$$
C @>{\tau}>> A @>>> \Cal A
\tag8.8
$$
is a twisting cochain. Likewise the composite 
$$
C @>>> \rbar \Cal A @>{\tau_{\rcob}}>> \rcob \rbar \Cal A
$$
of the adjoint of (8.8) with the universal cobar construction twisting cochain
$\tau_{\rcob}$ is a twisting cochain, and we can substitute $\rcob \rbar \Cal
A$ for $A$
and simply write $\tau \colon C \to \rcob \rbar \Cal
A$,
that is, we consider $N^\sharp$ as a differential graded 
left $\rcob \rbar \Cal A$-module
and $N$ as a differential graded 
right $\rcob \rbar \Cal A$-module
through the adjunction map
$\rcob \rbar \Cal A \longrightarrow \Cal A$.

The adjoint
$$
\overline{\tau_{\Pi}}\colon \rcob \rbar \Cal A \longrightarrow \Cal B
$$
of the twisting cochain $\tau_{\Pi}$
is a surjective morphism of differential graded algebras
which extends to a contraction
$$
\Nsddata{\Cal B}{\overline{\tau_{\Pi}}}{\nabla^\flat}{\rcob \rbar \Cal A}{h^\flat}.
$$

\proclaim{Lemma 8.9}
Suppose that the following data are given:

\noindent
--- augmented differential graded algebras ${\CC}$ and ${\BB}$;

\noindent
--- a contraction
$\Nsddata {\CC}{\nabla}{\pi}{\BB}h$ of chain complexes, $\pi$ being a
morphism of augmented differential graded algebras;

\noindent
--- a coaugmented differential graded coalgebra $C$;

\noindent
--- twisting cochains $t_1,t_2 \colon {C} \to {\CC}$;

\noindent
--- a homotopy $h^{\BB}\colon {C} \to {\BB}$ of twisting cochains $h^{\BB}\colon
\pi t_1 \simeq \pi t_2$, so that
$$
D(h^{\BB})=(\pi t_1)\cup h^{\BB} - h^{\BB}\cup(\pi t_2). 
\tag{8.9.1}
$$
Suppose that $C$ is cocomplete. Then the recursive rule
$$
h^{\CC}=\nabla h^{\BB} - h(t_1 \cup h^{\CC} - h^{\CC}\cup t_2) \tag{8.9.2}
$$
yields a homotopy $h^{\CC}\colon {C} \to {\CC}$ of twisting cochains
$h^{\CC}\colon t_1\simeq t_2$ such that $\pi h^{\CC} = h^{\BB}$.
\endproclaim

The rule (8.9.2) being recursive means that
$$
h^{\CC}=  \varepsilon \eta + h_1 + h_2 + \ldots
$$
where 
$h_1=\nabla h^{\BB} - h(t_1 -t_2)$, $h_2= - h(t_1 \cup h_1 - h_1\cup t_2)$, 
etc.

\demo{Proof} The identity $\pi h^{\CC} = h^{\BB}$ is obvious and,
since $t_1$ and $t_2$ are ordinary twisting cochains, the morphism
$t_1 \cup h^{\CC} - h^{\CC}\cup t_2$ is (easily seen to be) a cycle.
Furthermore, since $\pi$ is compatible with the algebra
structures,
$$
\align
\nabla \pi(t_1 \cup h^{\CC} - h^{\CC}\cup t_2) &=
\nabla ((\pi t_1) \cup (\pi h^{\CC}) - (\pi h^{\CC})\cup (\pi t_2))  
\\&=  \nabla ((\pi t_1) \cup h^{\BB} - h^{\BB}\cup (\pi t_2)).
\endalign
$$
Consequently
$$
\align
Dh^{\CC} &= \nabla (D(h^{\BB})) + Dh (t_1 \cup h^{\CC} - h^{\CC}\cup t_2)
\\&=
\nabla ((\pi t_1)\cup h^{\BB} - h^{\BB}\cup(\pi t_2))
 +
(t_1 \cup h^{\CC} - h^{\CC}\cup t_2)- \nabla
\pi(t_1 \cup h^{\CC} - h^{\CC}\cup t_2)
\\
&= t_1 \cup h^{\CC} - h^{\CC}\cup t_2
\endalign
$$
as asserted. \qed
\enddemo

\demo{Proof of Theorem {\rm 8.4} } 
Corollary 8.3 says that the condition is necessary.
To justify that the condition is also sufficient we note first that
the hypothesis that the spectral sequence collapses implies that
the twisting cochain that arises as the composite
$$
C @>{\tau}>> \rcob \rbar \Cal A @>{\overline{\tau_{\Pi}}}>> \Cal B
$$
of the twisting cochain $\tau\colon C \to \rcob \rbar \Cal A$
with the adjoint $\overline{\tau_{\Pi}}$ is the zero twisting cochain.
Lemma 8.9, applied with 
$\BB=\Cal B$, $\CC= \rcob \rbar \Cal A$,
$t_1= \tau$, $t_2$ the zero twisting cochain,
and $h^{\BB}$ the zero homotopy,
then yields a splitting $h \colon C \to \rcob \rbar \Cal A$
for $\tau$ so that
$Dh = \tau \cup h$,  $h \eta = \eta$ and $\varepsilon h =\varepsilon$,
cf. (8.1). \qed
\enddemo

\noindent{\smc Remark 8.10.\/} Theorem 8.4 does not contradict Example 5.2
in \cite\franpupp. To adjust the present notation to that example,
let $C=\roman H^*(BS^1)$, let $\Cal C Y$ denote the cone on the space $Y$,
let $X=S^2 \cup_{\phi}(\Cal C\Bbb R \Bbb P^2 \times S^1)$
be the space explored in  \cite\franpupp,
$\phi$ being the map
$\Bbb R \Bbb P^2 \times S^1 @>{\roman{pr}}>> \Bbb R \Bbb P^2 \to \Bbb R \Bbb
P^2\big/\Bbb R \Bbb P^1\cong S^2$, and let $N=C^*(X)$. 
(We use the font $\Cal C$ to distinguish the cone from our notation $C$ for
a differential graded coalgebra.) The space
$X$ has
$\roman H^0(X)$ and $\roman H^2(X)$ infinite cyclic,
$\roman H^4(X)$ cyclic of order 2, and its other cohomology groups are zero.
Plainly there is no way to contract $C_*(X)$ onto $\roman H_*(X)$, so there is
no contradiction. In fact, 
under such circumstances, the existence of the contraction
implies that the equivariant cohomology is an extended module
once the spectral sequence collapses. 
On the other hand, the Example 5.2
in \cite\franpupp\ shows that the collapse of the spectral sequence
from $\roman E_2$ alone does not guarantee that the equivariant 
cohomology is an extended module.
I am indebted to the referee for having
insisted that this point be clarified.

\smallskip\noindent
{\smc 8.11. Splitting homotopy and momentum mapping\/.}
We will now explain in which sense a splitting homotopy
generalizes a momentum mapping: As before, let $G$ be a
topological group of strictly exterior type, let $C = \rbar C_*G$,
$A= C_*G$, $\tau = \tau_G\colon \rbar C_*G \to C_*G$, let $X$ be a
left $G$-space, and let $N=C^*(X)$, viewed as a right
$(C_*G)$-module as above. In view of Lemma 8.3, Proposition 9.3 in
\cite\gorkomac\ shows that our notion of split twisted Hom-object
is consistent with the notion of split complex used in
\cite\gorkomac. Consider a 2-cocycle $\zeta$ of $X$. In the double
complex $(\roman E_2,d_2)$, the class $[\zeta] \in \roman H^2(X)$
lies on the fiber line, in fact, in $\roman E_2^{0,2}$.
Since $G$ is of striclty exterior type, 
$d_2[\zeta] \in \roman E_2^{2,1}$ is the obstruction to the
existence of an equivariantly closed extension of $\zeta$, that
is, the obstruction to the existence of a pre-image in
$\roman H_G^2(X)$ relative to the restriction mapping
$\roman H_G^2(X)\to \roman H^2(X)$.
This obstruction comes down 
to
the obstruction to the existence of a
$C^0(X)$-valued
2-cochain $\varphi$ on $\rbar C_*G$ (where $C^0(X)$ refers to the
singular zero cochains on $X$) such that
$$
d \circ \varphi = \zeta \cup \tau \colon \rbar C_*G @>>> C^1(X)
\tag8.11.1
$$
where $d$ is the differential on $C^*(X)$. Such a cochain
$\varphi$ is related to $\zeta$ in the same manner as a momentum
mapping to an equivariant  closed 2-form,
a momentum mapping
for the 2-form 
(not necessarily non-degenerate and {\it not\/} equivariantly closed)
being an equivariantly closed extension. 

To explain what this means,
suppose that $G$
is a compact Lie group, let $X$ be a $G$-manifold and, instead of
$C^*(X)$, take the de Rham complex $\Cal A(X)$.
The Cartan model
for the $G$-equivariant de Rham cohomology of $X$ has the form
$$
\roman{Hom}^{\tau^{\Sigm'}}(\Sigm'[s^2\fra g], \Cal A(X))^G;
\tag8.11.2
$$
here $\Sigm'[s^2\fra g]$ is the graded symmetric coalgebra on the
double suspension $s^2\fra g$, 
the algebra
$\Lambda [s\fra g]$ is the graded
exterior algebra on the suspension $s\fra g$, 
the twisting cochain
$\tau^{\Sigm'}
\colon \Sigm'[s^2\fra g] \to \Lambda [s\fra g]$ is the universal
twisting cochain, and $\Lambda [s\fra g]$ acts on $\Cal A(X)$ via
contraction. We now substitute a closed $G$-invariant 2-form
$\sigma$ for $\zeta$, a $G$-equivariant $\Cal A^0(X)$-valued
2-cochain $\Phi$ on $\Sigm'[s^2\fra g]$ for the $C^0(X)$-valued
2-cochain $\varphi$ on $\rbar C_*G$, that is, essentially a
$G$-equivariant linear map $\Phi$ from $\fra g$ to $\Cal A^0(X) =
C^{\infty}(X)$ and, likewise, we substitute the twisting cochain
$\tau^{\Sigm'}$ for $\tau$. The identity (8.11.1) then translates
to the identity
$$
d \circ \Phi = \sigma \cup \tau^{\Sigm'} \colon \Sigm'[s^2\fra g]
@>>> \Cal A^1(X) \tag8.11.3
$$
where $d$ is the de Rham differential on $X$. With the
notation $\mu \colon X \to \fra g^*$
for adjoint of $\Phi$, we may rewrite
(8.11.3) as
$$
\sigma(\xi_X, \cdot) = \xi \circ d \mu,\quad \xi \in \fra g,
\tag8.11.4
$$
where $\xi_X$ refers to the fundamental vector field
on $X$ coming from $\xi \in \fra g$,
and this is exactly the momentum mapping property.

This recovers the familiar fact that a momentum mapping
for a closed $G$-equivariant 2-form $\sigma$
(not necessarily non-degenerate and {\it not\/} equivariantly closed)
is an equivariantly closed extension of $\sigma$,
that is, $d_2[\sigma]$ is the obstruction to the existence
of a momentum mapping for $\sigma$.

These substitutions can be given an entirely rigorous meaning;
to this end, one has to develop the formalism
with $\Sigm'[s^2\fra g]$ instead of $\rbar C_*G$.
In particular, given a splitting homotopy
$h \colon \Sigm'[s^2\fra g] \to \roman{End}(\Cal A(X))$
for $\tau^{\Sigm'}$, the corresponding adjoint thereof includes
a morphism of the kind
$$
h^{\flat} \colon \fra g \otimes \Cal A^2(X) @>>> \Cal A^0(X)
$$
such that, given a closed 2-form $\sigma \in \Cal A^2(X)$,
the association 
$$
\fra g \ni \xi
\longmapsto
h^{\flat}(\xi\otimes \sigma) \in \Cal A^0(X)
\tag8.11.5
$$
yields a \lq\lq comomentum\rq\rq\ for $\sigma$, i.~e.
the adjoint thereof is an ordinary momentum mapping for $\sigma$.

Thus,
given a general topological group of strictly exterior type
and a $G$-action on a space $X$,
{\sl a splitting homotopy $h$ for $\tau$ generalizes
the concept of momentum mapping in a very strong sense\/}
since it provides, in particular,
{\sl a single object which yields, via the association\/} (8.11.5),
{\sl a momentum mapping for every closed 2-cocycle
and, furthermore, the correct replacement thereof
for cocycles of arbitrary degree\/}.

\bigskip

\widestnumber\key{999}
\centerline{References}
\smallskip\noindent

\ref \no  \atibottw
\by M. F. Atiyah and R. Bott
\paper The Yang-Mills equations over Riemann surfaces
\jour Phil. Trans. R. Soc. London  A
\vol 308
\yr 1982
\pages  523--615
\endref

\ref \no \borelboo \by A. Borel \book Seminar on transformations
groups \bookinfo Annals of Math. Studies vol. 46 \publ Princeton
University Press \publaddr Princeton, New Jersey \yr 1960
\endref

\ref \no \bottone
\by R. Bott
\paper On the Chern-Weil homomorphism and the continuous cohomology of
Lie groups
\jour Advances
\vol 11
\yr 1973
\pages  289--303
\endref

\ref \no \bottsega \by R. Bott and G. Segal \paper The cohomology
of the vector fields on a manifold \jour Topology \vol 16 \yr 1977
\pages  285--298
\endref

\ref \no \botshust \by R. Bott, H. Shulman, and J. Stasheff \paper
On the de Rham theory of certain classifying spaces \jour Advances
\vol 20 \yr 1976 \pages 43--56
\endref

\ref \no \ebrown
\by E. Brown
\paper Twisted tensor products.~I
\jour Ann. of Math.
\vol 69
\yr 1959
\pages  223--246
\endref

\ref \no \brownez
\by R. Brown
\paper The twisted Eilenberg--Zilber theorem
\jour Celebrazioni Archimedee del Secolo XX, Simposio di topologia
\yr 1964
\pages 33-37
\endref

\ref \no \cartanon
\by H. Cartan
\paper Notions d'alg\`ebre diff\'erentielle; applications aux groupes
de Lie et aux vari\'et\'es o\`u op\`ere un groupe de Lie
\jour Coll. Topologie Alg\'ebrique
\paperinfo Bruxelles
\yr 1950
\pages  15--28
\endref

\ref \no \cartantw
\by H. Cartan
\paper La transgression dans un groupe de Lie et dans un espace
fibr\'e principal
\jour Coll. Topologie Alg\'ebrique
\paperinfo Bruxelles
\yr 1950
\pages  57--72
\endref

\ref \no \chenone
\by K.T. Chen
\paper Iterated path integrals
\jour Bull. Amer. Math. Soc.
\vol 83
\yr 1977
\pages  831--879
\endref

\ref \no \chentwo
\bysame 
\paper Degeneracy indices and Chern classes
\jour  Adv. in Math.
\vol 45
\yr 1982
\pages 73--91
\endref

\ref \no \chenfou
\bysame 
\paper Extension of $C^{\infty}$
Function Algebra by Integrals and Malcev Completion of $\pi_1$
\jour Advances in Mathematics
\vol 23
\yr 1977
\pages 181--210
\endref

\ref \no \doldone
\by A. Dold
\paper Zur Homotopietheorie der Kettenkomplexe
\jour Math. Ann.
\vol 140
\yr 1960
\pages  278--298
\endref

\ref \no \doldpupp
\by A. Dold und D. Puppe
\paper Homologie nicht-additiver Funktoren. Anwendungen
\jour Annales de l'Institut Fourier
\vol 11
\yr 1961
\pages  201--313
\endref

\ref \no \drachman
\by B. Drachman
\paper  A note on principal constructions
\jour Duke Math. J.
\vol 39
\yr 1972
\pages 701-710
\endref

\ref \no \eilmactw
\by S. Eilenberg and S. Mac Lane
\paper On the groups ${\roman H(\pi,n)}$. I.
\jour Ann. of Math.
\vol 58
\yr 1953
\pages  55--106
\moreref
\paper II. Methods of computation
\jour Ann. of Math.
\vol 60
\yr 1954
\pages  49--139
\endref

\ref \no \eilmotwo
\by S. Eilenberg and J. C. Moore
\paper Limits and spectral sequences
\jour Topology
\vol 1
\yr 1961
\pages  1--23
\endref

\ref \no \eilmothr
\bysame 
\paper Foundations of relative homological algebra
\jour Memoirs AMS
\vol 55
\yr 1965
\publ Amer. Math. Soc.
\publaddr Providence, Rhode Island
\endref

\ref \no \eilenmoo
\bysame
\paper Homology and fibrations.~I.
Coalgebras, cotensor product and its derived functors
\jour Comm. Math. Helv.
\vol 40
\yr 1966
\pages  199--236
\endref

\ref \no \frankone \by T. Frankel \paper Fixed points and torsion
on K\"ahler manifolds \jour  Ann. Math. II. Ser. \vol 70 \yr 1959
\pages 1--8
\endref

\ref \no \franzone \by M. Franz \paper Koszul duality and
equivariant cohomology for tori \jour Int. Math. Res. Notices \vol
42 \yr 2003 \pages 2255--2303 \finalinfo{\tt math.AT/0301083}
\endref

\ref \no \franztwo
\bysame  
\paper Koszul duality and equivariant cohomology
\paperinfo{\tt math.AT/0307115}
\endref

\ref \no \franpupp \by M. Franz and V. Puppe
\paper Exact cohomology sequences
with integral coefficients
for torus actions
\jour Transformation groups \vol 12
\yr 2007 \pages 65--76 \finalinfo{\tt math.AT/0505607}
\endref

\ref \no \gorkomac
\by M. Goresky, R. Kottwitz, and R. Mac Pherson
\paper Equivariant cohomology,
Koszul duality and the localization theorem
\jour Invent. Math.
\vol 131
\yr 1998
\pages 25--83
\endref

\ref \no \gugenhtw
\by V.K.A.M. Gugenheim
\paper On the chain complex of a fibration
\jour Illinois J. of Mathematics
\vol 16
\yr 1972
\pages 398--414
\endref

\ref \no \gugenhon
\bysame 
\paper On a perturbation theory for the homology of the loop space
\jour J. of Pure and Applied Algebra
\vol 25
\yr 1982
\pages 197--205
\endref

\ref \no \gulstatw
\by V.K.A.M. Gugenheim, L. Lambe, and J.~D. Stasheff
\paper Perturbation theory in differential homological algebra. II.
\jour Illinois J. of Math.
\vol 35
\yr 1991
\pages 357--373
\endref

\ref \no \gugenmay
\by V.K.A.M. Gugenheim and J.~P. May
\paper On the theory and applications of differential
torsion products
\jour Memoirs of the Amer. Math. Soc.
\vol 142
\yr 1974
\endref

\ref \no \gugenmun
\by V.K.A.M. Gugenheim and H.~J. Munkholm
\paper On the extended functoriality of Tor and Cotor
\jour J. of Pure and Applied Algebra
\vol 4
\yr 1974
\pages  9--29
\endref

\ref \no \habili
\by J. Huebschmann
\paper Perturbation theory and small models for the chains of
certain induced fibre spaces
\paperinfo Habilitationsschrift, Universit\"at Heidelberg, 1984
\finalinfo {\bf Zbl} 576.55012
\endref

\ref \no \perturba
\bysame  
\paper Perturbation theory and free resolutions for nilpotent
groups of class 2
\jour J. of Algebra
\yr 1989
\vol 126
\pages 348--399
\endref

\ref \no \cohomolo
\bysame 
\paper Cohomology of nilpotent groups of class 2
\jour J. of Algebra
\yr 1989
\vol 126
\pages 400--450
\endref

\ref \no \modpcoho
\bysame 
\paper The mod $p$ cohomology rings of metacyclic groups
\jour J. of Pure and Applied Algebra
\vol 60
\yr 1989
\pages 53--105
\endref

\ref \no \intecoho
\bysame 
\paper Cohomology of metacyclic groups
\jour Trans. Amer. Math. Soc.
\vol 328
\yr 1991
\pages 1-72
\endref

\ref \no \kan
\bysame 
\paper Extended moduli spaces, the
Kan construction, and lattice gauge theory
\jour Topology
\vol 38
\yr 1999
\pages 555--596
\finalinfo{\tt dg-ga/9505005, dg-ga/9506006}
\endref

\ref \no \berikas
\bysame 
\paper Berikashvili's functor $\Cal D$ and the deformation equation
\paperinfo Festschrift in honor of N. Berikashvili's 70-th birthday;
{\tt math.AT/9906032}
\jour Proceedings of the A. Razmadze Mathematical Institute
\vol 119
\yr 1999
\pages 59--72
\endref

\ref \no \minimult
\bysame 
\paper Minimal free multi models for chain algebras
\paperinfo in: Chogoshvili Memorial
\jour Georgian Math. J.
\vol 11
\yr 2004
\pages 733--752
\finalinfo {\tt math.AT/0405172}
\endref

\ref \no \pertlie
\bysame 
\paper The Lie algebra perturbation lemma
\paperinfo in: Festschrift in honor of M. Gerstenhaber's 
80-th and J. Stasheff's
70-th birthday (to appear)
\finalinfo{\tt arxiv 0708:3977}
\endref

\ref \no \pertlitw
\bysame 
\paper The sh-Lie algebra perturbation lemma
\jour Forum math. (to appear)
\finalinfo{\tt arxiv 0710:2070}
\endref

\ref \no \origins
\bysame 
\paper Origins and breadth of the theory of higher homotopies
\paperinfo in: Festschrift in honor of M. Gerstenhaber's 
80-th and J. Stasheff's
70-th birthday (to appear)
\finalinfo{\tt arxiv 0710:2645}
\endref

\ref \no \duatwo 
\bysame 
\paper Relative homological
algebra, homological perturbations, and equivariant de Rham
cohomology \paperinfo{\tt math.DG/0401161}
\endref

\ref \no \equivg 
\bysame 
\paper Equivariant cohomology over groupoids and Lie-Rinehart algebras
\paperinfo preprint 2009 
\endref

\ref \no \tornike
\bysame 
\paper On the construction of $A_{\infty}$-algebras
\paperinfo {\tt arxive:0809.4791}
\jour Georgian Math. Journal (to appear)
\endref

\ref \no \huebkade
\by J. Huebschmann and T. Kadeishvili
\paper Small models for chain algebras
\jour Math. Z.
\vol 207
\yr 1991
\pages 245--280
\endref

\ref \no \huebstas
\by J. Huebschmann and J. D. Stasheff
\paper Formal solution of the master equation via HPT and
deformation theory
\paperinfo {\tt math.AG/9906036}
\jour Forum mathematicum
\vol 14
\yr 2002
\pages 847--868
\endref

\ref \no\husmosta
\by D. Husemoller, J. C. Moore, and J. D. Stasheff
\paper Differential homological algebra and homogeneous spaces
\jour J. of Pure and Applied Algebra
\vol 5
\yr 1974
\pages  113--185
\endref

\ref \key \kirwaboo \by F. Kirwan \book Cohomology of quotients in
symplectic and algebraic geometry \publ Princeton University Press
\publaddr Princeton, New Jersey \yr 1984
\endref

\ref \no \maclaboo
\by S. Mac Lane
\book Homology
\bookinfo Die Grundlehren der mathematischen Wissenschaften
 No. 114
\publ Springer Verlag
\publaddr Berlin $\cdot$ G\"ottingen $\cdot$
Heidelberg \yr 1963
\endref

\ref \no \mooretwo
\by J. C. Moore
\paper Cartan's constructions
\paperinfo Colloque analyse et topologie, en l'honneur de Henri
Cartan
\jour Ast\'erisque
\vol 32--33
\yr 1976
\pages  173--221
\endref

\ref \no \munkholm
\by H. J. Munkholm
\paper The Eilenberg--Moore spectral sequence and strongly homotopy
multiplicative maps
\jour J. of Pure and Applied Algebra
\vol 5
\yr 1974
\pages  1--50
\endref

\ref \no \munkhotw
\bysame
\paper Shm maps of differential algebras. I.
A characterization up to homotopy
\jour J. of Pure and Applied Algebra
\vol 9
\yr 1976
\pages  39--46
\moreref
\paper II. Applications to spaces with polynomial cohomology
\jour J. of Pure and Applied Algebra
\vol 9
\yr 1976
\pages  47--63
\endref

\ref \no \munkhthr
\bysame  
\paper
DGA algebras as a Quillen model category; relations to shm maps
\jour J. of Pure and Applied Algebra
\vol 13
\yr 1978
\pages  221--232
\endref

\ref \no \rothstee
\by M. Rothenberg and N. Steenrod
\paper The cohomology of classifying spaces of H-spaces
\jour Bull. Amer. Math. Soc.
\vol 71
\yr 1965
\pages  872--875
\endref

\ref \no \sourithr \by J. M. Souriau \paper Groupes
diff\'erentiels 
\inbook Diff. geom. methods in math.
Physics, Proc. of a conf., Aix en Provence and Salamanca, 1979,
Lecture Notes in Mathematics, No. 836 \publ Springer
Verlag
\publaddr Berlin $\cdot$ Heidelberg $\cdot$ New York $\cdot$
Tokyo \yr 1980 \pages 91--128
\endref

\ref \no \stashalp
\by J.~D. Stasheff and S. Halperin
\paper Differential algebra in its own rite
\jour Proc. Adv. Study Alg. Top. August 10--23, 1970, Aarhus, Denmark
\pages 567--577
\endref

\ref \no \steeepst
\by N.~E. Steenrod and D.~B.~A. Epstein
\book Cohomology Operations
\bookinfo Annals of Mathematics Studies, vol. 50
\publ Princeton University Press
\publaddr Princeton, N. J. 08540
\yr 1962
\endref

\enddocument